%% file: main.tex
\documentclass[12pt]{article}

\include{prematter}

\title{\Large Error Estimators for the Small-Biot Lumped Approximation\\for the Conduction Dunking Problem}

\author[1]{\bf{Kento Kaneko}}
\author[2]{\bf{Claude Le Bris}}
\author[1]{\bf{Anthony T Patera}}

\affil[1]{\small{Massachusetts Institute of Technology, Cambridge, MA 02139, USA. \newline\texttt{kaneko@mit.edu}, \texttt{patera@mit.edu}}}

\affil[2]{\'Ecole des Ponts and INRIA, Champs-sur-Marne, 77455 Marne La Vall\'ee, France. \newline\texttt{claude.le-bris@enpc.fr}}

\date{}

\begin{document}

\maketitle
\vspace{-0.9cm}
\centerline{\bf \small Acknowledgments}
\vspace{0.15cm}

{\small
\small We thank Prof. Masayuki Yano for his many contributions to the adaptive refinement finite element method software suite used extensively in this work.

This work is supported by ONR Grant N00014-21-1-2382, Grant Monitor Dr Reza Malek-Madani. The research of the second author is partially supported by ONR and EOARD.
}

\input{abs}      % Abstract

\input{intro}    % Section 1
\input{eqn}      % Section 2
\input{esb}      % Section 3
\input{lump}     % Section 4
\input{phi}      % Section 5
\input{comp}     % Section 5.4
\input{err}      % Section 6
\input{uavg1}    % Section 6.1
\input{uavg2}    % Section 6.2
\input{udelta}   % Section 6.3

\appendix
\input{appa}

\input{appb}

\input{appc}

\input{appd}

{
  \bibliographystyle{ieeetr}
  \addcontentsline{toc}{part}{Bibliography}
  \bibliography{sbp}
}

\end{document}

%% file: prematter.tex
\usepackage[left=1.in,right=1.in,top=1.in,bottom=1.in]{geometry}
\usepackage{tikz}
\usepackage{authblk}
\usetikzlibrary{arrows.meta, positioning, calc}
\usepackage{verbatim,setspace}
\usepackage{lmodern}
\usepackage{pdfcomment}
\usepackage{showlabels}

\usetikzlibrary{circuits.ee.IEC, positioning}

\usepackage[rflt]{floatflt}
\usepackage{float}
\usepackage{subcaption}

\usepackage{graphicx}
\graphicspath{{./Figures/}}

\usepackage{amsthm,amsmath,amssymb,amsfonts,amscd,amsbsy,thmtools}
\usepackage{epstopdf}
\usepackage{calrsfs}
\DeclareMathAlphabet{\pazocal}{OMS}{zplm}{m}{n}

\usepackage{textcomp}

\usepackage{hyperref}
\hypersetup{colorlinks=true,linkcolor=red,citecolor=red}

\usepackage{lmodern}

\usepackage{caption}
\def\LiLt{{L^{\infty}({\calT},L^2(\Omega))}}
\def\LtX{{L^2({\calT},X)}}

\DeclareMathAlphabet{\pazocal}{OMS}{zplm}{m}{n}

\def\u{{u}}

\def\OmMeas{{|\Omega |}}
\def\dOmMeas{{|\partial\Omega |}}

\def\Ti{{T_{\text{i}}}}
\def\tf{{t_{\text{final}}}}
\def\Tf{{T_{\text{final}}}}

\def\R0{{\text{RHE}^0}}

\def\calB{\mathcal{B}}

\def\calD{\mathcal{D}}

\def\calL{\mathcal{L}}

\def\calR{\mathcal{R}}

\def\calT{{\mathcal{T}}}

\def\II{\mathbb{I}}
\def\MM{\mathbb{M}}

\def\FF{\mathbb{F}}

\def\NN{\mathbb{N}}
\def\PP{\mathbb{P}}
\def\QQ{\mathbb{Q}}
\def\RR{\mathbb{R}}

\def\TT{\mathbb{T}}

\def\pcalO{\pazocal{O}}

\def\avg{{\text{avg}}}
\def\paavg{{\partial{\text{avg}}}}

\def\tf{{t_{\text{final}}}}

\def\Xint#1{\mathchoice
{\XXint\displaystyle\textstyle{#1}}%
{\XXint\textstyle\scriptstyle{#1}}%
{\XXint\scriptstyle\scriptscriptstyle{#1}}%
{\XXint\scriptscriptstyle\scriptscriptstyle{#1}}%
\!\int}
\def\XXint#1#2#3{{\setbox0=\hbox{$#1{#2#3}{\int}$ }
\vcenter{\hbox{$#2#3$ }}\kern-.6\wd0}}

\def\dashint{\Xint-}

\newcommand{\dv}[1]{{\color{blue} {#1}}}

\def\bnabla{\boldsymbol{\nabla}}

\usepackage{mdframed}

\definecolor{light-gray}{gray}{0.90}
\definecolor{light-blue}{RGB}{173,216,230}

\declaretheoremstyle[
spaceabove=2pt,
spacebelow=9pt,
notefont=\normalfont,
notebraces={}{},
headfont=\itshape,
bodyfont=\normalfont,
headformat=\NAME~{\em {\NUMBER}}{\normalfont }\NOTE,
mdframed={backgroundcolor=light-gray}
]{mystyle0}

\declaretheoremstyle[
spaceabove=2pt,
spacebelow=9pt,
notefont=\normalfont,
notebraces={}{},
headfont=\itshape,
bodyfont=\normalfont,
headformat=\NAME~{\em {\NUMBER}}{\normalfont }\NOTE,
mdframed={backgroundcolor=white}
]{mystyle1}

\declaretheoremstyle[
spaceabove=2pt,
spacebelow=9pt,
notefont=\normalfont,
notebraces={}{},
headfont=\itshape,
bodyfont=\normalfont,
headformat=\NAME~{\em {\NUMBER}}{\normalfont }\NOTE,
mdframed={backgroundcolor=light-blue}
]{mystyle3}

\declaretheoremstyle[
spaceabove=2pt,
spacebelow=9pt,
notefont=\normalfont,
notebraces={}{},
headfont=\itshape,
bodyfont=\normalfont,
headformat=\NAME~{\em {\NUMBER}}{\normalfont }\NOTE,
mdframed={backgroundcolor=light-gray}
]{mystyle}

\declaretheorem[style=mystyle,name={Proposition},numberwithin=section]{proposition}

\declaretheoremstyle[
spaceabove=2pt,
spacebelow=9pt,
notefont=\normalfont,
notebraces={}{},
headfont=\itshape,
bodyfont=\normalfont,
headformat=\NAME~{\em {\NUMBER}}{\normalfont }\NOTE,
mdframed={backgroundcolor=white}
]{mystyle}

\declaretheorem[style=mystyle,name={Lemma},numberlike={proposition}]{lemma}

\declaretheorem[style=mystyle0,name={Definition},numberlike={proposition}]{definition}

\def\OmMeas{|\Omega|}

\def\psiz{\psi^0}
\def\lambdaz{\lambda^0}
\def\psipz{\psi'^{\,0}}
\def\lambdapz{\lambda'^{\,0}}
\def\psip{\psi'}
\def\lambdap{\lambda'}
\def\psipp{\psi''}
\def\psippz{\psi''^{\,0}}
\def\lambdapp{\lambda''}
\def\lambdappz{\lambda''^{\,0}}
\def\lambdapppz{\lambda'''^{\,0}}

\def\Rzp{\RR_{+}}
\def\Rp{\RR_+^\star}

\newcommand{\uavgMID}{u_\avg^{\text{MID}}}
\newcommand{\uavgMIDp}{u_\avg^{\text{MID}'}}

\newcommand{\uavgLB}{u_\avg^{\text{LB}}}
\newcommand{\uavgLBp}{u_\avg^{\text{LB}'}}
\newcommand{\uavgUB}{u_\avg^{\text{UB}}}
\newcommand{\Bsup}{\calB_{\text{sup}}}
\newcommand{\Binf}{\calB_{\text{inf}}}

\newcommand{\Bbar}{\overline{\calB}}

\newcommand{\Bdunk}{\text{Bi}_\text{dunk}}
\newcommand{\Bdunkp}{\text{Bi}'_\text{dunk}}

\newcommand{\uavg}{u_\avg}
\newcommand{\tiuavg}{\tiu_\avg}
\newcommand{\uavgL}[1]{\tilde{u}_\avg^{#1}}
\newcommand{\UavgL}[1]{\widetilde{U}_\avg^{#1}}
\newcommand{\uavgLP}{\tilde{u}_\avg^{2\text{P}}}
\newcommand{\UavgLP}{\widetilde{U}_\avg^{2\text{P}}}
\newcommand{\uDeltaL}{\tilde{u}_\Delta^{2\text{P}}}
\newcommand{\UDeltaL}{\widetilde{U}_\Delta^{2\text{P}}}

\newcommand{\lambdaL}[1]{\tilde{\lambda}^{#1}}
\newcommand{\lambdaLP}{\tilde{\lambda}^{2\text{P}}}

\newcommand{\EavgL}[1]{E_\avg^{#1}}
\newcommand{\EavgLP}{E_\avg^{2\text{P}}}
\newcommand{\EavgLa}[1]{E_\avg^{{#1}\,\text{asymp}}}
\newcommand{\EavgLr}[1]{E_\avg^{{#1}\,\text{rel}}}
\newcommand{\EavgLUB}[1]{E_\avg^{{#1}\,\text{UB}}}
\newcommand{\EavgLra}[1]{E_\avg^{{#1}\,\text{rel asymp}}}
\newcommand{\EavgLrUB}[1]{E_\avg^{{#1}\,\text{rel UB}}}

\newcommand{\EDeltaL}{E_\Delta^{2}}
\newcommand{\EDeltaLr}{E_\Delta^{{2}\,\text{rel}}}

\newcommand{\Lcond}{\calL_\text{cond}}
\newcommand{\tiu}{u^{\star}}

\def\wm{{\tiu_{-}}}

\def\wm{{\tiu_{-}}}

\def\LiLt{{L^{\infty}({\calT},L^2(\Omega))}}

\def\LtX{{L^2({\calT},X)}}

\def\dhaus{{d_{\text{Hausdorff}}}}

\newcommand{\Feat}[1]{\text{Feature}_{#1}}

\definecolor{kk}{RGB}{191, 0, 255}

\DeclareMathOperator*{\argmax}{arg\,max}
\DeclareMathOperator*{\argsup}{arg\,sup}
\DeclareMathOperator*{\essinf}{ess\,inf}
\DeclareMathOperator*{\esssup}{ess\,sup}

\newcommand{\norm}[1]{\left\lVert#1\right\rVert}
\newcommand{\sci}[2]{#1\;\!\cdot\;\!10^{#2}}

%% file: abs.tex
\begin{abstract}
We consider the dunking problem: a solid body at uniform temperature $T_{\text i}$ is placed in a environment characterized by farfield temperature $T_\infty$ and spatially uniform time-independent heat transfer coefficient. We permit heterogeneous material composition: spatially dependent density, specific heat, and thermal conductivity. Mathematically, the problem is described by a heat equation with Robin boundary conditions. The crucial parameter is the Biot number --- a nondimensional heat transfer (Robin) coefficient; we consider the limit of small Biot number.

We introduce first-order and second-order asymptotic approximations (in Biot number) for several quantities of interest, notably the spatial domain average temperature as a function of time; the first-order approximation is simply the standard engineering `lumped' model. We then provide asymptotic error estimates for the first-order and second-order approximations for small Biot number, and also, for the first-order approximation, alternative non-asymptotic bounds valid for all Biot number. Companion numerical solutions of the heat equation confirm the effectiveness of the error estimates for small Biot number.

The second-order approximation and the first-order and second-order error estimates depend on several functional outputs associated to an elliptic partial differential equation; the latter is derived from Biot-sensitivity analysis of the heat equation eigenproblem in the limit of small Biot number. Most important is $\phi$, the only functional output required for the first-order error estimates; $\phi$ admits a simple physical interpretation in terms of conduction length scale. We investigate the domain and property dependence of $\phi$: most notably, we characterize spatial domains for which the standard lumped-model error criterion --- Biot number (based on volume-to-area length scale) small --- is deficient.
\end{abstract}


%% file: intro.tex
\section{Introduction}\label{sec:intro}

The dunking problem  is ubiquitous in heat transfer engineering education and professional practice: a body characterized by spatial domain $\dv{\Omega}$ and boundary $\dv{\partial\Omega}$, at initial uniform temperature $\dv{T_\text{i}}$, is immersed at time $\dv{t=0}$ in an environment (fluid and enclosure) at farfield temperature $\dv{T_\infty}$ \cite{LandL, Cengel, IandD}. The temperature distribution of the body, $\dv{T}$, evolves over the time interval $\dv{0 < t\le \tf}$. Note dimensional quantities shall be coded in blue, and nondimensional quantities (to be introduced later) in black. In this introduction, for simplicity of exposition, we largely consider homogeneous body composition, and hence uniform thermophysical properties --- in particular volumetric specific heat, $\dv{\rho c}$, and scalar thermal conductivity, $\dv{k}$; however, in subsequent development, we shall treat the important case of material heterogeneity (and hence non-uniform thermophysical properties), particularly in the form of several isotropic materials. The thermal environment is characterized, most simply, by an average heat transfer coefficient, $\dv {h}$. Finally, we associate to the body domain two length scales: a chosen extrinsic length scale, $\dv{\ell}$; an intrinsic length scale, $\dv{\calL \equiv |\Omega|/|\partial\Omega|}$, where $\dv{|\Omega|}$ and $\dv{|\partial\Omega|}$ denote the volume and surface area of the body, respectively. We note that the use of the notation ``$\equiv$'' in this work signifies ``is defined as'' (equivalent to the ``:='' notation).

The solution $\dv{T(\cdot,t)}$ to the dunking problem satisfies, for $\dv{t\in(0,\tf]}\,$,
\begin{align}
\dv{\rho c\, \partial_t T = \bnabla \cdot(k \bnabla T) \text{ in }\Omega\,,}\label{eq:dimivpstrong_1}
\end{align}
subject to boundary condition
\begin{align}
\dv{k\, \partial_n T + h (T-T_\infty) = 0 \text{ on }\partial\Omega}\,,\label{eq:dimivpstrong_2}
\end{align}
and initial condition
\begin{align}
\dv{T(\cdot,0) = T_\text{i} \text{ in } \Omega}\,. \label{eq:dimivpstrong_3}
\end{align}
Here, $\partial_t$ and $\partial_n$ refer to differentiation with respect to time and domain boundary outward normal direction, respectively.

We nondimensionalize the spatial coordinate by $\dv{\ell}$ and time by $\dv{\rho c \ell^2/k}$. The nondimensional temperature is defined in terms of the dimensional temperature $\dv{T}$ as $u \equiv \dv{(T-T_\infty)/(\Ti - T_\infty)}$. The mathematical formulation of the nondimensional dunking problem then corresponds to a heat equation characterized by the spatial domain $\Omega$ and a single parameter, the Robin coefficient. The Robin coefficient is given by an extrinsic Biot number, $B \equiv \dv{h\ell/k}$; we also define an intrinsic Biot number, $\Bdunk  \equiv \dv{h \calL/k} = B \calL$, where $\calL = \dv{\calL/\ell}$ (note $\Bdunk$ is typically indicated in textbooks simply by $\text{Bi}$ \cite{LandL,Cengel,IandD}). As we shall see, the general case of heterogeneous properties is considerably more complicated: the nondimensional volumetric specific heat and thermal conductivity will now be functions of the spatial coordinate.

In this work, we study the limit of small Biot number. The small-Biot limit arises quite often in practice in particular for ``everyday'' artifacts at modest temperatures subject to natural convection and radiation. In contrast, larger systems at very elevated temperatures subject to brisk forced convection or change-of-phase typically yield larger Biot numbers. An example of a small-Biot application is the process of annealing by natural convection --- for instance, to relieve residual stresses.

We consider two quantities of interest: the (spatial) domain average temperature as a function of time, $\dv{T_\avg (t)}$ — also directly related to the heat loss (or gain) of the body from (or to) the environment; the (spatial) domain-boundary average temperature as a function of time, $\dv{T_\paavg(t)}$. The nondimensional versions of these quantities of interest are given respectively by $u_\avg(t) \equiv \dv{(T_\avg(}t\dv{\rho c \ell^2/k)-T_\infty)/(T_{\text{i}} - T_\infty)}$ and $u_\paavg(t) \equiv \dv{(T_\paavg(}t\dv{\rho c \ell^2/k) - T_\infty)/(T_{\text{i}} - T_\infty)}$; it shall also prove convenient to define $u_\Delta(t) \equiv [({u_\avg} - {u_\paavg})/{u_\avg}](t)$, which measure the relative temperature variation within the body.

The first-order\footnote{Note in general ``order'' refers to the convergence rate in $B$ (or, equivalently, $\Bdunk$).} `classic' lumped approximation to $u_\avg$ for the small-$\Bdunk$ dunking problem is a simple exponential in time, $\uavgL{1}(t) \equiv \exp(-B \gamma t)$ \cite{LandL,Cengel,IandD}, where $\gamma \equiv 1/\calL$.  However, and despite the simplicity and utility of this small-$\Bdunk$ lumped approximation, there is relatively little analysis of the associated approximation error. For example, most textbooks \cite{LandL,Cengel,IandD} indicate only that $\Bdunk$ must be small; in some cases, a threshold might be provided, for example $\Bdunk \leq 0.1$ — but typically qualified by ``usually''. A notable exception is the important work of Gockenbach and Schmidtke \cite{ODD}, in which a rigorous asymptotic error estimate is provided for the particular case of a sphere with homogeneous thermophysical properties and uniform heat transfer coefficient: $|\uavgL{1}(t) - u_\avg(t)| \leq (3/5)\,\Bdunk/e$ for all time $t$. The analysis of \cite{ODD} is based on an eigenvalue expansion in Biot.

We enumerate here the specific contributions of this work:
\begin{enumerate}
\item We consider small-$\Bdunk$ approximation and error estimation based on asymptotic expansion of the first eigenfunction (for the eigenproblem associated with the heat equation), $\psi_1(\cdot;B) = \psi_1(\cdot;0) + B\frac{d\psi_1}{d B}(\cdot;0) + \pcalO(B^2)$, where $\frac{d\psi_1}{dB}(\cdot;0)$ is found following standard sensitivity-derivative approaches \cite{Joseph,Kato} as the solution to an elliptic partial differential equation (PDE). The key quantity in our analysis is $\phi\equiv \int_\Omega \kappa \nabla \frac{d\psi_1}{dB}(\cdot;0) \cdot\nabla \frac{d\psi_1}{dB}(\cdot;0)$, where $\kappa$ is the nondimensionalized thermal conductivity based on the infimum of $\dv{k}$ over $\dv{\Omega}$.  For canonical domains and selected triangles, we provide closed-form solutions for $\phi$; we also provide an efficient computational framework for evaluation of $\phi$ for general domains. 

\item We introduce a new second-order lumped approximation to the solution of the thermal dunking problem, $\uavgLP$, which is dependent on $\phi$. The new approximation has improved accuracy over the `classic' approximation, $\uavgL{1}$, for \textit{general geometry} and \textit{general thermophysical property distribution}. We also introduce an approximation for the spatial temperature variation, $\uDeltaL$.

\item We present an error analysis for all three approximations: $\uavgL{1}$, $\uavgLP$, and $\uDeltaL$.  The error analysis can be performed at marginal computational cost over the cost of computing $\phi$. In particular, the error estimate for the first-order approximation, $\phi B/(\gamma e)$, is directly computable once we are equipped with $\phi$. For the first-order approximation, we can provide two additional results: we demonstrate that $\uavgL{1}$ is in fact a lower bound of $u_\avg$ for all admissible $t$ and $B$; we provide a \textit{non-asymptotic} error bound for $\uavgL{1}$ which is valid for all $B$. 
\item We provide an explicit engineering interpretation of $\phi$ in terms of an improved engineering conduction length scale associated with a body thermal resistance. The latter generalizes in a rigorous way the usual definition of thermal resistance --- between two surfaces --- to the notion of thermal resistance between body effective center and body boundary.  Within this interpretation framework, we motivate $\uavgL{1}$, $\uavgLP$, and $\uDeltaL$.

    \item We provide an extensive discussion of the dependence of $\phi$ on spatial domain --- in fact, only body \textit{shape} matters --- and body composition. In some cases we can provide useful approximations and bounds for $\phi$ which require relatively little information about (say) the spatial distribution of the thermophysical properties, $\dv{\rho c}$ and $\dv{k}$. We define a distance between spatial domains which predicts well the stability of $\phi$ with respect to geometric perturbation. More generally, we identify two geometry features in terms of which we can characterize spatial domains for which $\phi$ is large --- and the textbook lumped error criterion potentially misleading.

    \item We provide a numerical study of $\phi$ with emphasis on geometry and thermophysical properties. The numerical results are based on an adaptive finite element framework \cite{Yano} with associated \textit{a posteriori} error estimators; the latter ensure that the finite element error does not compromise any of our conclusions.  We confirm (i) the effectiveness of our error estimators, and (ii) the relevance of our distance and geometry features in understanding the dependence of $\phi$ on spatial domain.
\end{enumerate} 
We provide more extensive references in context.

There is one aspect in which our analysis is incomplete: the treatment of the heat transfer coefficient, $\dv{h}$. We consider here almost exclusively the case in which the heat transfer coefficient is independent of time and furthermore uniform in space. For the case of natural convection and radiation, the heat-transfer coefficient is implicitly linearized about the initial temperature of the body. The more realistic case of temporally-dependent and spatially-dependent heat transfer coefficient is best addressed in the context of a more general discussion which takes as point of departure the ``truth” conjugate heat transfer formulation; this approach will be reported in a subsequent publication. However, we provide even here, in Appendix \ref{sec:htc_sketch} and Appendix \ref{sec:realflows}, a simple first result.

We now provide a roadmap for the manuscript. 
\begin{itemize}
\item[] In Section \ref{sec:equations} we present the $\Bdunk$ initial-value-problem formulation as well as the associated $\Bdunk$ eigenproblem; we also develop the connection between the former and the latter through a separation of variables expansion. 
\item[] In Section \ref{sec:smallb} we then develop the small-Biot asymptotic analysis of the eigenproblem: as $B \rightarrow 0$, (i) $\lambda_1$ (the first eigenvalue) $\sim B \gamma - \phi B^2$, and (ii) the spatial structure is well-approximated by $\psi_1$ (the first eigenmode). We also provide in Section \ref{sec:smallb} the sensitivity formulation for $\phi$. In Section \ref{sec:phi} we then analyze the sensitivity formulation for $\phi$ as regards general properties and particular illustrative solutions.  Note many of the results of Section \ref{sec:asymp} are particular instantiations of more general eigenproblem perturbation techniques \cite{Joseph,Kato}.

\item[] In Section \ref{sec:lumped} we develop, through the resistance formulation, the lumped approximations $\uavgL{1}(t)$, $\uavgLP(t)$, and $\uDeltaL$. In Section \ref{sec:qoierr} we then analyze the error in these lumped approximations and provide, for all approximations, asymptotic error estimates (as $B \rightarrow 0$) as well as strict error bounds, valid for all $B$, for $\uavgL{1}$.
\end{itemize}
Note we intercalate the discussion of $\phi$ and the lumped approximation for pedagogical purposes.

The main part of the manuscript, intended for a more general heat transfer audience, emphasizes the presentation and interpretation of the error estimates; we include the statement of the theoretical results but without proof. The appendices to the manuscript include all necessary mathematical details as well as sketches of the proofs of the main propositions. We note that many of the mathematical proofs are indeed elementary and many of the results are well known, however we provide them here for completeness and for the convenience of the reader.

%% file: eqn.tex
\section{Governing Equations}\label{sec:equations}

\subsection{Initial Value Problem}

A passive solid body at uniform temperature $\dv{T_{\text{i}}}$ is abruptly placed --- dunked --- at time $\dv{t = 0}$ in an  environment, fluid and enclosure, at initial and farfield temperature $\dv{T_\infty}$. Note that blue variables are dimensional, and black variables are nondimensional. Most of the analysis presented in this work can be extended to more general contexts, for example uniform heat generation in the body, $\dv{q^{\textsc{v}}(t)}$, or time-dependent farfield temperature, $\dv{T_\infty(t)}$. Note that, in general, units are SI, and temperature is in degrees Celsius (or Kelvin).

We first introduce the spatial domain of the body we consider, $\dv{\Omega} \subset \RR^d$, and associated boundary $\dv{\partial\Omega}$; we consider $d\in\{1,2,3\}$.  We shall require $\dv{\Omega}$ Lipschitz. A point in $\dv{\Omega}$ shall be denoted $\dv{x \equiv (x_1,\ldots,x_d)}$. We  choose for our length scale some characteristic dimension of $\dv{\Omega}$, $\dv{\ell}$, for example the diameter or the InRadius. We also introduce an intrinsic length scale $\dv{\mathcal{L}\equiv |\Omega|/|\partial\Omega|}$, and an intrinsic associated inverse length scale, $\dv{\gamma \equiv |\partial\Omega|/|\Omega|}$; here $|\,\cdot |$ denotes measure, and hence, in dimension $d = 3$, $\dv{\gamma}$ is the dimensional surface area to volume ratio. We shall denote time by $\dv{t}$; we shall restrict attention to the temporal interval of interest $\dv{t \in [0,\tf]}$. 

We next introduce two nondimensional quantities related to thermophysical properties:
\begin{align}
\sigma = \dfrac{\dv{\rho c}}{\dashint_\Omega \dv{(\rho c)}} \, ,
\end{align}
\begin{align}
\kappa = \dv{\dfrac{k}{k_{\inf}}} \, ,
\end{align}
where $\dv{\rho c}$ and $\dv{k}$ are, respectively, the volumetric specific heat and thermal conductivity, both strictly positive and scalar-valued (generalization to anisotropic heat-transfer with symmetric positive definite second-order $k$ tensor is not a subject of this work). Note $\dv{\rho}$ is the density and $\dv{c}$ is the mass specific heat; the product $\dv{\rho c}$ is the volumetric specific heat (which is bounded away from zero). Here $\dashint$ refers to the average of the integrand over the domain of integration, and $\dv{k_{\inf}}$ (more properly $\dv{k_{\essinf}}$) is the infimum of $\dv{k}$ over $\dv{\Omega}$. It follows from our definitions that $\sigma > 0$ and $\dashint_\Omega \sigma = 1$, (and $\kappa > 0$) and $\inf_{\Omega}\kappa = 1$.  Note also that $\sigma = 1$ corresponds to spatially uniform volumetric specific heat, and $\kappa = 1$ corresponds to spatially uniform thermal conductivity.

We may then define the extrinsic Biot number as
\begin{align}
B \equiv \dv{\frac{h \ell}{k_{\inf}}}\,, \label{eq:Bdef}
\end{align}
and the intrinsic Biot number as
\begin{align}
\mathrm{Bi}_\mathrm{dunk} = \dv{\dfrac{\dv{h \mathcal{L}}}{\dv{k_{\inf}}}}\; ; \label{eq:Bidef}
\end{align}
and hence $\mathrm{Bi}_\mathrm{dunk} = B(\dv{\mathcal{L}/\ell}) = B/(\ell\dv{\gamma})$. Note that
\begin{align}
    \gamma\equiv \dv{\ell\gamma}\label{eq:gammadef}
\end{align}
shall appear prominently in the exposition.
 As indicated in the introduction, we consider in the present work, the case in which the heat transfer coefficient $\dv{h}$  is taken as uniform in space and constant in time. We further assume that the heat transfer coefficient is non-negative, and hence also $B$ is non-negative. 

We now define the nondimensional spatial coordinate as $x \equiv \dv{x/\ell}$ and the nondimensional temporal coordinate as $t \equiv \dv{t/t_{\text{diff}}}$, where
\begin{align}
\dv{t_{\text{diff}} \equiv} \dfrac{\dv{\ell^2} \,\dashint_{\Omega}\dv{(\rho c)}}{\dv{k_{\inf}}} \, .\label{eq:def_tdiff}
\end{align}
We further introduce the nondimensional temperature
\begin{align}
u \equiv \dv{\dfrac{\dv{T} - \dv{T_\infty}}{\dv{T_{\text{i}}} - \dv{T_\infty}}} \, . \label{eq:def_und}
\end{align}
Our nondimensionalization reduces to the standard textbook form in the case in which the thermophysical properties are uniform; notably, $t$ is the Fourier number, $\text{Fo}_{\dv{\ell}}$.

The nondimensional temperature, $u(x,t;B;\sigma,\kappa)$, then satisfies, for any $B \in \Rzp\equiv\{x\in\RR\, |\, x \ge 0\}$,
\begin{align}
\sigma \partial_t u = \bnabla \cdot(\kappa \bnabla u) \text{ in }\Omega\,, \; 0 < t \le \tf\,, \label{eq:ivpstrong_1}
\end{align}
subject to boundary condition
\begin{align}
\kappa \partial_n u + B u = 0 \text{ on }\partial\Omega\,, \; 0 < t \le \tf \,, \label{eq:ivpstrong_2}
\end{align}
and initial condition
\begin{align}
u = 1 \text{ in } \Omega \text{ at } t = 0 \; . \label{eq:ivpstrong_3}
\end{align}
Here $\partial_t$ and $\partial_n$ refer to the partial derivatives with respect to time and the outward normal direction, respectively.  We note that, strictly speaking, the strong formulation \eqref{eq:ivpstrong_1} -- \eqref{eq:ivpstrong_3} requires $\sigma$ and $\kappa$ sufficiently smooth. However, the corresponding weak formulation, provided in Appendix \ref{App:ivp}, admits discontinuous data, and hence our analysis  --- based on the weak form --- shall directly apply to the practically important case of material heterogeneity.

We shall in general suppress the $\sigma$ and $\kappa$ parameters unless we are specifically focused on the dependence on these variables. We shall write $u(\cdot,t;B)$ or simply $u(t;B)$ to indicate the spatial nondimensional temperature function for any time $t$ and Biot number $B$ (and, implicitly, prescribed $\sigma$ and $\kappa$); we shall write $u(x,t;B)$ to refer to the evaluation of the spatial nondimensional temperature function at spatial coordinate $x \in \Omega$.  

\subsection{Eigenproblem}

We can associate to our initial value problem, in standard fashion, the corresponding eigenproblem: for given $B \in \Rzp$, find $(\psi_j,\lambda_j)_{j = 1,2,\ldots}$ solution to
\begin{align}
-\bnabla \cdot(\kappa \bnabla\psi_j) = \lambda_j \sigma \psi_j \text{ in } \Omega \,, \label{eq:evpstrong_1}
\end{align}
subject to boundary condition
\begin{align}
\kappa \partial_n \psi_j+ B \psi_j = 0 \text{ on }\partial\Omega\,,  \label{eq:evpstrong_2}
\end{align}
and normalization
\begin{align}
\int_\Omega \sigma \psi^2_j  = 1 \; . \label{eq:evpstrong_3}
\end{align}
The corresponding weak form is provided in Appendix \ref{App:evp}. We shall write $\psi_j(\cdot;B)$ or simply $\psi_j(B)$ to indicate the spatial eigenfunction for any given Biot number $B$ (and implicitly, prescribed $\sigma$ and $\kappa$); we shall write $\psi_j(x;B)$ to refer to the evaluation of the spatial eigenfunction function at spatial coordinate $x \in \Omega$.

It is a standard result that the eigenvalues are non-negative for $B \ge 0$. We shall enumerate the eigenvalues in increasing magnitude, hence $0 \le \lambda_1(B) <  \lambda_2(B) \le \lambda_3(B) \le \ldots$. It is also standard to note that the eigenfunctions are orthonormal:
\begin{align} \int_\Omega \sigma \psi_i(B) \, \psi_j(B) = \delta_{ij}\; , \; i,j = 1,2,\dots \; ; \label{eq:efortho}
\end{align}
here $\delta_{ij}$ is the Kronecker-delta symbol.

Considering $\lambda_1$, it can be shown that $\lambda_1$ is of multiplicity unity, and that (i) $\psi_1$ has constant sign, (ii) $\psi_1$ is (up to a change of sign) strictly positive, and (iii) $\psi_1$ is (up to a change of sign) unique. As regards (ii) and (iii), see for example \cite{Arendt}.  We shall choose $\psi_1(B)$ positive, which, with \eqref{eq:evpstrong_3}, uniquely specifies $\psi_1(B)$.
We shall henceforth adopt the convention that $\psi(B)$ (respectively, $\lambda(B)$) without subscript shall refer to the first eigenfunction, $\psi_1(B)$ (respectively, the first eigenvalue, $\lambda_1(B)$), and that furthermore $\psiz$ (respectively, $\lambdaz$) shall refer to $\psi(B = 0)$ (respectively, $\lambda(B = 0))$.

\subsection{Separation of Variables Representation}

Much of our analysis shall rely on the classical separation of variables representation of our solution. We begin with a lemma which shall serve subsequently in several contexts:
\begin{lemma}[Representation of Unity]\label{lem:rou}
We may write, for any $B \in \Rzp$,
\begin{align}
\lim_{N\rightarrow \infty}\, \left|\left| \; 1 - \sum_{k = 1}^N \OmMeas M(\psi_k(B)) \psi_k(\cdot;B)\;  \right|\right|_{L^2(\Omega)}\! \rightarrow 0 \; , \label{eq:Mpsi_one0}
\end{align} 
or, in abbreviated form (to be interpreted as \eqref{eq:Mpsi_one0})
\begin{align}
1 = \sum_{k = 1}^\infty \OmMeas M(\psi_k(B)) \psi_k(\cdot;B) \; ; \label{eq:Mpsi_one}
\end{align}
here $\| w \|_{L^2(\Omega)} \equiv (\int_\Omega w^2)^{1/2}$.
Furthermore,
\begin{align}
\lim_{N \rightarrow \infty}\,  \OmMeas \sum_{k = 1}^N (M(\psi_k(B)))^2 = 1 \; , \label{eq:Mpsisquare0}
\end{align}
or, in abbreviated form (to be interpreted as \eqref{eq:Mpsisquare0})
\begin{align}
1 = \OmMeas \sum_{k = 1}^\infty (M(\psi_k(B)))^2  \; . \label{eq:Mpsisquare}
\end{align}
Here $|\Omega|$ is the measure of $\Omega$ and $M:X \rightarrow \RR$ is given by
\begin{align}
M(v;\sigma) \equiv \dashint_{\Omega} \sigma v \; . \label{eq:Mdef}
\end{align}
We note that for $\sigma = 1$ (uniform volumetric specific heat) $M$ evaluates to the domain mean of $v$. Recall also that $\psi(B)$ refers to the function $\psi(\cdot;B) \in X$, whereas $\psi(x;B)$ refers to the evaluation of $\psi(\cdot;B)$ at spatial coordinate $x \in \Omega$.

Hence $X\equiv H^1(\Omega)$, where $H^1(\Omega)$ is imbued with norm $||\cdot||_{H^1}\equiv||\cdot||_X=\left(\int_\Omega \cdot\right)^{1/2}$
\end{lemma}
The sketch of the proof of Lemma \ref{lem:rou} is provided in Appendix \ref{App:sketch_rou}.

We can then write our solution as a standard separation of variables sum:

\begin{proposition}[Separation of Variables: Field]\label{prop:sov} For any $B \in \Rzp$,
\begin{align}
u^N(\cdot,t;B) \equiv \sum_{k = 1}^N \OmMeas M(\psi_k(B)) \psi_k(\cdot;B) \exp(-\lambda_k(B) t) \; . \label{eq:sov_sum0}
\end{align}
Then, for all $t \in (0,\tf]$,
\begin{align}
\lim_{N \rightarrow \infty} \| \; u^N(\cdot,t;B) - u(\cdot,t;B) \; \|_X \rightarrow 0 \;, \label{eq:sov_lim}
\end{align}
or in abbreviated form (to be interpreted as \eqref{eq:sov_lim}) 
\begin{align}
u(\cdot,t;B) \equiv \sum_{k = 1}^\infty \OmMeas M(\psi_k(B)) \psi_k(\cdot;B) \exp(-\lambda_k(B) t) \; .\label{eq:sov_sum}
\end{align}
Recall that the $X$ norm is the $H^1(\Omega)$ norm.
\end{proposition}
The sketch of the proof of Proposition \ref{prop:sov} is provided in Appendix \ref{App:sketch_sov}. We note here only that $u(\cdot,t;B)$ of \eqref{eq:sov_sum} for $t = 0$ reduces to $\sum_{k = 1}^\infty\OmMeas M(\psi_k(B))\psi_k(\cdot;B) = 1$ from Lemma \ref{lem:rou}.

\subsection{Quantities of Interest}

Quantity of Interest (QoI) shall refer to a functional of the temperature $u(\cdot,t;B)$. We shall consider in the present work three QoI: the domain mean, denoted $u_\avg$; the boundary mean, denoted $u_\paavg$; and the relative domain mean-boundary mean difference, $u_\Delta$. In our analysis and results, we focus on $u_\avg$ and $u_\Delta$; $u_\paavg$ is a necessary intermediate definition. The domain mean QoI is defined as
\begin{align}
u_\avg(t;B) \equiv M(u(\cdot,t;B)) = \dashint_{\Omega} \sigma u(\cdot,t;B) \; ;\label{eq:uavg}
\end{align}
the boundary mean QoI is defined as
\begin{align}
u_\paavg(t;B) \equiv H(u(\cdot,t;B)) = \dashint_{\partial\Omega} u(\cdot,t;B) \; ;\label{eq:upaavg}
\end{align}
the relative domain mean-boundary mean difference QoI is defined as
\begin{align}
u_\Delta(t;B) \equiv \dfrac{u_\avg(t;B) - u_\paavg(t;B)}{u_\avg(t;B)} \; .\label{eq:udelta}
\end{align}
Here, $M$ is defined in \eqref{eq:Mdef}, and $H: X \rightarrow \RR$ is given by
\begin{align}
H(v) \equiv \dashint_{\partial\Omega} v \; ; \label{eq:Hdef}
\end{align}
note $H$ evaluates to the boundary mean of $v$. These QoIs are defined for almost all times $t$ in our interval $[0,\tf]$. Note it follows from the maximum principle \cite{Evans} that $u_\avg$ and $u_\paavg$ take on values in the range $[0,1]$.

We can now express our QoI in terms of our separation of variables representation:
\begin{proposition}[Separation of Variables: QoI]\label{prop:sov_QoI}
We may write
\begin{align}
u_\avg(t;B) = \sum_{k = 1}^\infty \OmMeas (M(\psi_k(B)))^2 \exp(-\lambda_k(B) t) \; \label{eq:sov_avgQoI}
\end{align}
and
\begin{align}
u_\paavg(t;B) = \sum_{k = 1}^\infty \OmMeas\, M(\psi_k(B)) H(\psi_k(B)) \exp(-\lambda_k(B) t)  \label{eq:sov_paavgQoI_1}
\end{align}
or
\begin{align}
u_\paavg(t;B) = \sum_{k = 1}^\infty \OmMeas (M(\psi_k(B)))^2 \,\dfrac{\lambda_k(B)}{B \gamma} \exp(-\lambda_k(B) t) \; , \label{eq:sov_paavgQoI}
\end{align}
where $M$ and $H$ are defined in \eqref{eq:Mdef} and \eqref{eq:Hdef}, respectively; recall also that $\gamma \equiv |\partial\Omega|/\OmMeas$ is the (nondimensional) surface area to volume ratio as given in \eqref{eq:gammadef}.
\end{proposition}
We provide a sketch of the proof in Appendix \ref{App:sketch_sov_QoI}.

We note that the QoI, in particular the domain-mean QoI, can serve in either forward or inverse mode. In forward mode, given $B$, we wish to evaluate $u_\avg(t;B), 0 \le t \le \tf$. In inverse mode, given a target domain-mean QoI value, $u^*_\avg$, we wish to evaluate $\tau_\avg(B)$ such that 
\begin{align}
u_\avg(\tau_\avg(u_\avg^*;B);B) = u_\avg^* \; .
\end{align}
For example, if $u_\avg^* = \exp(-1)$, then we might interpret $\tau_\avg$ as the time constant of the system. In some cases $\tau_\avg$ can be a more meaningful and sensitive output than $u_\avg$. The error estimators for $\tau_\avg$ follow directly from the error estimators for $u_\avg$; we shall not further consider the inverse mode in this work.

%% file: esb.tex
\section{Eigenproblem for Small $B$}\label{sec:smallb}

\subsection{Asymptotic Analysis}\label{sec:asymp}

We first introduce some notation. As before, $\psi$ and $\lambda$ shall refer respectively to $\psi_1$ and $\lambda_1$, and $\psiz$ and $\lambdaz$ shall refer respectively to $\psi_1$ and $\lambda_1$ evaluated at $B = 0$. We now further define $\psip$ as the first derivative of $\psi_1$ with respect to $B$, $d\psi/dB (B)$; similarly, we define $\lambdap$ as the first derivative of $\lambda_1$ with respect to $B$, $d\lambda/dB(B )$. We define $\psipz$ as the first derivative of $\psi_1$ with respect to $B$ evaluated at $B = 0$, $\psip(B = 0)$; similarly, we define $\lambdapz$ as the first derivative of $\lambda_1$ with respect to $B$ evaluated at $B = 0$, $\lambdap(B = 0)$. Finally, second derivatives with respect to $B$ are denoted by $''$ and also inherit $^0$ for evaluation at $B = 0$. The existence of the derivatives of $\psi_1$ and $\lambda_1$ with respect to $B$, in fact to all orders, follows from the uniqueness of the first eigenfunction and continuity considerations.

As we shall perturb about $B = 0$, we first summarize the case $B = 0$:
\begin{lemma}[{(}Eigenfunction, Eigenvalue{)} for $B = 0$]\label{lem:Bzero} For $B = 0$ we obtain $\psiz (\equiv \psi_1(B = 0)) = \OmMeas^{-1/2}$ and $\lambdaz (\equiv \lambda_1(B = 0)) = 0$.
\end{lemma} 
Note that for $B = 0$ the Robin problem reduces to the Neumann problem. The proof of Lemma \ref{lem:Bzero} is provided in Appendix \ref{App:sketch_Bzero}.

The small-$B$ analysis will depend on the predominance of the first mode. In preparation, we provide the following lemma:

\begin{lemma}[Amplitude of First Eigenfunction] \label{lem:eigfcnamp} For $B \rightarrow 0$, the square of the first eigenfunction ($\sigma$--weighted) domain mean admits expansion
\begin{align}
\OmMeas (M(\psi(B)))^2 = 1 - B^2 \Upsilon  + \pcalO(B^3) \; , \label{eq:eigamp}
\end{align}
where
\begin{align}
\Upsilon \equiv \int_\Omega \sigma\, \psipz\,\psipz \;  \label{eq:Upsilondef}
\end{align}
is non-negative and independent of $B$.
\end{lemma}
We may thus conclude that, indeed, as $B \rightarrow 0$, the spatial structure is well-approximated by the first mode. We provide a sketch of the proof in Appendix \ref{App:sketch_eigfcnamp}. 

We now provide the main result of this section:
\begin{proposition}[First Eigenvalue Asymptotic Expansion]\label{prop:eigvasymp} The first eigenvalue admits expansion
\begin{align}
\lambda(B;\kappa,\sigma) = B \gamma - \phi B^2 + \pcalO(B^3) \text{ as } B \rightarrow 0\; , \label{eq:lambda_exp}
\end{align}
where $\phi \in \Rzp$ is given by
\begin{align}
\phi(\sigma,\kappa) \equiv \int_\Omega \kappa \bnabla\psipz(\cdot;\sigma,\kappa) \cdot \bnabla \psipz(\cdot;\sigma,\kappa) \;;\label{eq:phidef1}
\end{align}
we recall the definition $\psipz\equiv d\psi_1/dB\text{ evaluated at }B=0$.
\end{proposition} We provide in Appendix \ref{App:sketch_eigvasymp} a sketch of the proof of Proposition \ref{prop:eigvasymp}.

We can then readily develop a first-order approximation to $\lambda(B) (\equiv \lambda(B;\kappa,\sigma))$ as $B \rightarrow 0$:
\begin{proposition}[First Eigenvalue: First-Order and Second-Order Approximations in $B$]\label{prop:lamapps} We introduce the approximations
\begin{align}
\lambdaL{1}(B) &\equiv B \gamma \;,\label{eq:lamapp1}\\
\lambdaL{2}(B) &\equiv B \gamma - \phi B^2 \; , \label{eq:lamapp2}\\
\lambdaLP(B) &\equiv \frac{B \gamma}{1+\phi B / \gamma} \;. \label{eq:lamapp2p}
\end{align}
Here, `P' in $\lambdaLP(B)$ refers to Pad\'e as $\lambdaLP(B)$ is the $[1/1]$ Pad\'e approximant \cite{BandGM} to $\lambda(B)$. These approximations satisfy
\begin{align}
\left| \dfrac{\lambda(B) - \lambdaL{1}(B)}{\lambda(B)}   \right| = \phi \underbrace{(B/\gamma)}_{\Bdunk}  + \pcalO(B^2) \text{ as } B \rightarrow 0\; , \label{eq:lamerr1}
\end{align}
\begin{align}
\left| \dfrac{\lambda(B) - \lambdaL{2}(B)}{\lambda(B)}   \right| =  \pcalO(B^2) \text{ as } B \rightarrow 0 \; , \label{eq:lamerr2}
\end{align}
\begin{align}
\left| \dfrac{\lambda(B) - \lambdaLP(B)}{\lambda(B)}   \right| =  \pcalO(B^2) \text{ as } B \rightarrow 0 \; , \label{eq:lamerr2}
\end{align}
where we recall that $\text{Bi}_{\text{dunk}} \equiv B/\gamma$ is the standard engineering definition of the Biot number in the dunking context \cite{LandL,Cengel,IandD}. The superscript numerals in $\lambdaL{1}$, $\lambdaL{2}$, and $\lambdaLP$ refer to the leading-order of the relative approximation error.
\end{proposition}
The proof of Proposition \ref{prop:lamapps} follows directly from Proposition \ref{prop:eigvasymp}: evaluating the Taylor series around $B = 0$ of the left-hand sides of \eqref{eq:lamerr1} -- \eqref{eq:lamerr2}.

\subsection{Sensitivity Formulation for $\phi$}\label{sec:phi_sensitivity}

We observe that $\phi$ plays a key role in the error estimate for $\lambdaL{1}$ --- and, in the next section, the error estimate for the lumped approximation to the initial-value problem --- and also in the definition of $\lambdaLP$.  In turn, $\phi$ is determined by the function $\psipz$, which we shall denote the sensitivity derivative. In fact, the definition of $\psipz$ as the derivative of the (first) eigenfunction at $B=0$ is not directly amenable to either easy analysis or efficient/stable computation. We thus present an alternative formulation:

\begin{proposition}[Elliptic Equation for Sensitivity Derivative]\label{prop:EEforSensitivityDerivative}
The sensitivity derivative $\psipz$ is the solution to
\begin{align}
- \bnabla \cdot (\kappa \bnabla \psipz) = \OmMeas^{-1/2}\gamma \sigma \text{ in } \Omega \; \label{eq:xieqn_strong1}
\end{align}
subject to boundary condition
\begin{align}
\kappa \partial_n \psipz = -\OmMeas^{-1/2} \text{ on } \partial\Omega \;  \label{eq:xieqn_strong2}
\end{align}
and zero-mean condition
\begin{align}
\int_\Omega \sigma \psipz = 0 \; . \label{eq:xieqn_strong3}
\end{align}
We note that \eqref{eq:xieqn_strong1} -- \eqref{eq:xieqn_strong3} is solvable --- $\sigma$ is of mean unity over $\Omega$ --- and admits a unique solution. We recall that $\partial_n$ denotes the derivative with respect to the outward normal.
\end{proposition}
We shall refer to \eqref{eq:xieqn_strong1} -- \eqref{eq:xieqn_strong3} as the (strong form of the) ``Sensitivity Equation'', or more elaborately ``Equation for Sensitivity of the Eigenfunction to Biot''. We provide in Appendix \ref{App:phieqn} the weak form of the Sensitivity Equation, and in Appendix \ref{App:sketch_EEforSensitivityDerivative} a sketch of the proof of Proposition \ref{prop:EEforSensitivityDerivative}.

%% file: lump.tex
\section{Lumped Formulation}\label{sec:lumped}

In this section we apply engineering thermal resistance concepts to (i) understand the physical significance of $\phi$, and (ii) develop our lumped approximations. Note in this section there is no pretense of rigor: error estimates for the lumped approximations — which serve as justifications for our choices — are provided in Section \ref{sec:qoierr} . In Section \ref{sec:circuit}, we shall provide the approximate quantities introduced in the derivations with $\hat\cdot$ to emphasize the formal nature of the discussion.  The main result of this section is the engineering motivation for the three quantities of interest: (i) the (first-order) `classical' lumped approximation to $u_\avg$, (ii) a second-order lumped approximation to $u_\avg$, and (iii) an approximation to $u_\Delta$ which is null under the (first-order) `classical' theory. These approximations are summarized in Definition \ref{def:lumpedapprox}.

In this section, we first introduce the body average resistance $\dv{\calR_\avg}$: the effective resistance between the body average temperature and the surface temperature. In this resistance analogy, we derive an effective conduction length scale $\Lcond$ based on the standard intrinsic length scale $\calL$ but now corrected for second-order effects. We then move to a thermal circuit analogy to derive the effective (corrected) Biot number $\Bdunkp$. Note both corrections are based on $\phi$. Finally, these corrected quantities are then incorporated in a second-order approximation of the solution, which improves upon the (first-order) `classic' approximation.

\subsection{Body Average Resistance} \label{sec:body-average-resistance}

We first introduce a conduction problem in which (non-uniform) heat generation balances (uniform) heat flux at the boundary:
\begin{align}
\dv{-\bnabla\cdot(k\bnabla T^*)}\ &\dv{=}\ \sigma \dv{q \gamma} \text{ in } \dv{\Omega}  \label{eq:strong1}, \\
\dv{k \partial_n T^*}\ &\dv{= -q} \text{ on } \dv{\partial\Omega} \label{eq:strong2}, \\
\dashint_{\Omega} \sigma\dv{(T^* - T_{\text{ref}})}\ &\dv{= 0} \label{eq:strong3}\,;
\end{align}
here $\dv{T_\text{ref}}$ is a reference temperature. We can demonstrate that equations \eqref{eq:strong1} -- \eqref{eq:strong3} are solvable: we integrate \eqref{eq:strong1} over $\dv{\Omega}$, apply \eqref{eq:strong2}, and recall the definition of $\dv{\gamma}$ \eqref{eq:gammadef} of Section \ref{sec:equations}. We now define the associated body average thermal resistance as
\begin{align}
\dv{\mathcal{R}_{\text{avg}}} &\ \dv{\equiv}\ \dfrac{ \dashint_\Omega \sigma \dv{T^* -} \dashint_{\partial\Omega} \dv{T^*}}  { \dv{q |\partial\Omega |} } \\
&\ \dv{=}\ \dfrac{ \dashint_\Omega \sigma \dv{(T^* - T_{\text{ref}})} - \dashint_{\partial\Omega} \dv{(T^* - T_{\text{ref}})} } { \dv{q |\partial\Omega |} } \label{eq:seceq},
\end{align}
in which the second equality \eqref{eq:seceq} follows from the normalization $\dashint_\Omega \sigma = 1$ and the definition of $\dashint$. 

We now introduce the nondimensional temperature
\[ u^* \equiv \dv{\dfrac{ T^* - T_{\text{ref}} }{ (\ell q)/k_{\inf}} } \]
substitution of which into equations \eqref{eq:strong1} -- \eqref{eq:strong3} yields
\begin{align}
-\bnabla\cdot(\kappa\bnabla u^*) &= \sigma \gamma \text{ in } \Omega \label{eq:weakeq}, \\
\kappa \partial_n u^* &= -1 \text{ on } \partial\Omega, \\
\int_\Omega \sigma u^* &= 0 \label{eq:zeromean}.
\end{align}
The average resistance then takes the form 
\begin{align}
\dv{\mathcal{R}_{\text{avg}}} &\dv{\equiv \dfrac{ \ell } {k_{\inf} |\partial\Omega |}} \left(-\dashint_{\partial\Omega} u^*\right),
\end{align}
in which we have taken advantage of \eqref{eq:zeromean} to eliminate the integral over $\Omega$. We now choose to write
\begin{align}
\dv{\mathcal{R}_{\text{avg}} = \dfrac{\mathcal{L}_{\text{cond}}}{k_{\inf} |\partial\Omega|}},\label{eq:ravgdef2}
\end{align}
which implicitly defines the nondimensional conduction resistance length scale
\begin{align}
\mathcal{L}_\text{cond}\equiv\dv{\mathcal{L}_{\text{cond}}/\ell} &= -\dashint_{\partial\Omega} u^* \label{eq:lscale}.
\end{align}

We now note from inspection of equations \eqref{eq:weakeq} -- \eqref{eq:zeromean} that $u^* = |\Omega|^{1/2} \psi'^0$ for $\psi'^0$, the solution to \eqref{eq:xieqn_strong1} -- \eqref{eq:xieqn_strong3}. It then further follows from the weak form of \eqref{eq:xieqn_strong1} -- \eqref{eq:xieqn_strong3} (provided in Appendix \ref{App:phieqn}), definition of $\phi$ \eqref{eq:phidef1}, and definition of $\Lcond$ \eqref{eq:lscale} that
\begin{align}
\phi &= - |\Omega|^{-1/2}\int_{\partial\Omega} \psi'^0 \nonumber\\
&= -|\Omega|^{-1}|\partial\Omega|\,\dashint_{\partial\Omega} u^* \nonumber\\
&= \Lcond / \calL\,,\label{eq:phiLcond}
\end{align}
and hence we may write
\begin{align}
\mathcal{L}_{\text{cond}} &= \phi\, \calL.\label{eq:lcond}
\end{align}
We thus observe that $\mathcal{L}_{\text{cond}}$, the correct conduction length scale to represent the body average thermal resistance, is the usual length scale $\mathcal{L}$ amplified by $\phi$. Note that $\mathcal{L}_{\text{cond}}$ will depend, through $\phi$, on $\sigma$, $\kappa$, and of course $\Omega$. 
Our discussion here thus extends, in a quantitative fashion, the usual notion of resistance between two surfaces to the case of resistance from body ``center'' to body boundary.  Furthermore, we note that the definition yields a positive length scale:

\begin{proposition}[Positive Conduction Length Scale, $\mathcal{L}_\text{cond}$]
For any given dunking problem with admissible parameters $\sigma$, $\kappa$, $B$, and $\Omega$; the corrected conduction length scale 
\begin{align}
\Lcond > 0.
\end{align}
\end{proposition}
This proposition can be proven by considering \eqref{eq:phiLcond}, positivity of $\calL$ through positive $\OmMeas$ and $\dOmMeas$, and $\phi$ being a strictly positive quantity ---  which we shall later demonstrate.

\begin{figure}[H]
\centering
\input{fin}
\caption{Fin attached to a square block.}
\label{fig:fin_eg}
\end{figure}
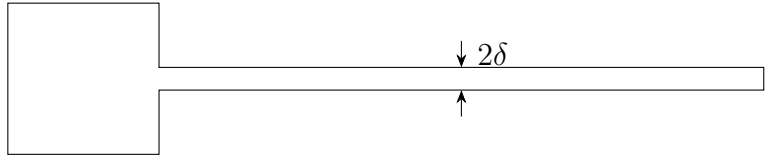

As we shall discuss below, for many ``typical'' domains (slab, cylinder, sphere), $\phi\approx\pcalO(1)$; however, for other domains, for example the finned block shown in Figure \ref{fig:fin_eg}, $\phi$ and thus $\mathcal{L}_\text{cond}/\mathcal{L}$ may be large. This large value of $\phi$ reflects a conjoined domain with disparate length scales associated with each of the component domains --- often   found in engineering applications. We investigate the dependence of $\phi$ on geometry and composition in Section \ref{sec:phi}.

\subsection{Thermal Circuit}\label{sec:circuit}

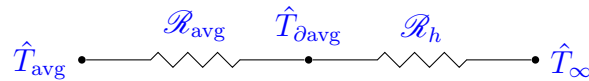
\begin{figure}[H]
  \centering
  \begin{tikzpicture}[circuit ee IEC, set resistor graphic=var resistor IEC graphic]
  \tikzset{point/.style={circle,fill,inner sep=1pt}}
  
  \node (Tavg) at (0,0) [point,label=left:$\dv{\hat{T}_\text{avg}}$] {};
  \node (Tpaavg) at (3,0) [point,label=above:$\dv{\hat{T}_\paavg}$] {};
  \node (Tinf) at (6,0) [point,label=right:$\dv{\hat{T}_\infty}$] {};
  
  \draw (Tavg) to [resistor={info={$\dv{\mathcal{R}_{\text{avg}}}$}}] (Tpaavg)
        (Tpaavg) to [resistor={info={$\dv{\mathcal{R}_h}$}}] (Tinf);
  
  \end{tikzpicture}
  \caption{3-node thermal circuit.}
  \label{fig:thermal_circuit}
\end{figure}

In this section, we continue in order to motivate approximations which will be rigorously analyzed: the main result being a `corrected' ordinary differential equation for $u_\avg$, leading to the second-order lumped approximation. To that effect, we now introduce a thermal circuit, shown in Figure \ref{fig:thermal_circuit}, with three nodes and associated nodal temperatures, $\dv{\hat{T}_{\text{avg}}}$, $\dv{\hat{T}_{\paavg}}$, and $\dv{\hat{T}_{\infty} (= T_{\infty})}$. The nodes associated with $\dv{\hat{T}_{\text{avg}}}$ and $\dv{\hat{T}_{\paavg}}$ are connected by a body average resistance $\dv{\mathcal{R}_{\text{avg}}}$, and the nodes $\dv{\hat{T}_{\paavg}}$ and $\dv{\hat{T}_{\infty}}$ are connected by a standard ``heat transfer coefficient'' resistance, $\dv{\mathcal{R}_h \equiv 1/(h|\partial\Omega|)}$. It is important to note that $\hat\cdot$ refers to approximation to the standard quantities: $\dv{\hat{T}_\text{avg}}$ and $\dv{\hat{T}_\paavg}$ are representative of $\dv{T_\text{avg}}$ and $\dv{T_\paavg}$, but in general not equal to those quantities.

We may now apply standard linear circuit theory to this network. In particular, we can evaluate the heat transfer rate through our circuit as 
\begin{align}
\dv{\hat{Q} = \dfrac{\hat{T}_{\text{avg}} - \hat{T}_{\infty}}{\mathcal{R}_{\text{eq}}} }
\end{align}
and the temperature $\dv{\hat{T}_{\paavg}}$ as
\begin{align}
\dv{\dfrac{\hat{T}_{\paavg} - \hat{T}_{\infty}}{\hat{T}_{\text{avg}} - \hat{T}_{\infty}} = \dfrac{\calR_h}{\mathcal{R}_{\text{eq}}}}, \label{eq:Tpaavg}
\end{align}
where $\dv{\mathcal{R}_{\text{eq}}}$ is the (total, or) equivalent thermal resistance of the circuit,
\begin{align}
\dv{\mathcal{R}_{\text{eq}} \equiv \mathcal{R}_{\text{avg}} + \mathcal{R}_h =}\, \dfrac{\phi \dv{\calL}}{\dv{k_{\inf} |\partial\Omega|}} \dv{\,+ \dfrac{1}{h |\partial\Omega|}} . \label{eq:req_sum}
\end{align}
We thus observe that large $\phi$ increases the body average resistance and thus decreases the heat transfer rate. From another perspective, $\phi$ large decreases the fraction of the temperature difference $\dv{(\hat{T}_\avg - \hat{T}_{\infty})}$ associated with $\dv{\mathcal{R}_h}$ which thus again decreases the heat transfer rate. Note our analysis is perforce approximate since our average resistance assumes uniform source and boundary flux; however, we anticipate that our approximation will be quite good since for small Biot number, the temperature distribution — for example, our first eigenmode — is quite uniform.

We now note that
\begin{align}
\Bdunkp \equiv \phi \Bdunk = \dv{\dfrac{\mathcal{R}_{\text{avg}}}{\mathcal{R}_h}} \label{eq:Bi_dunkp}
\end{align}
which is precisely the usual definition of the Biot number in terms of resistance ratio — but now with the improved conduction length scale provided by \eqref{eq:lcond}. We can then express $\dv{\calR_{\text{eq}}}$ as
\begin{align}
\dv{\mathcal{R}_{\text{eq}} = \dfrac{1}{h|\partial\Omega|}}(1 + \Bdunkp) . \label{eq:R_eq}
\end{align}
Our interest in \eqref{eq:R_eq} is in the limit of small Biot number. In the limit of large Biot number we would write instead $\dv{\mathcal{R}_{\text{eq}}} = \phi\dv{\left[\calL/(k_{\inf}|\partial\Omega|)\right]} \left[1 + (\Bdunkp)^{-1}\right]$.

We now turn to the time-dependent problem. We associate to the $\dv{\hat{T}_{\text{avg}}}$ node a heat capacitance $\dv{(\overline{\rho c})|\Omega|}$ which then directly yields
\begin{align}
\dv{(\overline{\rho c})|\Omega| \partial_t \hat{T}_{\text{avg}} = -\dfrac{\hat{T}_{\text{avg}} - \hat{T}_{\infty}}{\mathcal{R}_{\text{eq}}}} .\label{eq:lumped_circuit_evo}
\end{align}
We now define $\hat{u}^*_{\text{avg}} \equiv \dv{(\hat{T}_{\text{avg}} - \hat{T}_{\infty})/(T_{\text{i}} - \hat{T}_{\infty})}$ and nondimensionalize \eqref{eq:lumped_circuit_evo} to obtain
\begin{align}
\partial_t \hat{u}^*_{\text{avg}} = -\dfrac{B \gamma}{1 + \phi B/\gamma} \hat{u}^*, \quad t > 0 \label{eq:lumped_circuit_ivp}
\end{align}
subject to $\hat{u}^*_{\text{avg}}(t = 0) = 1$ from our uniform initial condition. We can also identify the eigenproblem associated with \eqref{eq:lumped_circuit_ivp}:
\begin{align}
\dfrac{B \gamma}{1 + \phi B/\gamma} \hat{\psi} = \hat{\lambda} \hat{\psi} , \label{eq:lumped_circuit_eig}
\end{align}
where $(\hat{\psi},\hat{\lambda})$ is the single eigenpair associated with our scalar first-order problem.

\subsection{Lumped Approximations}

We begin with the eigenproblem in order to arrive at the lumped approximation. We directly identify from \eqref{eq:lumped_circuit_eig} our first-order approximation — neglecting the $\pcalO(B^2)$ contributions — as $\lambdaL{1} = B \gamma$, and our second-order approximation — including now the $\pcalO(B^2)$ terms — as $\lambdaLP = \frac{B \gamma}{1 + \phi B/\gamma}$. We have already presented and analyzed these approximations in Proposition \ref{prop:lamapps}, respectively. We can now also include a physical interpretation in particular for the second-order approximation: for larger $\phi$, hence larger conduction length scale, the heat transfer coefficient resistance is relatively smaller compared to the conduction (average) resistance, and thus $\dv{\hat{T}_{\paavg}}$  is closer to $\dv{T_\infty}$; thus the heat transfer rate to the environment, $\dv{h |\partial\Omega|(T_{\text{avg}} - T_\infty)}/(1+\Bdunkp)$ is reduced; and hence the time constant — $1/\lambdaLP$ — is increased.

Turning now to the initial value problem, we can directly identify from \eqref{eq:lumped_circuit_ivp} our first-order approximation for the domain average — neglecting the $\pcalO(B^2)$ contributions — as $\uavgL{1} \equiv \exp(-B\gamma t)$, the `classic' approximation — and then the second-order approximation for the domain average — now including the $\pcalO(B^2)$ terms — as $\uavgL{2} \equiv \exp\left(-\frac{B\gamma t}{1 + \phi B/\gamma}\right)$; note again that larger $\phi$, and larger conduction length scale, leads to longer time scales for equilibration. Finally, it follows from \eqref{eq:Tpaavg} that $\uDeltaL \equiv \frac{\phi B}{\gamma}\left(1 + \frac{\phi B}{\gamma}\right)^{-1}$ (independent of time) --- we use superscript $2$ in reference to its derivation from a second-order approximation. We summarize these three results in Definition \ref{def:lumpedapprox}. We also note that alternative approaches exist to derive the lumped approximation $\uavgL{1}$; one such approach, Galerkin projection of the PDE, is presented in Appendix \ref{App:sketch_galerkinlumped}: closely related to the Rayleigh quotient (Appendix \ref{App:sketch_Bzero}) for a constant function.
\begin{definition}[Lumped Approximations: Summary]\label{def:lumpedapprox}  The first and second-order lumped approximations to the QoI, $u_\avg$ and $u_\Delta$, are given by
\begin{align}
\uavgL{1}(t;B) &\equiv \exp(-B\gamma t) \; , \label{eq:u1avg}\\
\uavgL{2}(t;B) &\equiv \exp\left(-\dfrac{B\gamma t}{1 + \phi B/\gamma}\right) \; , \label{eq:u2avg}\\
\uDeltaL(B) &\equiv \dfrac{\phi B}{\gamma}\left(1 + \frac{\phi B}{\gamma}\right)^{-1} \; . \label{eq:u2Delta}
\end{align}
\end{definition} 
The associated error analysis for these approximations is provided in Section \ref{sec:qoierr}.

We note that a second-order lumped approximation is proposed already in \cite{ODD} based on the eigenvalue expansion (in \cite{ODD}, for the particular case of a homogeneous sphere). This solution based on second-order Taylor approximation of the principal eigenvalue is different from our solution, $\uavgL{2}$, based on the second-order $[1,1]$ Pad\'e approximation \cite{BandGM} to the principal eigenvalue: the two proposals agree to second order in $B$ as $B$ tends to zero; however, the Taylor approximation is only meaningful for $\frac{\phi B}{\gamma} < 1$ whereas the Pad\'e approximation yields stable results for any $B$. It is interesting that the engineering resistance arguments directly yield the Pad\'e approximation; see also \cite{Ostro2008}. In a similar fashion, the Pad\'e approximation for $\uDeltaL$ approaches $1$, the correct limit, as $B\to \infty$.

It is also important that calculation of the dimensional quantity $\dv{u_\avg}$ requires only $\u_\avg$, $\dv{\ell}$, $\dv{\overline{\rho c}}$, and $\dv{\overline{h}}$. Hence again we need only the average, not local, heat transfer coefficient. It is the average heat transfer coefficient, $\dv{\overline{h}}$, and not the local heat transfer coefficient, $\dv{h}$, which is most readily measured and most easily archived. Furthermore, for heterogeneous bodies, $\dv{\overline{\rho c}}$ can be evaluated solely on the basis of constituent properties and volume fractions --- the details of the material distribution are {\em not} required.

%% file: fin.tex
\begin{tikzpicture}
  \draw (0,0) -- (0,2) -- (2,2) -- (2,1.15) -- (10,1.15) -- (10,0.85) -- (2,0.85) -- (2,0) -- cycle;
  \draw[-{Stealth[length=5pt]}] (6,1.5) -- (6,1.15) node[midway, right, xshift=2pt] {$2\delta$};
  \draw[-{Stealth[length=5pt]}] (6,0.5) -- (6,0.85);
\end{tikzpicture}

%% file: phi.tex
\section{Solutions of the Sensitivity Equation: $\phi$}\label{sec:phi}

We now turn to a detailed analysis of $\psipz$ and $\phi$. The quantity $\phi$ appears in our second-order approximations and in all our error estimates.

We first introduce some notation relevant to our study of the geometric dependence of $\phi$. Let 
\begin{align}
\Xi(x_{\text{pg}}^1\in \RR^2, x_{\text{pg}}^2\in \RR^2,\ldots,x_{\text{pg}}^{n_{\text{pg}}}\in \RR^2)  \label{eq:polydef}
\end{align}
denote the polygon $\Omega_{\text{pg}} \subset \RR^2$ with boundary
\begin{align}
\partial\Omega_{\text{pg}} \equiv \bigcup_{i = 1}^{n_{\text{pg}}} \;\overline{x_{\text{pg}}^i\; x_{\text{pg}}^{i+1}} \; . \label{eq:polybndry}
\end{align}
We shall assume that the vertices $x_{\text{pg}}^i, 1 \le i \le n_{\text{pg}}$, are chosen such that the interiors of the segments $\overline{x_{\text{pg}}^i\; x_{\text{pg}}^{i+1}}, 1 \le i \le n_{\text{pg}}$, are mutually non-intersecting. The polygons are simply connected but not necessarily convex or star-shaped.

\subsection{General Properties of $\phi$}\label{sec:phiprops}

We summarize here several properties of $\phi(\sigma,\kappa;\Omega)$ which can be gainfully applied in the evaluation of $\phi$ and the application of $\phi$ to error analysis:
\begin{proposition}[$\phi$ is Strictly Positive]\label{prop:phiprop1} For any given admissible $\sigma$ and $\kappa$, $\phi(\sigma,\kappa) > 0$.
\end{proposition}
\begin{proposition}[$\phi$ is Scale-Invariant]\label{prop:phiprop2}  We are given $\Omega_1 \subset \RR^d$ and a translation-rotation-uniform dilation map of $\Omega_1$, ${\Omega}_2 \equiv \{ \MM(x) \, |\, x \in \Omega_1 \,\}$; here $\MM(x) = \TT + \QQ\, x + \alpha\,\II\,x$, where $\alpha$ is a positive real scalar, $\TT\in\mathbb{R}^d$ is a translation vector, $\QQ \in \RR^{d\times d}$ is a rotation matrix, and $\II \in \RR^{d\times d}$ is the identity. For given $\sigma \equiv \sigma_1, \kappa \equiv \kappa_1$, $\Omega \equiv \Omega_1$ we evaluate $\phi_1 \equiv \phi(\Omega_1)$ from Proposition \ref{prop:EEforSensitivityDerivative}; for given $\sigma \equiv \sigma_2 \equiv \sigma_1 \circ \MM^{-1}, \kappa \equiv \kappa_2 \equiv \kappa_1 \circ \MM^{-1}$, $\Omega \equiv \Omega_2$ we evaluate $\phi_2 \equiv \phi(\Omega_2)$. Then, for any $\alpha$, $\TT$, and $\QQ$, $\phi_1 = \phi_2$.  
\end{proposition}
\begin{proposition}[$\phi$ is Monotonic in $\kappa$]\label{prop:phiprop4} For $\kappa_1$ and $\kappa_2$ such that $\kappa_1 \le \kappa_2, \forall x \in \Omega$, $\phi(\sigma,\kappa_1) \ge \phi(\sigma,\kappa_2)$. 
\end{proposition}
We henceforth indicate the dependence of $\phi$ as $\phi(\sigma,\kappa;\Omega_{\text{ref}})$ and, from Proposition \ref{prop:phiprop2}, we can apply this value of $\phi$ to any domain $\Omega$ which is translation-rotation-uniform dilation map of $\Omega_{\text{ref}}$. We can further state, from Proposition \ref{prop:phiprop4}, that $\phi(\sigma,\kappa;\Omega_{\text{ref}}) \le  
\phi(\sigma,1;\Omega_{\text{ref}})$ for any $\kappa \ge 1$; given the latter bound, we shall often restrict attention to $\kappa = 1$ even if the result can be developed for general $\kappa$. We provide in Appendix \ref{App:sketch_phiprops} sketches of the proofs of Propositions \ref{prop:phiprop1} -- \ref{prop:phiprop4}.

Proposition \ref{prop:phiprop4} will permit us to develop our error estimators  even for the case --- quite common in engineering practice --- of heterogeneous material composition. In the next proposition we provide also some insight into the dependence of $\phi$ on $\sigma$ and in particular on deviations of $\sigma$ from uniform:

\begin{proposition}[$\phi$ Dependence on $\sigma$]\label{prop:phisigma} We consider here the case in which $\kappa = 1$ and we furthermore assume that the domain $\Omega$ is fixed; we may thus abbreviate $\phi$ as simply $\phi(\sigma)$. Let us introduce
\begin{align}
\mu \equiv \inf_{w \in Z_0(1)} \dfrac{a_0(w,w;1)}{m(w,w;1)} \; . \label{eq:mudef}
\end{align}
We further define
\begin{align}
\delta(\sigma) \equiv \left( \dfrac{\gamma^2}{\mu} \dashint_{\Omega} (\sigma - 1)^2 \right)^{1/2} \; , \label{eq:deltadef}
\end{align}
in terms of which we may then bound $\phi(\sigma)$ as
\begin{align}
\phi(\sigma) \le  \phi_\text{UB}(\sigma)\equiv \left( \left(\phi(1)\right)^{1/2} + \delta(\sigma)\right)^2  \; , \label{eq:phisigmabound}
\end{align}
where $\phi(1)$ refers to $\phi$ for $\sigma=1$. We recall that $\dashint$ refers to the average of the integrand over the domain.
\end{proposition} 
We provide in Appendix \ref{App:sketch_phisigma} a proof of Proposition \ref{prop:phisigma}. 

We note that $\gamma^2/\mu$ is scale invariant. We further note that $\mu$ is simply the second (or any second) eigenvalue of the Neumann Laplacian eigenproblem on the domain $\Omega$. It thus follows that $\mu$ may be estimated very efficiently, and often with sufficient precision, by application of the Payne-Weinberger \cite{PandW} lower bound for convex domains 
\begin{align}
  \mu \ge \mu_\text{LB}\equiv \pi^2\calD_\Omega^{-2} \; \label{eq:muPW}
\end{align}
where
\begin{align}
  \calD_\Omega = \sup_{x\in\Omega,y\in\Omega}|x - y| \; \label{eq:Ddef} 
\end{align}
requires calculation of only the diameter of our domain. For later numerical results in Section \ref{sec:comp_results}, we shall numerically evaluate $\mu_h(=\mu)$; in general, $(\cdot )_h$  shall refer to finite element approximation of $(\cdot)$.  In actual practice, the lower bound $\mu_\text{LB}$ is more convenient.

Proposition \ref{prop:phisigma} demonstrates that in fact there are many physical situations in which $\phi(\kappa = 1,\sigma = 1)$ --- easily determined without detailed knowledge of the heterogeneous material composition --- can serve as a surrogate for $\phi(\kappa = 1,\sigma)$. We note in particular that either if the two materials are similar in thermophysical properties {\em or} if the two materials occupy commensurate volume (fraction), then $ \dashint_{\Omega} (\sigma - 1)^2$ is either small or of order unity. As an example of the latter, we consider two materials of very different properties each of which occupies half of $\Omega$:

\begin{proposition}[Equi-Volume-Fraction Two-Phase Composite]\label{prop:evfsigma}
We consider two materials (of different thermophysical properties) each of which constitutes one half of the volume of our domain $\Omega$. Then $\dashint_{\Omega}(\sigma -1)^2  \le 1$. 
\end{proposition}
Note that even in this equi-volume-fraction case, $\delta(\sigma)$ may still be large due to geometry effects through $\mu$.
We provide in Appendix \ref{App:sketch_evfsigma} a proof of Proposition \ref{prop:evfsigma}.

\subsection{Closed-Form Results}\label{sec:closedform}

We first present in Table \ref{table:canonicalphi} the values of $\phi(\kappa = 1,\sigma = 1;\Omega_{\text{ref}})$ for the canonical reference domains: interval in $\RR^1$, disk in $\RR^2$, and sphere in $\RR^3$. We choose $\kappa = 1$ and $\sigma = 1$ for two reasons: first, often in practice we encounter homogeneous properties; second, thanks to Proposition \ref{prop:phiprop4} and Proposition \ref{prop:phisigma}, we can bound the effect of deviations from uniform properties in terms of readily measured quantities. Note the sphere calculation is presented in \cite{ODD} for the case of homogeneous properties ($\sigma = \kappa = 1$) by direct eigenvalue perturbation expansion; in Appendix \ref{App:phispherecalc} we reproduce this homogeneous sphere result based on the (strong form of the) Sensitivity Equation (Proposition \ref{prop:EEforSensitivityDerivative}). 

\begin{table}[h!]
\begin{center}
\caption{Values of $\phi(\sigma=1,\kappa=1;\Omega_{\text{ref}})$ for Canonical Domains in Dimension $d\in\{1,2,3\}$. Note that we include a nominal length scale for each geometry, however it follows from scale invariance, Proposition \ref{prop:phiprop1}, that $\phi(\sigma=1,\kappa=1;\Omega_{\text{ref}})$ is independent of  $\dv{\ell}$ (Proposition \ref{prop:phiprop2}).}  \label{table:canonicalphi}
\begin{tabular}{l|c}
\multicolumn{1}{c|}{$\dv{\Omega_{\text{ref}} \equiv \{x \in \RR^d\,|\, |x| \le \ell\}}$} & $\phi$ \\
\hline\hline
$\quad\quad d = 1$: Interval & 1/3\\
$\quad\quad d = 2$: Disk  & 1/2\\
$\quad\quad d = 3$: Sphere & 3/5
\end{tabular}
\end{center}
\end{table}

It is also possible to analyze selected triangles in closed form. We present here a limiting case for a right triangle which is algebraically simple but illustrates the key point:
\begin{proposition}[Evaluation of $\phi$ for Selected Triangles]\label{prop:phirt} We consider the right triangle $\Omega^{\text{right triangle}} \subset \RR^2$ given by 
\begin{align}
\Xi^{\text{right triangle}}(W)(x_{\text{pg}}^1 \equiv (0,0), x_{\text{pg}}^2 \equiv (W,0), x_{\text{pg}}^{3} \equiv (0,1))  \label{eq:rtriangledef} \; ,
\end{align}
for $W \in \RR_+$. Then
\begin{align}
\phi(\sigma=1,\kappa=1;\Xi^{\text{right triangle}}(W)) \sim \dfrac{2}{3} W^{-2} \text{ as } W \rightarrow 0 \; . \label{eq:rtrianglephi}
\end{align}
As by-products of this analysis, we find that
\begin{align}
\phi(\sigma=1,\kappa=1;\Xi^{\text{right triangle}}(1)) =\frac{4}{3}
\end{align}
and
\begin{align}
\phi(\sigma=1,\kappa=1;\Xi^{\text{equilateral triangle}}) = 1 \; .
\end{align}
Recall that, as always, these results are scale invariant (Proposition \ref{prop:phiprop2}).
\end{proposition}
We provide a proof of Proposition \ref{prop:phirt} in Appendix \ref{App:sketch_phirt}.

We note that, whereas the canonical domains of Table \ref{table:canonicalphi} all yield $\phi$ on the order of unity, a triangular domain can, in the limit $W \rightarrow 0$, produce a very large $\phi$. In the next section we will encounter other domains (and thermophysical property limits) which exhibit large $\phi$. And we subsequently demonstrate, in the context of error analysis, that the standard engineering criterion for applicability of the lumped approximation --- $\Bdunk$ small --- implicitly assumes $\phi \approx \pcalO(1)$.  This relates to the physical (thermal circuit) analysis of Section \ref{sec:circuit} (i.e., $\phi \approx \pcalO(1) \rightarrow \Bdunk \approx \Bdunkp$).

Tensorized domains are important and prevalent for the dunking problem class and we show that$\phi$ may be readily evaluated for such domains. We provide here a particular but illustrative case (and for simplicity, we consider $\kappa = \sigma = 1)$:
\begin{proposition}[$\phi$ for Tensorized Domains]\label{prop:phitendom} We introduce, for a given two-dimensional domain $\Omega^{\text{2D}} \subset \RR^2$, the tensorized (or extruded) three-dimensional domain $\Omega \subset \RR^3 \equiv \Omega^{\text{2D}}\times (0,L_3)$.
Then
\begin{align}
\phi(1,1;\Omega) = \phi_{\text{2D}}(1,1;\Omega^{\text{2D}}) + \dfrac{1}{3} \; , \label{eq:phitendom1}
\end{align}
where $\phi_{\text{2D}}$ refers to the value of $\phi$ for a two-dimensional domain (in our case here, $\Omega^{\text{2D}}$). 
As always, our result is scale invariant with respect to $\Omega^{\text{2D}}$ (Proposition \ref{prop:phiprop2}).
\end{proposition}
We provide a proof of Proposition \ref{prop:phitendom} in Appendix \ref{App:sketch_phitendom}.  Note that $\phi$ of Proposition \ref{prop:phitendom} does not depend on $L_3$, and in particular there is no distinction between small $L_3$ and large $L_3$. This result is counter-intuitive, however it is important to recall that $\phi$ is a coefficient of the second-order contribution to the eigenvalue. Recall also our assumption that the heat transfer coefficient (Robin coefficient) is uniform over the entire boundary $\partial\Omega$.

Finally, we consider two illustrative applications of Proposition \ref{prop:phitendom} in concert with respectively Table \ref{table:canonicalphi} and Proposition \ref{prop:phirt}: for a cylinder of radius $R$ and length $L_3$, $\phi = 1/2 + 1/3 = 5/6$; for a prismatic right triangle $\Xi^{\text{right triangle}}(W)\times (0,L_3)$, $\phi \sim \frac{2}{3}W^{-2}$ as $W \rightarrow 0$. A small variation of Proposition \ref{prop:phitendom} (in conjunction with Table \ref{table:canonicalphi}) readily reveals that, for a rectangle of dimensions $L_1 \times L_2$ and an orthogonal parallelepiped of dimensions $L_1 \times L_2 \times L_3$, $\phi = 2/3$ and $\phi = 1$, respectively.  

\subsection{Foundations for Alternative Evaluation Strategies}\label{sec:altstrat}
As discussed in the introduction, our lumped approximation, a simple exponential, has many attractive features: transparency — explicit parameter dependence; flexibility — ready incorporation into design or optimization or system analysis; rapid evaluation — just a few standard elementary function evaluations. 

However, our error estimates shall require the new quantity $\phi$: we already observe in Proposition \ref{prop:eigvasymp} that the relative error estimate for the first-order lumped approximation to the first eigenvalue is given by $\Bdunkp(\equiv\phi\, \Bdunk)$; we shall similarly observe in the Section \ref{sec:qoierr} that the error in the first-order lumped approximation to $u_\mathrm{avg}$ is similarly given by $\Bdunkp/e$ (Proposition \ref{Prop:N1}). We now understand that $\phi$ may, in fact, be large — and hence the implicit engineering assumption (say, $\phi = 1$) can not, in general, be trusted: $\phi$ must be evaluated. We further know that, in general, $\phi$ cannot be obtained in closed form except for a few simple geometries. We conclude that we must, explicitly or implicitly, resort to computational treatment of the partial differential equation \eqref{eq:phiellp_0} for $\psi'^0$ in order to evaluate $\phi$ from \eqref{eq:phiellp_a}. Do we thus lose the advantages of the lumped approximation? In the next paragraphs we show that, in fact, we do not lose these advantages.

We first note that although the set of geometries amenable to closed-form evaluation of $\phi$ (Table \ref{table:canonicalphi}, Proposition \ref{prop:phirt}) is small, these geometries are of practical interest as they often arise in both engineered and natural systems, and furthermore, these geometries can serve to approximate more general shapes.  We also emphasize that this set of geometries is further enriched by Proposition \ref{prop:phitendom}, providing closed-form $\phi$ expressions for finite circular cylinder, orthogonal parallelepiped, and certain triangular prisms.

However, complex geometries not amenable to simple approximations often arise in engineering applications. In practice, we do not lose the advantages of the lumped approximation: a sufficiently accurate approximation to our single elliptic partial differential equation can be very rapidly obtained by adaptive finite element methods (see Section \ref{sec:comp_results}); furthermore, more recent software advances, in particular containers and microservices, permit cloud-based deployment to virtually any user — with no issues of licenses or installation or compilation or access. Nevertheless, the lumped approximation requires only the volume and surface area of the body, whereas our partial differential equation requires a full geometry (and material composition) description — often an inconvenience. We may thus wish to supplement direct evaluation of $\phi$ with a more convenient alternative strategy in which $\phi$ is replaced by an inexpensive surrogate.

As regards geometry, there are two options. In both cases, we are given the desired domain, $\Omega$, and a dictionary of $(\Omega,\phi)^k$ pairs. In the first option, we estimate $\phi$ as $\phi^k$ for the element $k$ of the dictionary for which $\text{dist}(\Omega,\Omega^k)$ (with suitable registration) is smallest. In the second option, we identify a few features which permit approximate separation of our dictionary entries for $\phi$ smaller and $\phi$ larger; in the offline stage we train the classifier, and then in the online stage we query the classifier for any given desired domain, $\Omega$. We are thus presented with two mathematical questions: (for the first option) How should we choose $\mathrm{dist}$, in particular to provide Lipschitz continuity for $\phi$ — $\mathrm{dist}(\Omega^i,\Omega^j)$ small implies $|\phi^i - \phi^j|$ small \cite{FaberKrahn}?; (for the second option) How should we choose our features to provide adequate separation between the $\phi$ smaller and $\phi$ larger subsets of the dictionary (or, more generally, training set)?

We cannot yet provide rigorous answers. However, we have conjectures which we will then support (in the next section) with computational evidence. A subsequent manuscript will elaborate on the material presented here and provide (and assess) specific surrogate techniques.
\begin{itemize}
\item[]
For the distance we propose 
\begin{align}
\text{dist}(\Omega^i,\Omega^j) \equiv C_1 d_{\mathrm{Hausdorff}}(\Omega^i,\Omega^j) + C_2 \left| |\partial\Omega^i| - |\partial\Omega^j| \right|,\label{eq:distdef}
\end{align}
where $d_{\mathrm{Hausdorff}}$ is the Hausdorff distance \cite{RockWets98} and $C_1$ and $C_2$ are positive real constants. 
\item[]
The geometric features should, in general, increase with $\phi$. For the first feature, we propose
\begin{align}
\Feat 1 &\equiv \text{Modified Weighted Inverse Angle (MWIA)},\label{eq:feat1def}
\end{align}
$\Feat 1$ is complicated and we thus restrict attention to $\Omega \equiv \Omega^{\text{2D}} \subset \mathbb{R}^2$ (or associated extruded domains $\Omega \times (0,W)$ in $\mathbb{R}^3$) and furthermore polygonal domains $\Omega^{\text{2D}}$. In $\Feat 1$ we consider all polygonal domains which are sufficiently close to $\Omega^{\text{2D}}$ in $\text{dist}$; for each such sufficiently close domain, we sum the inverse of the squared angles over all vertices weighted by the relative perimeter of the edges which share the vertex. The `modified' qualifier indicates that a pre-processing procedure in which small segments (compared to its neighbors) are removed to account for `smoothed' sharp angles which may fool an `unmodified' procedure. In short, $\Feat 1$ will be large if some good fraction of the boundary corresponds to small interior angles.
\item[]
For the second feature, we propose
\begin{align}
\Feat 2 &\equiv \text{InRadius $\gamma$}.\label{eq:feat2def}
\end{align}
$\Feat 2$ is much simpler than $\Feat 1$, since the InRadius, the radius of the largest inscribed ball, is readily evaluated or estimated. In short, $\Feat 2$ will be large if there are extended regions which contribute substantially to $|\partial\Omega|$ but not to $|\Omega|$.
\end{itemize}
Examples will be provided in the next section.

\begin{figure}[H]
\begin{center}
\includegraphics[width=0.6\textwidth]{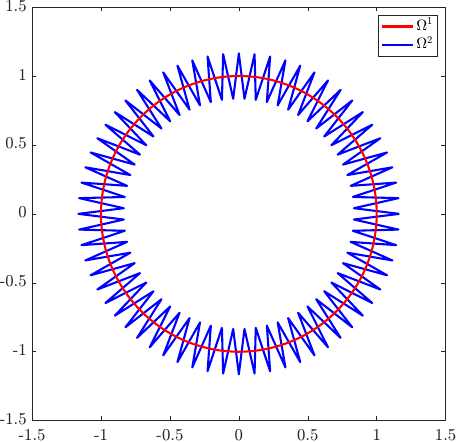}
\caption{Illustrative domains $\Omega^1$ and $\Omega^2$.}
\label{fig:omega12}
\end{center}
\end{figure}

The choice of distance $\text{dist}$ and also $\Feat 1$ and $\Feat 2$ can be understood in terms of the simple resistance arguments presented in Section \ref{sec:lumped}. As an example, consider two domains illustrated in Figure \ref{fig:omega12}: $\Omega^1$ is the unit-radius disk; $\Omega^2$ is the unit-radius disk but with a highly oscillatory (in space) boundary --- a version of the latter, a gear, shall be studied in the next section. We observe that for $\Omega^2$, $\mathcal{L} = 1/\gamma$ will be very small, even though clearly the relevant conduction length scale is the radius (= 1): to reflect the latter, $\phi$ for $\Omega^2$ must be large compared to $\phi = 1/2$ for our disk, $\Omega^1$. And, as we would hope, $\text{dist}(\Omega^1,\Omega^2)$ will be large — thanks to the term associated with the difference in boundary measure. And, again as we would hope, for $\Omega^1$, both $\Feat 1$ and $\Feat 2$ are not large, where as for $\Omega^2$, both $\Feat 1$ and $\Feat 2$ are large. We shall observe other examples in which only one of the two features will be large — hence the need for both features.

Furthermore, we provide here a lower bound for $\phi$ based on $\Feat 2$:
\begin{proposition}[Lower Bound for $\phi$]\label{prop:philb}
For a dunking problem in $\Omega$, with material properties $\sigma$ and $\kappa$, $\phi$ is bounded below by a function $\FF(\sigma,\kappa,\Omega)$: $\phi(\sigma,\kappa,\Omega)\ge\FF(\sigma,\kappa,\Omega)$. For a two dimensional domain $\Omega^\text{2D}$,
\begin{align}
\FF(\sigma,\kappa,\Omega^\text{2D}) &= \frac{\underline{\sigma}^2}{\overline{\kappa}}\dfrac{\pi}{8}\frac{\left(\text{InRadius}(\Omega^\text{2D})\right)^2}{|\Omega^\text{2D}|}\left(\Feat{2}(\Omega^{\text{2D}})\right)^2,\label{eq:ffdef}
\end{align}
and for a three dimensional domain $\Omega^\text{3D}$,
\begin{align}
\FF(\sigma,\kappa,\Omega^\text{3D}) &= \frac{\underline{\sigma}^2}{\overline{\kappa}}\dfrac{4 \pi}{45}\frac{\left(\text{InRadius}(\Omega^\text{3D})\right)^3}{|\Omega^\text{3D}|}\left(\Feat{2}(\Omega^\text{3D})\right)^2;\label{eq:ff3def}
\end{align}
where $\underline{\sigma} \equiv \essinf_{x\in\Omega}\sigma(x)$ and $\overline{\kappa} \equiv \sup_{x\in\Omega}\kappa(x)$.
\end{proposition}

We provide in Appendix \ref{App:sketch_philb} sketch of the proofs of Proposition \ref{prop:philb}.

To address material heterogeneity --- non-uniform thermophysical properties --- our Proposition \ref{prop:phisigma} suggests very simple bounds that in practice can yield sufficiently sharp results; we shall provide numerical examples in the next section demonstrating the effectiveness of the proposed bounds.

%% file: comp.tex
\subsection{Representative Computational Results} \label{sec:comp_results}

In this section we provide some illustrative computational results. We exploit standard $\PP_2$ finite element approximation with adaptive refinement \cite{Yano} and extrapolation based error estimation. In the case of curved geometries we consider isoparametric approximation and uniform refinement. In the below, we report the finite element approximation $\phi_h$ and also the associated finite element error estimator $e^\phi_h$; in  practice (and in theory as $h\rightarrow 0$), $|\phi - \phi_h| \lessapprox e^\phi_h$. The error estimator is based on the solution convergence rate for a successively mesh refinement \cite{BandLR,SandB}. In all cases we consider domains $\Omega_{\mathrm{2D}} \subset \RR^2$, but recall we can readily extend these results to extruded three-dimensional domains thanks to Proposition \ref{prop:phitendom}. Note the finite element approximation is based on the saddle formulation of Appendix \ref{App:phisaddle}: the latter effectively incorporates the constraint associated with $Z_0(\sigma)$ and furthermore admits a sparse solution.

\begin{figure}[H]
  \centering
  \includegraphics[width=.6\textwidth]{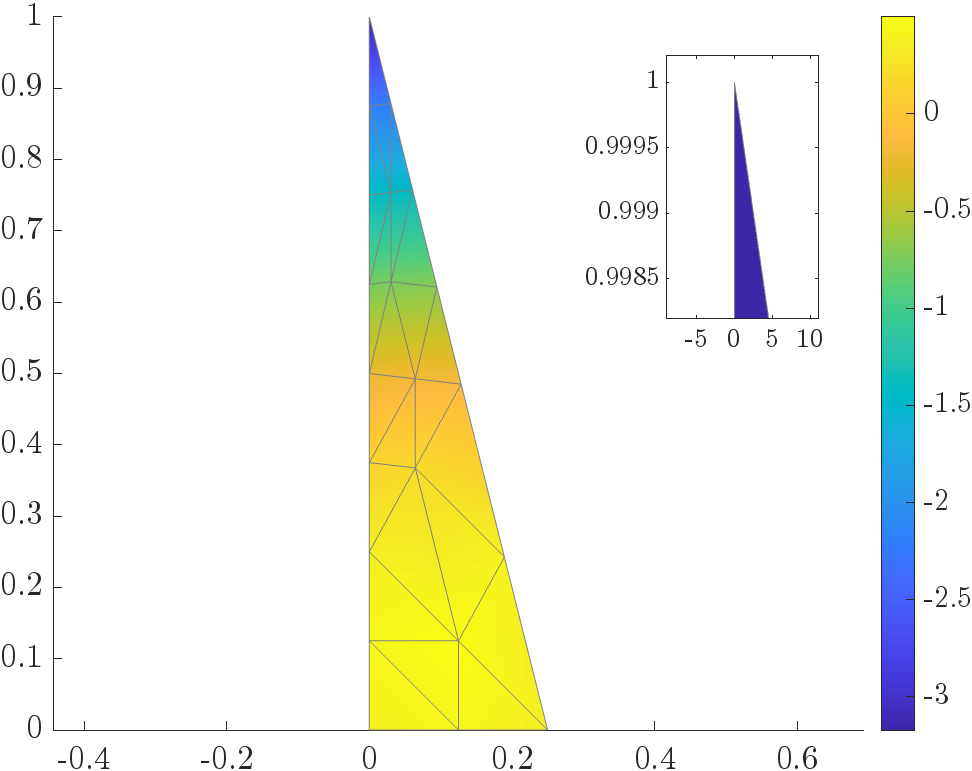}
  \caption{Geometry Small-Angle Right Triangle, denoted in short as SART: solution $\psi'^0$ and final (adaptively refined) finite element mesh. The geometry is given by \eqref{eq:rtriangledef} with $W=1/4$.   Note we may construct from $\Omega^{2\mathrm{D}}$ associated extruded domains in $\RR^3$.}\label{fig:sart}
\end{figure}

\begin{figure}[H]
  \centering
  \includegraphics[width=.6\textwidth]{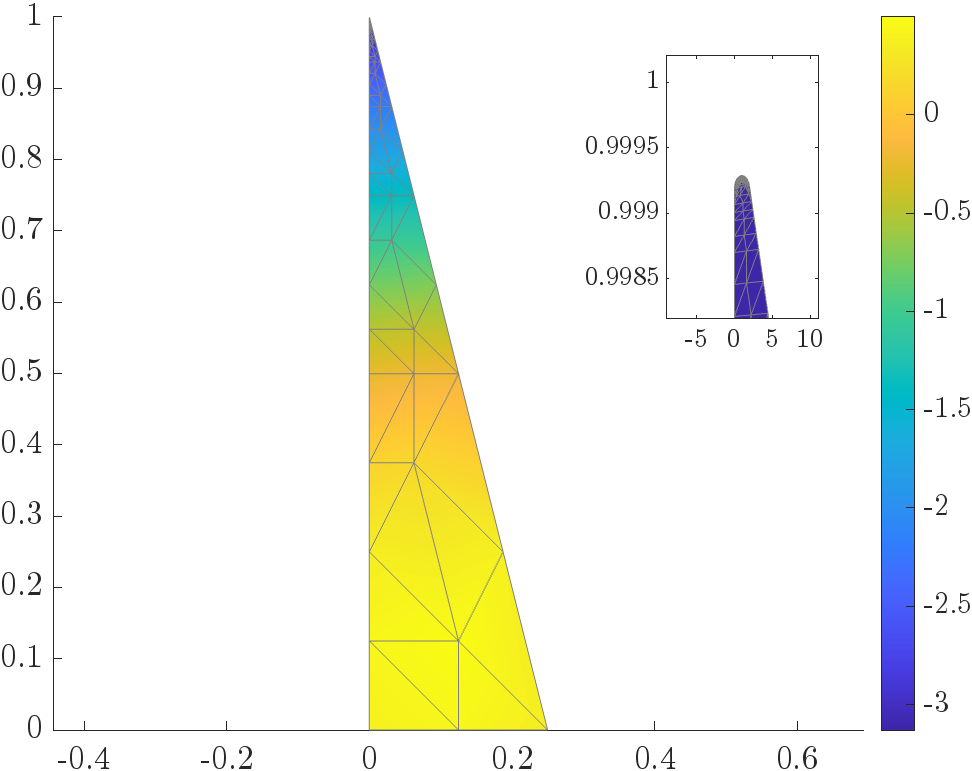}
  \caption{Geometry Small-Angle Right Triangle with Fillet (SART with a $10^{-4}$ radius fillet), denoted in short as SARTF: solution $\psi'^0$ and final (adaptively refined) finite element mesh. Note we may construct from $\Omega^{2\mathrm{D}}$ associated extruded domains in $\RR^3$.}\label{fig:sartf}
\end{figure}

\begin{figure}[H]
    \centering
  \includegraphics[width=.6\textwidth]{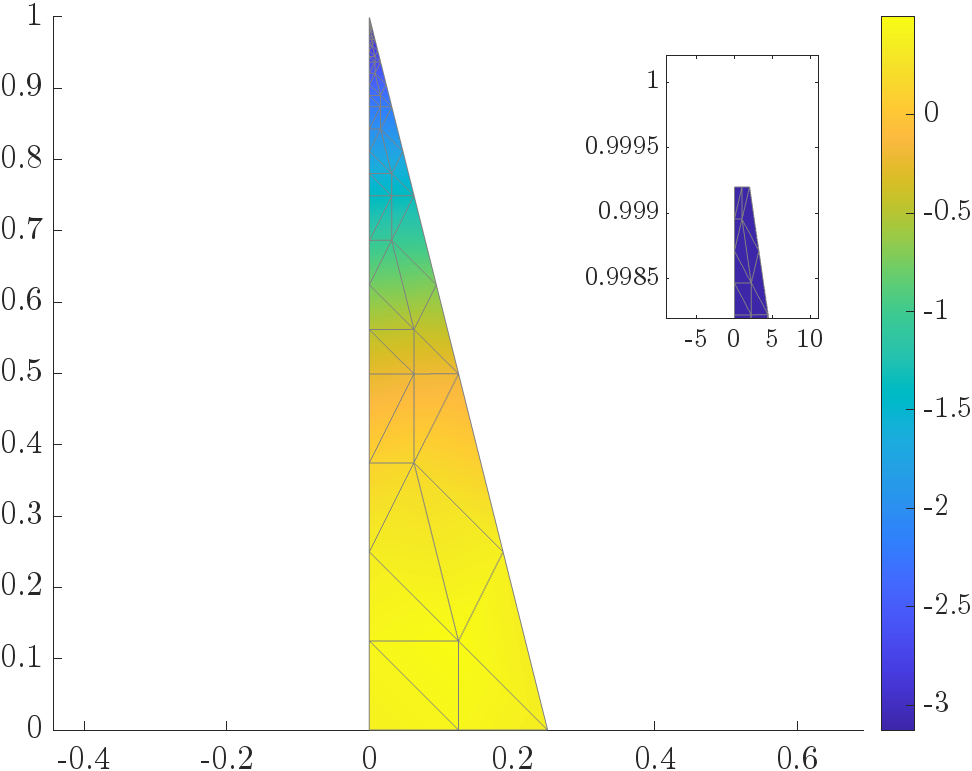}
  \caption{Geometry Small-Angle Right Triangle with Cut ($2\times 10^{-4}$ horizontal cut modification to SART), denoted in short as SARTC: solution $\psi'^0$ and final (adaptively refined) finite element mesh. Note we may construct from $\Omega^{2\mathrm{D}}$ associated extruded domains in $\RR^3$.}\label{fig:tznrt}
\end{figure}

\begin{figure}[H]
    \centering
  \includegraphics[width=.6\textwidth]{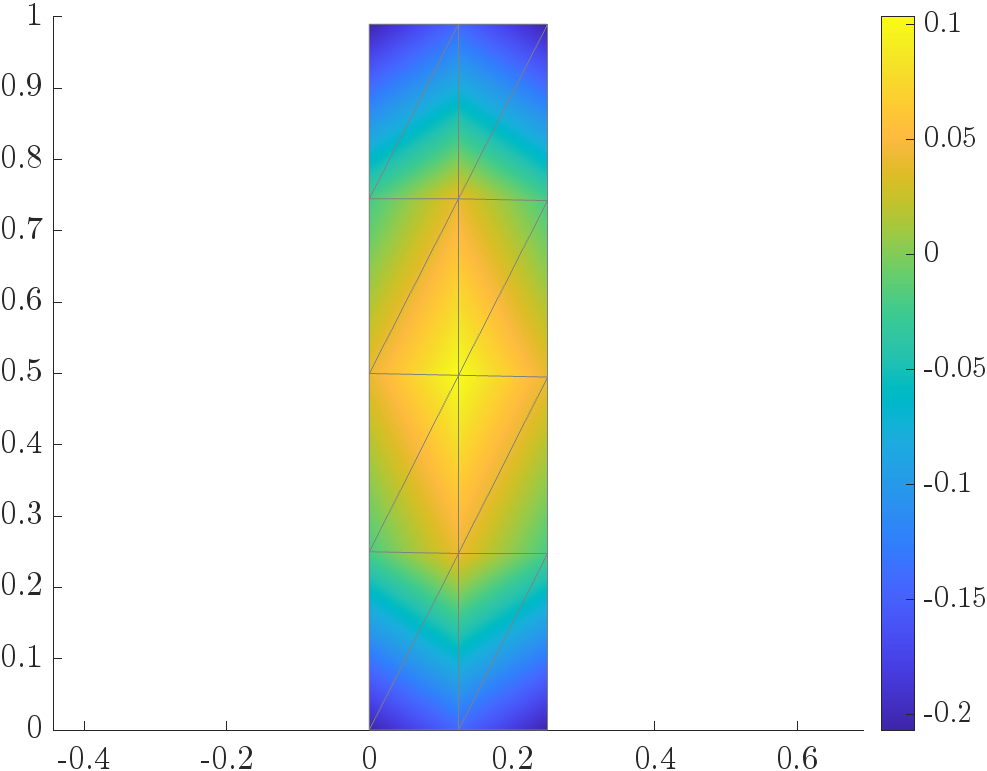}
  \caption{Geometry Rectangle (width $1/4$ and height $99/100$), denoted in short as RECT: solution $\psi'^0$ and final (adaptively refined) finite element mesh. Note we may construct from $\Omega^{2\mathrm{D}}$ associated extruded domains in $\RR^3$.}\label{fig:rect}
\end{figure}

\begin{figure}[H]
  \centering
  \includegraphics[width=.6\textwidth]{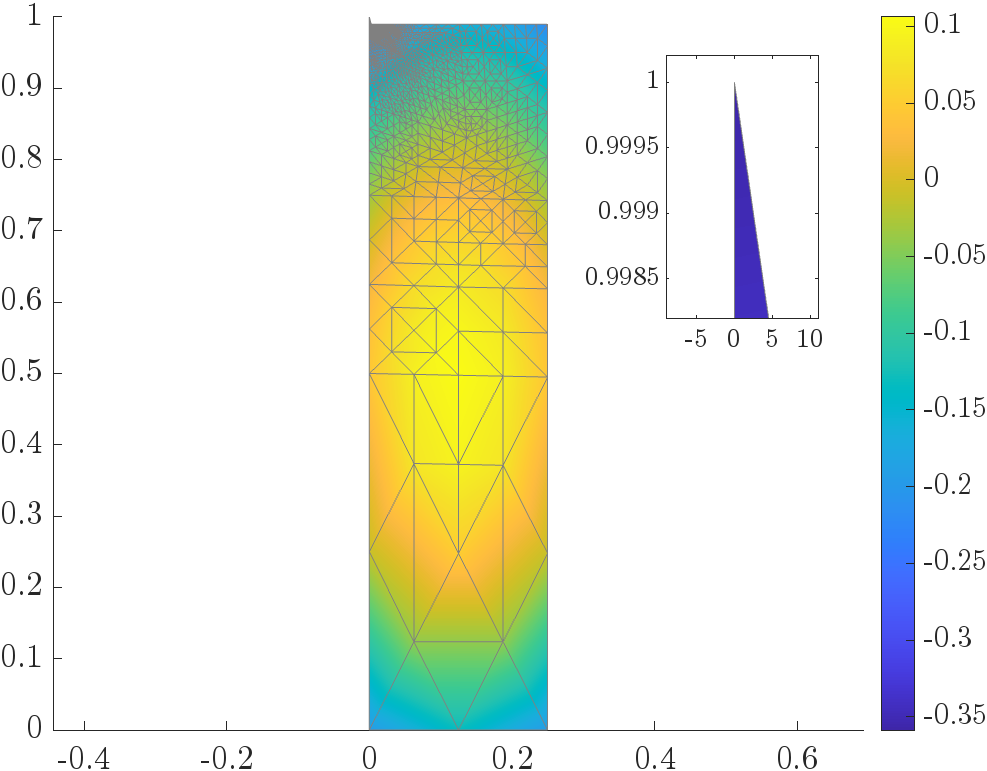}
  \caption{Geometry RECT with a small SART perturbation (union of RECT and SART), denoted in short as RECTSART: solution $\psi'^0$ and final (adaptively refined) finite element mesh. Note we may construct from $\Omega^{2\mathrm{D}}$ associated extruded domains in $\RR^3$.}\label{fig:rectsart}
\end{figure}

\begin{figure}[H]
  \centering
  \includegraphics[width=.6\textwidth]{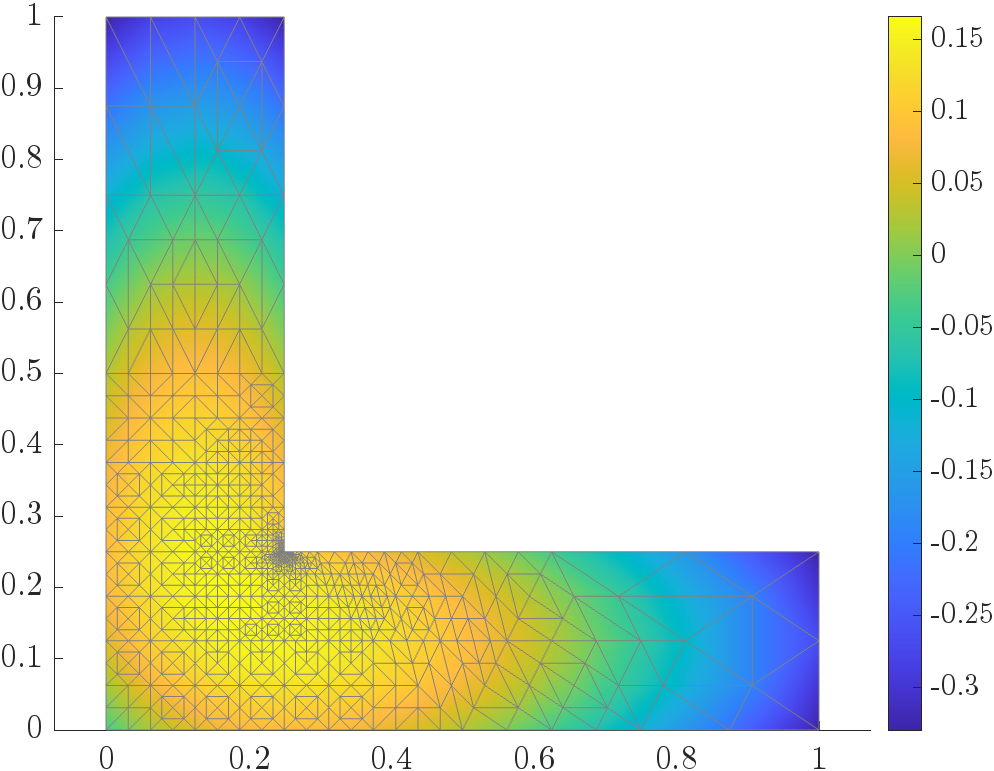}
  \caption{Geometry Bracket Equi-Thickness Arms, denoted in short as BETA: solution $\psi'^0$ and final (adaptively refined) finite element mesh. Note we may construct from $\Omega^{2\mathrm{D}}$ associated extruded domains in $\RR^3$.}\label{fig:beta}
\end{figure}

\begin{figure}[H]
  \centering
  \includegraphics[width=.6\textwidth]{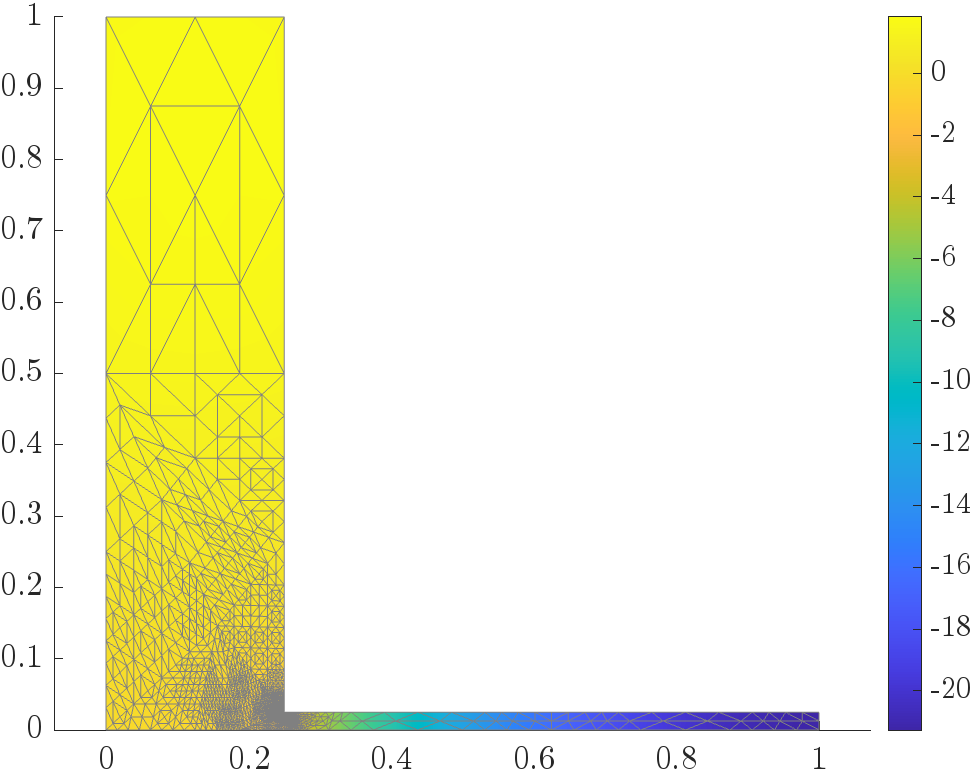}
  \caption{Geometry Bracket with Disparate Thickness Arms (ratio 10), denoted in short as BDTA: solution $\psi'^0$ and final (adaptively refined) finite element mesh. Note we may construct from $\Omega^{2\mathrm{D}}$ associated extruded domains in $\RR^3$.}\label{fig:bdta}
\end{figure}

\begin{figure}[H]
  \centering
  \begin{subfigure}[b]{0.49\textwidth}
    \includegraphics[width=\textwidth]{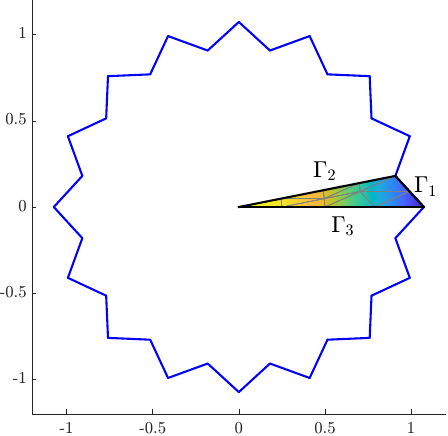}
    \caption{Full domain and half-tooth domain with boundaries.}
    \label{subfig:full_gear}
  \end{subfigure}
  \begin{subfigure}[b]{0.49\textwidth} 
    \includegraphics[width=\textwidth]{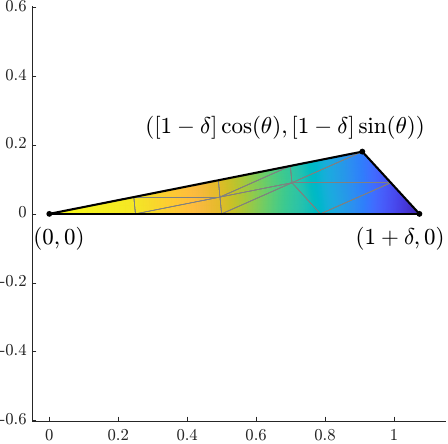}
    \caption{Half-tooth computational domain, vertices indicated with $\theta$ and $\delta$ .}
    \label{subfig:gear_param}
  \end{subfigure}
  \caption{Geometry Gear: full domain and extracted half-tooth ($q = 1.6$, $n = 32$). We present for the half-tooth the solution $\psi'^0_h$ and final (adaptively refined) finite element mesh.}
  \label{fig:gear}
\end{figure}

We shall first focus on the dependence of $\phi$ on $\Omega$; we set $\sigma=\kappa=1$. We summarize our results, in Table \ref{tab:phiGeos} and associated Figures \ref{fig:sart} -- \ref{fig:gear}, for a variety of geometric configurations: labels in Table \ref{tab:phiGeos} will be used to reference test domains and the acronyms will be expanded as a pneumonic in associated figures. We observe, consistent with the asymptotic theory of Proposition \ref{prop:phirt}, that $\phi$ for SART (Figure \ref{fig:sart}) is quite large (and indeed very close to the asymptotic result). We also observe in Figure \ref{fig:sart} the rapid decrease of $\psi'^0$ and associated large gradients of $\psi'^0$ into the corner; the latter gives rise to the large value of $\phi$. But we then can conclude from SARTF (Figure \ref{fig:sartf}) that the sharp angle is not required for large $\phi$ (less than 1\% $\phi$ difference compared to SART); note the fillet relieves the sharp corner but only modifies an extremely small portion of the domain. We further conclude from SARTC (Figure \ref{fig:tznrt}) that neither an acute angle nor high curvature is required for large $\phi$. We also observe that RECT (Figure \ref{fig:rect}) yields the analytical $\phi$ for a rectangle ($2/3$) with numerical error on the order of machine-precision, and furthermore that a small SART perturbation to RECT, RECTSART (Figure \ref{fig:rectsart}), does not lead to a significant increase in $\phi$.  Finally we note that the two bracket cases, BETA (Figure \ref{fig:beta}) and BDTA (Figure \ref{fig:bdta}), superficially somewhat similar, lead to very different values of $\phi$; demonstrating that a composite domain of small-$\phi$ subdomains can be associated with a large value of $\phi$. We also remark that in most of the cases, $\phi_h/\gamma$ is on the order of $10^{-1}$, meaning $\phi$ and $\gamma$ are correlated; however, we observe from configuration BDTA that $\phi_h$ can be disproportionately large compared to $\gamma$.

\input{Tables/phigeos}

These results which may appear to defy simple summary are, however, consistent with our conjecture as regards $\text{dist}$ \eqref{eq:distdef}. In fact, for the cases of Table \ref{tab:phiGeos}, $\dhaus$ suffices for distance. The results of Table \ref{tab:phiGeos} are also consistent with our two features. $\Feat{1}(\equiv\text{MWIA})$ suffices for all the triangle-related cases; for example it distinguishes between SART and RECTSART. However, for these triangle-related cases, $\Feat 2(\equiv\text{InRadius }\gamma)$ is not helpful. In contrast, for discrimination between BETA and BDTA, $\Feat 1$ is not relevant and instead $\Feat{2}$ is required. We can interpret the latter in terms of the thermal resistance concept (Section \ref{sec:lumped}): for BDTA, $\calL$ is not representative of the length scale associated with the conduction thermal resistance $\Lcond$, and must be corrected by $\phi$; in temporal terms, the small-thickness arm will equilibrate quickly to ambient  --- $T_\paavg \approx T_\infty$ --- meaning less heat transfer and hence a  \textit{larger} in amplitude $-\phi$ (negative) contribution to $\lambda_1$.




We now consider the gear shown in Figure \ref{subfig:full_gear} for a more comprehensive review of the effect of $\phi$ as well as further justification of $\text{dist}$ \eqref{eq:distdef}; note our extrinsic length scale $\dv{\ell}$ is the radius. We consider $ m_{\text{teeth}} \geq 4 $, even, and define the number of half teeth as $ n = 2m_{\text{teeth}} $. We then consider a sequence of geometries indexed by the number of half teeth, $ n $: $ \delta(n;q) \equiv \left(\theta(n)\right)^q $, $ \theta(n) \equiv \pi/n $; here $ \delta $ and $ \theta $ determine the halftooth dimensions in Figure \ref{subfig:gear_param}. \footnote{Our gear example is not realistic for $ q < 1 $ and $ n \rightarrow \infty $, not only in terms of the geometry, but also in the assumption that $ B $ is independent of the number of teeth: in practice, convection and radiation heat transfer coefficients would decrease as the exterior tooth angle closes. Our purpose is only to illustrate that our distance $\text{dist}$ and features are relevant to very rough geometries such as the gear for finite $ n $.} Note that by symmetry arguments, $ \phi $ for the full gear is given by $ \phi $ for the half tooth with homogeneous Neumann conditions on $ \Gamma_2 $ and $ \Gamma_3 $ (Figure \ref{subfig:full_gear}). The formulation for $ \psi'^{0} $ and $ \phi $ on the full gear requires only minor modification to treat the half tooth and associated homogeneous Neumann conditions: the same geometric factors in equation \eqref{eq:phiellp_0} yield solvable systems either on the full gear with Robin conditions or on the half tooth with mixed Robin-Neumann conditions. Of course the interest in the half tooth is computational expediency. 

We present tabulated computational results in Table \ref{tab:gear} for the cases of $(q,n)\in\{1.6,0.4\}\times\{32,256,2048\}$. For $q = 1.6 > 1$, our gear approaches a disk (of unit radius) in $\dhaus$, in domain measure $|\Omega^{\mathrm{2D}}|$ and in boundary measure $|\partial\Omega^\mathrm{2D}|$; we observe that $\phi$ for the gear approaches $\phi$ for the disk while $\text{dist}(gear,disk)$ approaches 0. For $q =0.4 < 1$,  our gear approaches a disk in $\dhaus$ and also in domain measure $|\Omega^\mathrm{2D}|$ but not in boundary measure $|\partial\Omega^\mathrm{2D}|$ (indicated by $\gamma$ progression) --- and hence not in $\text{dist}$; we observe that $\phi$ does not converge to $\phi$ for a disk. For these geometries, we observe that our $\Feat{2}$ is highly predictive of large $\phi$; a consequence of effective $\FF$ lower bound, especially as the number of teeth increases.

As anticipated, a large change in $\phi$ for two domains is reflected in a large distance $\text{dist}$ between the two domains. Note also that $\dhaus$ alone is not sufficient: $\phi$ diverges even as $\dhaus$ between successive gears tends to zero. We emphasize that as we increase the number of teeth, the geometry changes but also the geometric factors in the equation for $\phi$ \eqref{eq:phiellp_0} change.\footnote{We also note that $\phi/\gamma$, and hence our error bound $\Bdunkp$, slowly diverges — and hence so will our bounds for, say, the relative error in the eigenvalue. Indeed, it can be shown that, for $q < 1$, $\lambda_1$ converges to a value independent of $B$, observed computationally to be the Dirichlet eigenvalue for the disk.} Similarly to the BDTA shape class, we can interpret the geometry-dependence of $\phi$ in terms of the thermal resistance concept (Section \ref{sec:lumped}): for this gear shape class for low $q$, $\calL$ is not representative of the length scale associated with the conduction thermal resistance $\Lcond$, and must be corrected by $\phi$; in temporal terms, the small-thickness arm will equilibrate quickly to ambient  --- $T_\paavg \approx T_\infty$ --- meaning less heat transfer and hence a \textit{larger} in amplitude $-\phi$ (negative) contribution to $\lambda_1$. 
The gear can also provide insight into our features. We consider here the case $q < 1$. As expected from our discussion of Section \ref{sec:altstrat}, both features (\ref{eq:feat1def}) -- (\ref{eq:feat2def}) diverge as $n$ increases: the gears for larger $n$ are clearly separated from members of a dictionary which correspond to small $\phi$.

\input{Tables/gear}

\begin{figure}[H]
  \centering
  \includegraphics[width=.2\textwidth]{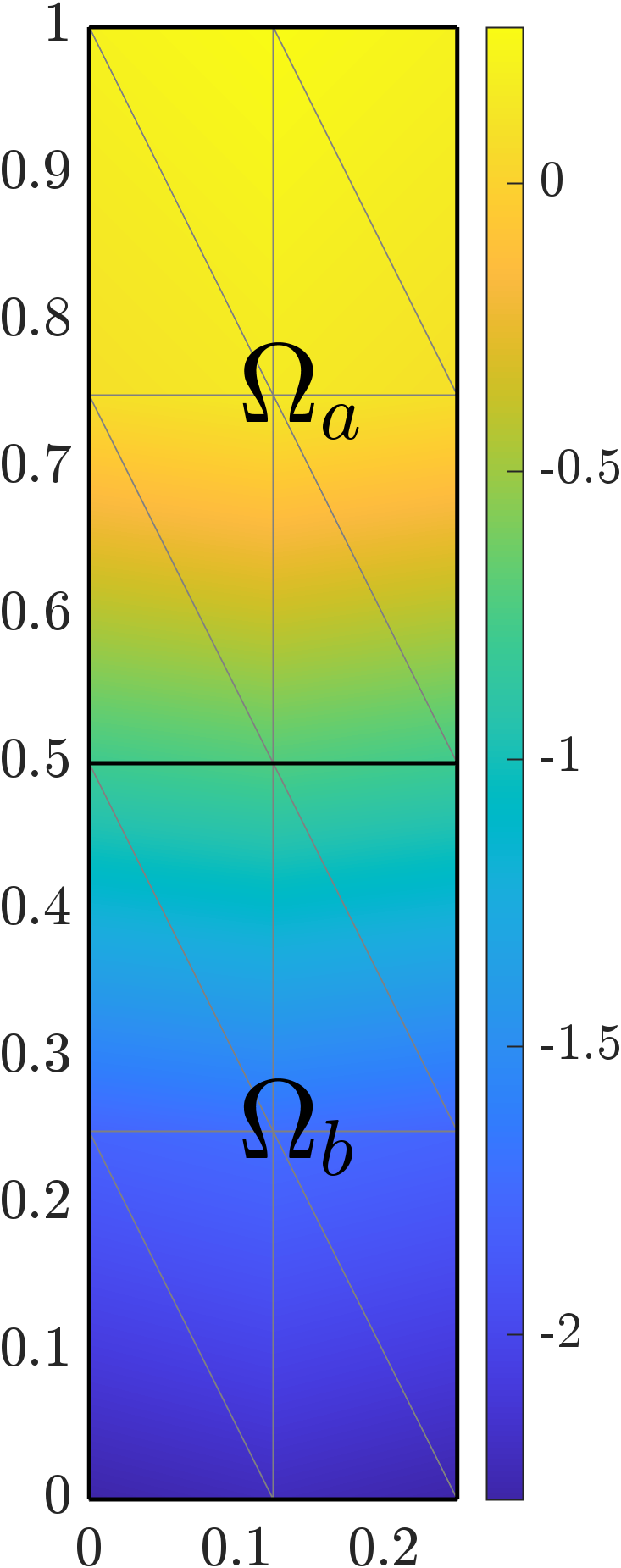}
  \caption{Geometry Rectangle with Horizontal Interface, denoted in short as RECTHI: solution $\psi'^0$ and final (adaptively refined) finite element mesh. The geometry is that of RECT with increased height of 1 and material interface at $y=1/2$.  Note we may construct from $\Omega^{2\mathrm{D}}$ associated extruded domains in $\RR^3$.}\label{fig:recthi}
\end{figure}

\begin{figure}[H]
  \centering
  \includegraphics[width=.5\textwidth]{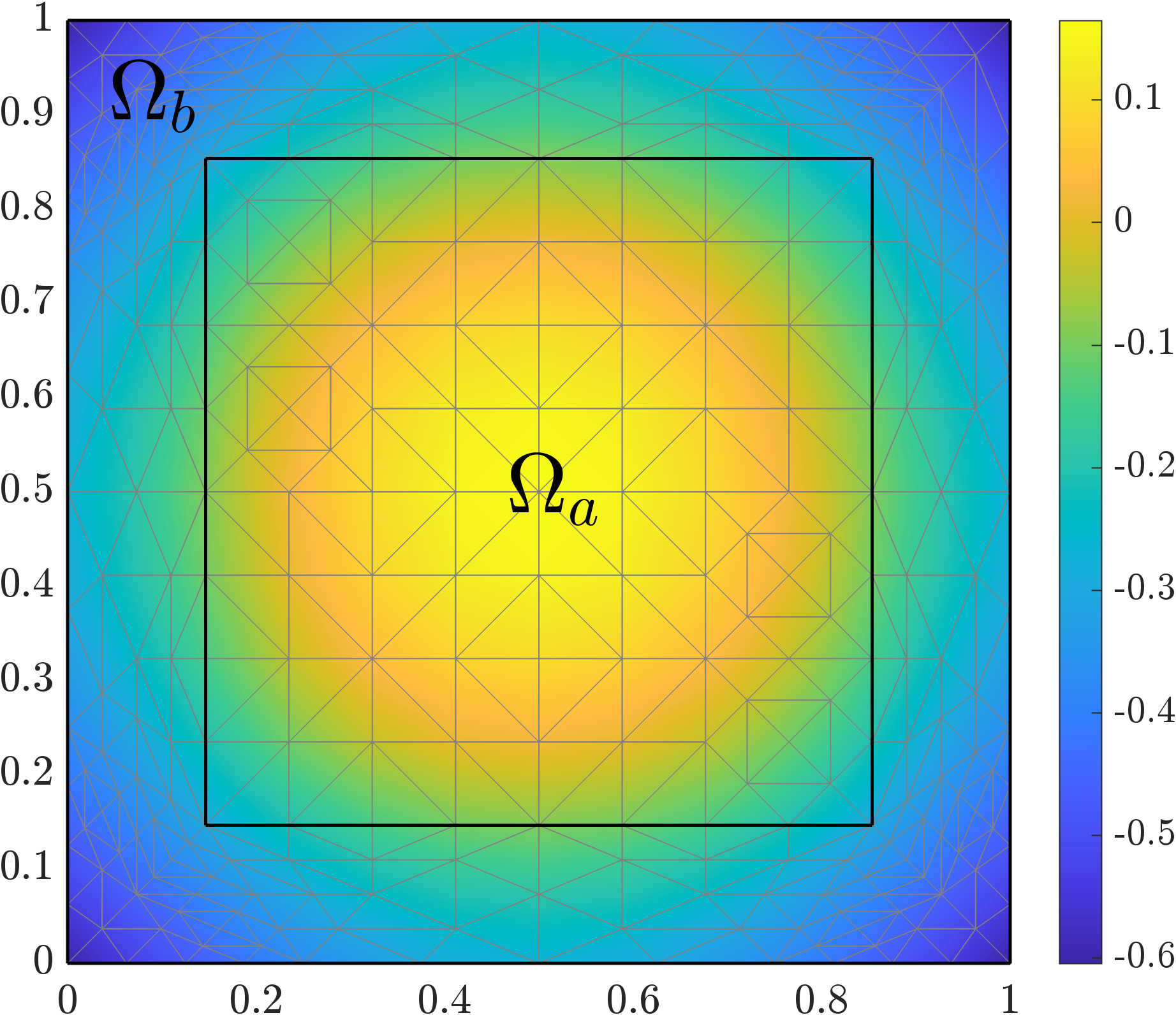}
  \caption{Geometry Equi Volume-Fraction Concentric Squares, denoted in short as EVFCS: solution $\psi'^0$ and final (adaptively refined) finite element mesh.  Note we may construct from $\Omega^{2\mathrm{D}}$ associated extruded domains in $\RR^3$.}\label{fig:evfcs}
\end{figure}

\begin{figure}[H]
  \centering
  \includegraphics[width=.5\textwidth]{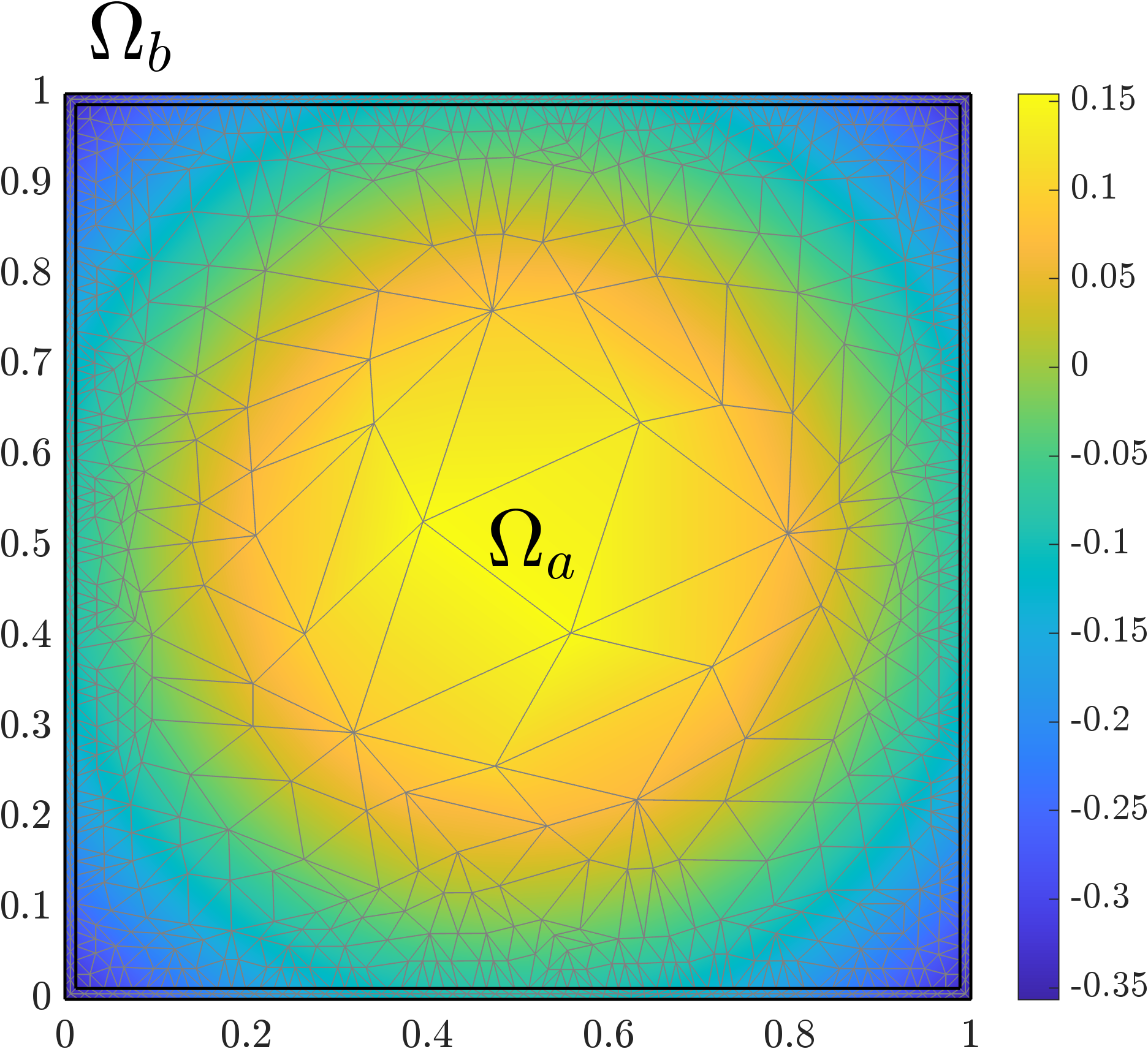}
  \caption{Geometry Disparate Volume-Fraction Concentric Squares with Light Film, denoted in short as DVFCSLF: solution $\psi'^0$ and final (adaptively refined) finite element mesh.  Note we may construct from $\Omega^{2\mathrm{D}}$ associated extruded domains in $\RR^3$.}\label{fig:dvfcslf}
\end{figure}

\begin{figure}[H]
  \centering
  \includegraphics[width=.5\textwidth]{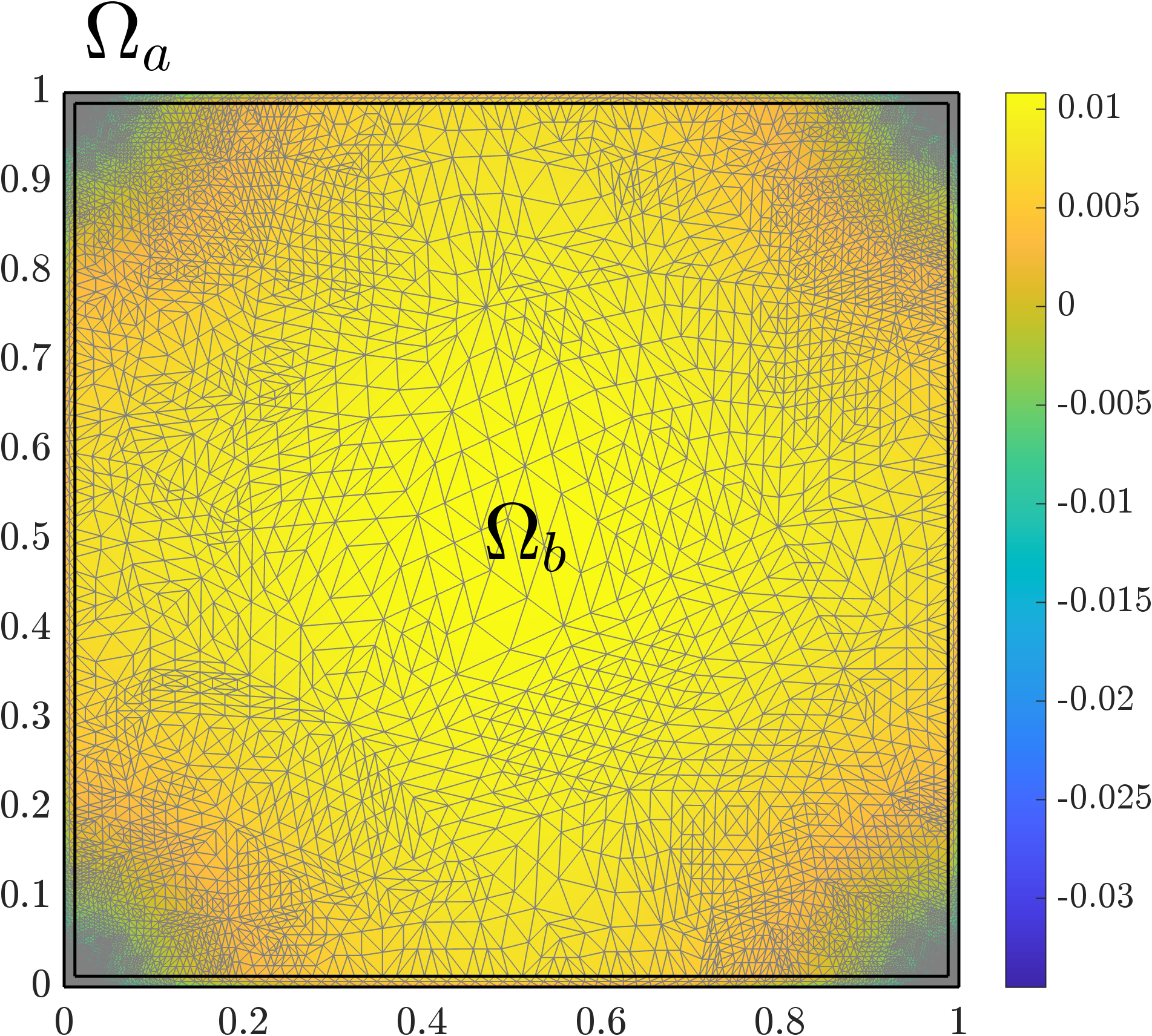}
  \caption{Geometry Disparate Volume-Fraction Concentric Squares with Heavy Film, denoted in short as DVFCSHF: solution $\psi'^0$ and final (adaptively refined) finite element mesh.  Note we may construct from $\Omega^{2\mathrm{D}}$ associated extruded domains in $\RR^3$.}\label{fig:dvfcshf}
\end{figure}

We now turn to a few examples in which the thermophysical properties are non-uniform, and in particular $\sigma \neq 1$ (we continue to take $\kappa = 1$). The relevant theoretical results are Proposition \ref{prop:phisigma} and Proposition \ref{prop:evfsigma}. In this section we consider a domain $\RR^2$ ($d=2$) comprised of two materials, $\sigma_a$ over $\Omega_a$ and $\sigma_b$ over $\Omega_b$, with $\sigma_a/\sigma_b = 1000$ and $|\Omega_a|/|\Omega_b| = r$ (to be varied). We consider four configurations: a horizontal interface rectangle RECTHI (Figure \ref{fig:recthi}), an equal volume-fraction concentric squares EVFCS (Figure \ref{fig:evfcs}), and two disparate volume-fraction concentric squares DVFCSLF and DVFCSHF (Figures \ref{fig:dvfcslf} -- \ref{fig:dvfcshf}) corresponding to light and heavy `film' configurations, respectively. We present the results for all 4 configurations in Table \ref{table:phibound}. Overall, we observe that the $\sigma$ contrast has a significant impact on $\phi$, and that the bounds of Proposition \ref{prop:phisigma} and \eqref{eq:muPW} are effective in some cases but not in others.  We observe that volume fraction has a significant impact on $\phi$ as demonstrated by the two concentric square examples EVFCS and DVFCSLF. (Note we do not consider small $\sigma$ contrast since the effect of non-uniformity, as predicted by Proposition \ref{prop:phisigma}, is quite small.) We also note that for the two $r=1$ cases RECTHI and EVFCS, Proposition \ref{prop:evfsigma} is applicable and we do indeed observe less than unity $\sigma$ variance (5th column of Table \ref{table:phibound}).

\input{Tables/twodomains}

We start with a rectangular domain (almost) equivalent to that of RECT shown in Figure \ref{fig:recthi}. The domain is split into two equal parts along a horizontal interface ($r = 1$). The upper subdomain $\Omega_a$ corresponds to $\sigma_a = 1.998$ and the lower subdomain $\Omega_b$ corresponds to $\sigma_b = 1.998\times 10^{-3}$. We observe that the bound for $\phi$ is within a factor 2 of $\phi_h$. We can provide a physical interpretation of this high-$\phi$ configuration: the lumped approximation assumes uniform decay to equilibrium over the entire domain; however, it is clear that, for this configuration, the two subdomains are coupled only along a short interface -- hence effectively decoupled -- with disparate $\sigma$-based time scales. Therefore, a large correction to $\Bdunk$ is required.

Next we consider EVFCS (shown in Figure \ref{fig:evfcs}) for which the two subdomains, an inner-square ($\Omega_a$) and border ($\Omega_b$), each have an area of $1/2$ ($r=1$). The outer region $\Omega_b$ corresponds to $\sigma_b = 1.998\times 10^{-3}$ and the inner region $\Omega_a$ corresponds to $\sigma_a = 1.998$.  We observe a decrease in the sharpness of the bound for $\phi$ compared to RECTHI. More importantly, we observe that $\phi_h$ is $\pcalO(1)$, as might be expected on physical grounds: the two subdomains interact in a coupled fashion, since the thermal flux from $\Omega_a$ \textit{must} pass through $\Omega_b$ prior to leaving the domain. Therefore, compared to RECTHI, the $\Bdunk$ correction is less severe.

For DVFCSLF (shown in Figure \ref{fig:dvfcslf}), we start with the same domain and topology as EVFCS, but we increase $r$ to 20, for the same $\sigma$ ratio. Due to this disparity in the volume-fraction, the $\sigma$ values are significantly altered to $\sigma_a=1.05$ and $\sigma_b=1.05\times 10^{-3}$. We do see an increase in sharpness of the bound for $\phi$ in this case, a result of dramatically reduced $\sigma$ variance. We also note that, for this configuration, we expect to see $\phi(\sigma)$ approach that of $\phi(1)=2/3$ as $r$ increases: $\sigma$ variance diminishes while $\gamma$ and $\mu$ are fixed. In effect, since
\begin{align}
    \dv{\frac{(\rho c)_a|\Omega_a|}{(\rho c)_b|\Omega_b|}} = \frac{\sigma_a|\Omega_a|}{\sigma_b|\Omega_b|} = 1000\cdot 20,
\end{align}
so material $b$ is effectively absent.

For DVFCSHF (shown in Figure \ref{fig:dvfcshf}), we start with the same domain and topology as DVFCSLF, but we exchange $\Omega_a$ and $\Omega_b$ so that now $r=1/20$. The $\sigma$ values are now closer to those from EVFCS, $\sigma_a=2.06$ and $\sigma_b=2.06\times 10^{-3}$, but due to the proportions of the subdomains,  we see a remarkable increase in $\sigma$ variance --- leading to large $\phi_\text{UB}$. The contrast with DVFCSLF demonstrates the importance of topology of the composite material with respect to the environment in determining $\phi(\sigma)$.

%% file: Tables/phigeos.tex
\begin{table}[ht]
\centering
\caption{Numerical Values of $\phi$ and Related Quantities for Select Domains.}
\label{tab:phiGeos}
\begin{tabular}{l|l|l|l} \multicolumn{1}{c|}{Shape} &\multicolumn{1}{c|} {$\phi_h$} & \multicolumn{1}{c|}{$\phi_h / \gamma$} & \multicolumn{1}{c}{$e_h^{\phi}$} \\ \hline \hline
SART     & $\sci{9.13}{ 0}$ & $\sci{5.01}{-1}$ & $\sci{1.64}{-14}$ \\ 
SARTF    & $\sci{9.07}{ 0}$ & $\sci{4.97}{-1}$ & $\sci{1.56}{-7 }$ \\ 
SARTC    & $\sci{9.06}{ 0}$ & $\sci{4.97}{-1}$ & $\sci{2.17}{-7 }$ \\ 
RECT     & $\sci{6.67}{-1}$ & $\sci{6.65}{-2}$ & $\sci{5.90}{-16}$ \\ 
RECTSART & $\sci{7.04}{-1}$ & $\sci{6.98}{-2}$ & $\sci{2.63}{-7 }$ \\ 
BETA     & $\sci{6.95}{-1}$ & $\sci{7.60}{-2}$ & $\sci{7.45}{-7 }$ \\ 
BDTA     & $\sci{7.55}{ 1}$ & $\sci{5.07}{ 0}$ & $\sci{8.90}{-4 }$
\end{tabular}
\end{table}

%% file: Tables/gear.tex
\begin{table}[H]
\centering
\caption{Numerical results for half-tooth geometries $(q,n)\in\{1.6,0.4\}\times\{32,256,2048\}$.}
\label{tab:gear}
\begin{tabular}{c||l|l|l||l|l|l}
 & \multicolumn{3}{c||}{$q=1.6$} & \multicolumn{3}{|c}{$q=0.4$}              \\ \hline
$n$        & \multicolumn{1}{|c}{$32$} & \multicolumn{1}{|c}{$256$} & \multicolumn{1}{|c||}{$2048$} & \multicolumn{1}{|c}{$32$} & \multicolumn{1}{|c}{$256$} & \multicolumn{1}{|c}{$2048$} \\ \hline
$\gamma$   & $\sci{2.53}{0} $  & $\sci{2.05}{0}  $ & $\sci{2.00}{0}  $  & $\sci{1.49}{ 1 }$ & $\sci{3.91}{ 1 }$ & $\sci{1.30}{ 2 }$ \\ \hline
$\phi_h$   & $\sci{7.94}{-1}$  & $\sci{5.23}{-1 }$ & $\sci{5.02}{-1 }$  & $\sci{3.00}{ 1 }$ & $\sci{1.94}{ 2 }$ & $\sci{2.12}{ 3 }$ \\ \hline
$e^\phi_h$ & $\sci{1.22}{-14}$ & $\sci{4.30}{-13}$ & $\sci{2.92}{-11}$  & $\sci{1.36}{-11}$ & $\sci{1.84}{-10}$ & $\sci{9.66}{-8 }$\\ \hline
$\FF$   & $\sci{5.94}{-1}$  & $\sci{5.18}{-1 }$ & $\sci{5.02}{-1 }$  & $\sci{2.00}{ 0 }$ & $\sci{7.18}{ 1 }$ & $\sci{1.41}{ 3 }$ 
\end{tabular}
\end{table}
%$\phi_h / \gamma$ & 0.314 & 0.256 & 0.250 & 2.02 & 4.97 & 16.3 \\ \hline

%phis =
%
%     7.944082548203217e-01
%     5.234186355550829e-01
%     5.019290356087089e-01
%
%K>> FF
%
%FF =
%
%     5.939701311050799e-01
%     5.179905012460405e-01
%     5.017393634479493e-01
%
%
%
%FF =
%
%     2.004304044685498e+00
%     7.179984352595658e+01
%     1.409248162287651e+03
%
%K>> close all
%K>> phis
%
%phis =
%
%     3.002352104318972e+01
%     1.939486202522156e+02
%     2.122485947203602e+03

%% file: Tables/twodomains.tex
\begin{table}[H]
\begin{center}
\caption{Two-material $\phi_h$ results with bounds and associated quantities. Note: all geometries correspond to $\mu=\pi^2$.}  \label{table:phibound}
\begin{tabular}{l|c|c|c|c|c|l}
\multicolumn{1}{c|}{Shape} & $r$ & $\mu_h$ & $\dashint_\Omega \left(\sigma-1\right)^2$ & $\phi_{\mathrm{UB}}$ & $\phi_h(\sigma)$  & \multicolumn{1}{|c}{$e_h^\phi$}\\
\hline\hline
RECTHI & 1 &  9.87 & 0.996 & 15.9 & 8.97 & $\sci{1.66}{-14}$\\
EVFCS & 1 &  9.87 & 0.996 & 4.36 & 1.58 &         $\sci{7.17}{-6}$\\
DVFCSLF & 20 &  9.87 & $0.0499$ & 1.21 & 0.732 &  $\sci{1.77}{-9}$ \\
DVFCSHF & 1/20 &  9.87 & $19.2$ & 40.9 & 0.0181 & $\sci{6.07}{-10}$
\end{tabular}
\end{center}
\end{table}

%% file: err.tex
\section{Error Estimation}\label{sec:qoierr}

In this section we establish QoI bounds for the solution to the initial value problem \eqref{eq:ivpstrong_1}~--~\eqref{eq:ivpstrong_3} and errors associated with QoI approximations introduced in Definition \ref{def:lumpedapprox}: $\uavgL{1}$, $\uavgLP$, and $\uDeltaL$. We shall consider primarily two domains: SART-1 and SART-2. Both are instances of our Small-Angle Right Triangle introduced in Proposition \ref{prop:phirt} ; SART-1 corresponds to a relatively small angle (equivalent to SART in Section \ref{sec:comp_results}); SART-2 corresponds to an even smaller angle ($W = 1/16$). These domains serve in the numerical examples presented throughout this section.

Also in this section it shall be convenient to introduce a ``slow” time scale:
\begin{definition}[Slow-Time Representation]\label{def:slowtime} Let us define
\begin{align}
T \equiv B\gamma t,\ T_{\text{final}} \equiv B\gamma t_\text{final}
\end{align}
and associated field and quantities of interest,
\begin{align}
U(x,T) \equiv u(x,t = T/(B\gamma)),
\end{align}
and
\begin{align}
U_\avg(T) &\equiv u_\avg(t = T/(B\gamma)),\\
U_\paavg(T) &\equiv u_\paavg(t = T/(B\gamma)),\\
U_\Delta(T) &\equiv u_\Delta(t = T/(B\gamma)),
\end{align}
respectively. The arguments $(x,T)$ (respectively, $T$) are an abbreviated form for $(x,T;B)$ or even $(x,T;B;\sigma,\kappa)$ (respectively,  $(T;B)$ or even $(T;B;\sigma,\kappa)$).

In a similar fashion,
\begin{align}
\UavgL{1}(T) &\equiv \uavgL{1}(t = T/(B\gamma)),\label{eqn:defUavgL1}\\
\UavgLP(T) &\equiv \uavgLP(t = T/(B\gamma)),\text{ and}\label{eqn:defUavgL2P}\\
\UDeltaL(T) &\equiv \uDeltaL(t = T/(B\gamma))\label{eqn:defUDeltaL}
\end{align}
are the slow-time representations of our lumped approximations of Definition \ref{def:lumpedapprox}. Again, $T$ is an abbreviated argument list.
\end{definition}
Note that $\UavgL{1} (T) = \exp(-T)$, which demonstrates the convenience of the slow time representation — in which time is effectively normalized with respect to the small-$B$ time constant.

Finally, we make a remark on nomenclature: $T = B\gamma t = \dv{h |\partial\Omega| t/(\rho c |\Omega|)}$ is our nondimensional slow time; in contrast, $\dv{T}$ is the dimensional temperature, with nondimensional counterpart $u \equiv \dv{(T - T_\infty)/(T_{\text{i}} - T_\infty)}$. We shall not invoke dimensional variables in the remainder of this section.

%% file: uavg1.tex
\subsection{First-Order Lumped Approximation: Domain-Average}\label{sec:uavg1}

Before we proceed to the main domain-average analysis, we first consider the positivity of the relative domain-boundary difference $U_\Delta(T)$ for all $T> 0$:
\begin{lemma}[Strictly Positive $U_\Delta$]\label{lem:UDeltaL} The relative difference between the domain-mean and the boundary-mean temperature satisfies
\begin{align}
U_\Delta(T;B) > 0,\ \forall\,  T,B >0.\label{eq:UDeltaL}
\end{align}
For $T = 0$ or $B = 0$, we have equality of the domain and boundary averages such that
\begin{align}
U_\Delta(0;B) &= 0,\label{eq:UDeltaLT0}\\
U_\Delta(T;0) &= 0.\label{eq:UDeltaLB0}
\end{align}
\end{lemma}
We provide a sketch of the proof for Lemma \ref{lem:UDeltaL} in Appendix \ref{App:UDeltaL}.  Following this result, we have:
\begin{proposition}[Global Domain-Average Lower Bound]\label{prop:UavgBound} The first-order lumped approximation $\UavgL{1}$ is a lower bound for the domain-average $U_\avg(T)$ for all $T, B\in \Rzp$. More specifically, for $B>0$ and $T>0$, we have a strict inequality
\begin{align}
U_\avg(T)>\UavgL{1}(T).
\end{align}
\end{proposition}
We note that this lower bound result has a simple physical interpretation: the (first-order) lumped approximation neglects the conduction resistance, $\dv{\calR_{\text{avg}}}$ of \eqref{eq:ravgdef2} in Section \ref{sec:lumped}, and hence the associated heat transfer rate is greater than the actual heat transfer rate, and, as a result, the lumped temperature approaches equilibrium faster than the actual temperature, and thus --- given our nondimensionalization --- the nondimensional lumped temperature will always be less than the actual temperature.  We provide a sketch of the proof for Proposition \ref{prop:UavgBound} in Appendix \ref{App:UavgBound}.

Now We define our errors associated with the first-order lumped approximation $\UavgL{1}$:
\begin{definition}[First-Order Lumped Domain-Average: Absolute Error]
We define
\begin{align}
\EavgL{1}(\Tf) \equiv \max_{T \in [0,T_\text{final}]} U_\avg(T) - \UavgL{1}(T),\label{eq:e1avg}
\end{align}
As always, $T$ is an abbreviated argument list.
\end{definition}

We can then prove, in similar fashion to \cite{ODD}, but for general geometry,
\begin{proposition}[Asymptotic Error Estimates for $\UavgL{1}$]\label{Prop:N1}
The absolute error, for any final time, satisfies
\begin{align}
\EavgL{1}(\Tf) \le \frac{\phi B}{\gamma\,e} + \pcalO(B^2)\text{ as }B \rightarrow 0\label{eq:easymp}\,.
\end{align}
We note the leading error term
\begin{align}
\EavgLa{1}(B) \equiv \underbrace{\frac{\phi B}{\gamma}}_{\Bdunkp}e^{-1}\,\label{eq:eavgla}
\end{align}
is proportional to the corrected intrinsic Biot number, $\Bdunkp$. We also note that
\begin{align}
T_\text{max} = 1 + \dfrac{1}{2} \underbrace{\frac{\phi B}{\gamma}}_{\Bdunkp} +\, \pcalO(B^2)\text{ as }B \rightarrow 0\label{eq:Tmax}
\end{align}
is the (slow) time $T$ at which $U_\avg(T) - \UavgL{1}(T)$ attains the maximum.
For both expansions about $B=0$, the high-order terms encapsulated in $\pcalO(B^2)$ are independent of $\Tf$. 
\end{proposition}
We note that in practice slow times $T$ of order unity are of greatest interest.  We provide in Appendix \ref{App:PropN1} a sketch of the full proof of Proposition \ref{Prop:N1}.  We also note that the relative error associated with $\UavgL{1}$ based on Proposition \ref{Prop:N1} can be defined (as well as its associated error bound).\footnote{ We can define $\EavgLr{1}(T) \equiv \left(U_\avg(T) - \UavgL{1}(T)\right)/\UavgL{1}(T)$ as the relative error associated with the first-order lumped approximation $\UavgL{1}$. We then define the relative error bound $\EavgLra{1}(T) \equiv \Bdunkp T$ such that $\EavgLr{1}(T) \le \EavgLra{1} + \pcalO(B^2)\text{ as }B \rightarrow 0$.
}

Now we consider non-asymptotic error bounds for $\UavgL{1}$:
\begin{proposition}[Non-Asymptotic Error Bounds for $\UavgL1$]\label{Prop:N2}
For all positive $B$, the absolute error satisfies
\begin{align}
\EavgL{1}(\Tf) \,\le\, \frac{1}{2}\sqrt{\frac{\phi B}{\gamma}}\,,\ \forall\, \Tf \in [0,\infty).
\end{align}
We note that for this non-asymptotic result, the error bound
\begin{align}
\EavgLUB{1} \equiv \frac{1}{2} \biggl(\underbrace{\frac{\phi B}{\gamma}}_{\Bdunkp}\biggr)^{1/2}\;
\end{align}
also depends on the corrected intrinsic Biot number $\Bdunkp$.
\end{proposition}
We note from Proposition \ref{Prop:N2} and Proposition \ref{prop:UavgBound} respectively that the domain average can be bounded by $\UavgL{1}$ and $\UavgL{1}+\EavgLUB{1}$. We provide a sketch of the proof for Proposition \ref{Prop:N2} in Appendix \ref{App:PropN2}. Similarly to Proposition \ref{Prop:N1}, we can define a relative error bound for $\UavgL{1}$.\footnote{
We can define $\EavgLrUB{1}(T) \equiv \frac{1}{2} \sqrt{\Bdunkp}\exp(T)$, where $\EavgLr{1}(T) \le \EavgLrUB{1}(T)$.
}

The non-asymptotic bound of Proposition \ref{Prop:N2} eliminates any uncertainty in the error estimate. However, it is clear from Proposition \ref{Prop:N1} and the result for $\EavgLa{1}$ that the non-asymptotic bound $\EavgLUB{1}$ is not sharp as $B \rightarrow 0$. On the other hand, for large $B$, the non-asymptotic bound $E_\avg^{1\,\text{UB}}$ will often be irrelevant — greater than unity. Note given the maximum principle, the error $E_\avg^{1}$ must always be less than one. Ironically, the non-asymptotic bound $E_\avg^{1\,\text{UB}}$ may be most relevant for very small $B$, confirming without any doubt that the error is less than some acceptable threshold.

We show in Table \ref{tab:NumRes1a} (respectively, Table \ref{tab:NumRes2a}) the absolute error and associated asymptotic error estimates and non-asymptotic error bounds for the first-order lumped approximation of $U_{\text{avg}}$ for domain SART-1 (respectively, SART-2). Note when $\EavgLUB{1} > 1$, no useful information is obtained and we hence indicate — in the tables. Note for SART-1, $\phi = 9.136$, and for SART-2, $\phi = 161.1$. The results of these tables are obtained by finite difference in time (BDF1 for the first time-step and BDF2 for subsequent time-steps) and finite element in space ($\mathbb{P}_2$) discretizations of the Robin heat equation. We integrate to $T_\text{final}=2$ and present $E^1_\avg(2)$. For any given result, we consider a sequence of uniform refinements in space and time until an error estimator (the $H^1$ norm of the discretization error in $U$ at $T=T_\text{final}=2$) is suitably small compared to the predicted error $\EavgL{1}$.

We make the following observations:
\begin{enumerate}
\item We confirm, in both Table \ref{tab:NumRes1a} and Table \ref{tab:NumRes2a}, that $\EavgL{1}$ indeed tends to zero linearly in B: the first-order approximation is indeed first-order. 
\item We confirm — by comparison of Table \ref{tab:NumRes1a} and Table \ref{tab:NumRes2a} — that, for a fixed $B$, or $\Bdunk$, the error increases substantially with increased $\phi$. 
\item We confirm that $\EavgLUB{1}$ is indeed an upper bound, albeit --- and as expected --- not sharp as $B \rightarrow 0$.
\item We observe that $\EavgLa{1}/\EavgL{1} \rightarrow 1$ and from above as $B \rightarrow 0$: the asymptotic error estimate is a sharp upper bound for small $B$. It is crucial to note that, absent the correct value of $\phi$, $\EavgLa{1}/\EavgL{1}$ would tend to a constant less than unity — underestimation of the error — or a constant greater than unity — overestimation of the error — as $B$ tends to zero. 
\end{enumerate}
We first note that the numerical results are illustrations of the provided mathematical proofs with regard to these observations. We also note that if $\EavgLa{1}$ is constructed in terms of an upper bound for $\phi$, for example as in Section \ref{sec:phiprops} Proposition \ref{prop:phisigma}, then $\EavgLa{1}$ will still provide an upper bound for $\EavgL{1}$ as $B\rightarrow 0$, however $\EavgLa{1}/\EavgL{1}$ will approach a constant greater than unity. We have not provided a rigorous proof for this statement, however.

\input{Tables/sart1a}

\input{Tables/sart2a}

The approximation $\UavgL{1}(T;B)$ is agnostic to the details of material properties, but the approximation is not restricted to problems characterized by uniform material properties. Note, however that although the evaluation of the approximation is independent of $\sigma$ and $\kappa$ --- we require only $B$ and geometric information ---  the error estimate must be defined relative to the (possibly) non-uniform problem data; the latter appear in the equation for $\phi$.

%% file: Tables/sart1a.tex
\begin{table}[H]
\centering
\caption{Absolute error and associated asymptotic error estimates and non-asymptotic error bounds for the first-order lumped approximation of $U_\avg(T)$, $T\in[0,2]$, for domain SART-1.}
\label{tab:NumRes1a}
\begin{tabular}{l|c|c|c|c|c|c}
\multicolumn{1}{c|}{$B$} & $\Bdunk\equiv B/\gamma$ & $\Bdunkp \equiv \phi B / \gamma$ & $E^1_\text{avg}(2)$ & $E_\text{avg}^\text{1 asymp}$ & $E_\text{avg}^\text{1 asymp} / E^1_\text{avg}(2)$ & $E_\text{avg}^\text{1 UB}$ \\ 
\hline\hline
$\sci{1}{-3}$ & $\sci{5.48}{-5}$ & $\sci{5.01}{-4}$ & $\sci{1.84}{-4}$ & $\sci{1.84}{-4}$ & $1.00$ & $\sci{1.12}{-2}$ \\ 
$\sci{2}{-3}$ & $\sci{1.10}{-4}$ & $\sci{1.00}{-3}$ & $\sci{3.67}{-4}$ & $\sci{3.68}{-4}$ & $1.00$ & $\sci{1.58}{-2}$ \\ 
$\sci{5}{-3}$ & $\sci{2.74}{-4}$ & $\sci{2.50}{-3}$ & $\sci{9.10}{-4}$ & $\sci{9.21}{-4}$ & $1.01$ & $\sci{2.50}{-2}$ \\ 
$\sci{1}{-2}$ & $\sci{5.48}{-4}$ & $\sci{5.01}{-3}$ & $\sci{1.80}{-3}$ & $\sci{1.84}{-3}$ & $1.03$ & $\sci{3.54}{-2}$ \\ 
$\sci{2}{-2}$ & $\sci{1.10}{-3}$ & $\sci{1.00}{-2}$ & $\sci{3.51}{-3}$ & $\sci{3.68}{-3}$ & $1.05$ & $\sci{5.00}{-2}$ \\ 
$\sci{5}{-2}$ & $\sci{2.74}{-3}$ & $\sci{2.50}{-2}$ & $\sci{8.17}{-3}$ & $\sci{9.21}{-3}$ & $1.13$ & $\sci{7.91}{-2}$ \\ 
$\sci{1}{-1}$ & $\sci{5.48}{-3}$ & $\sci{5.01}{-2}$ & $\sci{1.47}{-2}$ & $\sci{1.84}{-2}$ & $1.26$ & $\sci{1.12}{-1}$ \\ 
$\sci{2}{-1}$ & $\sci{1.10}{-2}$ & $\sci{1.00}{-1}$ & $\sci{2.44}{-2}$ & $\sci{3.68}{-2}$ & $1.51$ & $\sci{1.58}{-1}$ \\ 
$\sci{5}{-1}$ & $\sci{2.74}{-2}$ & $\sci{2.50}{-1}$ & $\sci{4.10}{-2}$ & $\sci{9.21}{-2}$ & $2.25$ & $\sci{2.50}{-1}$ \\ 
$\sci{1}{0}$  & $\sci{5.48}{-2}$ & $\sci{5.01}{-1}$ & $\sci{5.55}{-2}$ & $\sci{1.84}{-1}$ & $3.32$ & $\sci{3.54}{-1}$ \\ 
\end{tabular}
\end{table}

%% file: Tables/sart2a.tex
\begin{table}[H]
\centering
\caption{Absolute error and associated asymptotic error estimates and non-asymptotic error bounds for the first-order lumped approximation of $U_\avg(T)$, $T\in[0,2]$, for domain SART-2.}
\label{tab:NumRes2a}
\begin{tabular}{l|c|l|c|c|c|c}
\multicolumn{1}{c|}{$B$} & $\Bdunk\equiv B/\gamma$ & $\Bdunkp\equiv \phi B/\gamma $ & $E^1_\text{avg}(2)$ & $E_\text{avg}^\text{1 asymp}$ & $E_\text{avg}^\text{1 asymp} / E^1_\text{avg}(2)$ & $E_\text{avg}^\text{1 UB}$ \\ 
\hline\hline
$\sci{1}{-3}$ & $\sci{1.51}{-5}$ & $\quad\sci{2.44}{-3}$ & $\sci{8.89}{-4}$ & $\sci{8.97}{-4}$ & $1.01$ & $\sci{2.47}{-2}$\\ 
$\sci{2}{-3}$ & $\sci{3.03}{-5}$ & $\quad\sci{4.88}{-3}$ & $\sci{1.76}{-3}$ & $\sci{1.79}{-3}$ & $1.02$ & $\sci{3.49}{-2}$\\ 
$\sci{5}{-3}$ & $\sci{7.57}{-5}$ & $\quad\sci{1.22}{-2}$ & $\sci{4.27}{-3}$ & $\sci{4.49}{-3}$ & $1.05$ & $\sci{5.52}{-2}$\\ 
$\sci{1}{-2}$ & $\sci{1.51}{-4}$ & $\quad\sci{2.44}{-2}$ & $\sci{8.14}{-3}$ & $\sci{8.97}{-3}$ & $1.10$ & $\sci{7.81}{-2}$\\ 
$\sci{2}{-2}$ & $\sci{3.03}{-4}$ & $\quad\sci{4.88}{-2}$ & $\sci{1.49}{-2}$ & $\sci{1.79}{-2}$ & $1.21$ & $\sci{1.10}{-1}$\\ 
$\sci{5}{-2}$ & $\sci{7.57}{-4}$ & $\quad\sci{1.22}{-1}$ & $\sci{2.93}{-2}$ & $\sci{4.49}{-2}$ & $1.53$ & $\sci{1.75}{-1}$\\ 
$\sci{1}{-1}$ & $\sci{1.51}{-3}$ & $\quad\sci{2.44}{-1}$ & $\sci{4.31}{-2}$ & $\sci{8.97}{-2}$ & $2.08$ & $\sci{2.47}{-1}$\\ 
$\sci{2}{-1}$ & $\sci{3.03}{-3}$ & $\quad\sci{4.88}{-1}$ & $\sci{5.59}{-2}$ & $\sci{1.79}{-1}$ & $3.21$ & $\sci{3.49}{-1}$\\ 
$\sci{5}{-1}$ & $\sci{7.57}{-3}$ & $\quad\sci{1.22}{0} $ & $\sci{6.83}{-2}$ & $\sci{4.49}{-1}$ & $6.57$ & $\sci{5.52}{-1}$\\ 
$\sci{1}{0}$  & $\sci{1.51}{-2}$ & $\quad\sci{2.44}{0}$  & $\sci{7.46}{-2}$ & $\sci{8.97}{-1}$ & $12.0$ & $\sci{7.81}{-1}$\\ 
\end{tabular}
\end{table}

%% file: uavg2.tex
\subsection{Second-Order Lumped Approximation: Domain-Average}\label{sec:uavg2}

We now consider the second-order lumped approximation, $\UavgLP$ of \eqref{eqn:defUavgL2P}, which is based on a second-order Pad\'e-approximant \eqref{eq:lamapp2p} of $\lambda$. Note the Pad\'e choice is preferred over the Taylor series: in the former, unlike in the latter, $\lambda$ is bounded as $B\rightarrow\infty$, with an attendant improvement in accuracy for large $B$. Before proceeding, we introduce the third-order expansion for $\lambda$:
\begin{lemma}[First Eigenvalue Asymptotic Expansion Revisited]\label{lem:OB3} The first eigenvalue admits expansion
\begin{align}
\lambda = \gamma B - \phi B^2 + (\chi -\gamma\,\Upsilon) B^3 + \pcalO(B^4)\text{ as }B\rightarrow 0.\label{eq:lambdaO3}
\end{align}
We recall from \eqref{eq:phidef1} that $\phi\equiv\int_\Omega\kappa\, \nabla\psipz\cdot\nabla\psipz$, and from \eqref{eq:Upsilondef} that $\Upsilon\equiv\int_\Omega \sigma\,\psipz\psipz$; we now further introduce  $\chi\in \Rzp$ given by
\begin{align}
\chi &\equiv \int_{\partial\Omega}\psipz\psipz.\label{eqn:chidef}
\end{align}
we recall the definition $\psipz \equiv d\psi_1/dB\text{ evaluated at }B=0$.
\end{lemma}
We provide in Appendix \ref{App:OB3} a sketch of the proof of Proposition \ref{lem:OB3}.

\begin{proposition}[Scale-Invariance of $\gamma\chi$ and $\gamma^2\Upsilon$] \label{prop:chiscale}
 We are given $\Omega_1 \subset \RR^d$ and a translation-rotation-uniform dilation map of $\Omega_1$, ${\Omega}_2 \equiv \{ \MM(x) \, |\, x \in \Omega_1 \,\}$; here $\MM(x) = \TT + \QQ\, x + \alpha\,\II\,x$, where $\alpha$ is a positive real scalar, $\TT\in\mathbb{R}^d$ is a translation vector, $\QQ \in \RR^{d\times d}$ is a rotation matrix, and $\II \in \RR^{d\times d}$ is the identity. For given $\sigma \equiv \sigma_1, \kappa \equiv \kappa_1$, $\Omega \equiv \Omega_1$ we evaluate $\chi_1 \equiv \chi(\Omega_1)$, $\Upsilon_1 \equiv \Upsilon(\Omega_1)$, and $\gamma_1 \equiv |\Omega_1| / |\partial\Omega_1|$; for given $\sigma \equiv \sigma_2 \equiv \sigma_1 \circ \MM^{-1}, \kappa \equiv \kappa_2 \equiv \kappa_1 \circ \MM^{-1}$, $\Omega \equiv \Omega_2$ we evaluate $\chi_2 \equiv \chi(\Omega_2)$, $\Upsilon_2\equiv\Upsilon(\Omega_2)$, and $\gamma_2\equiv|\Omega_2|/|\partial\Omega_2|$. Then, for any $\alpha$, $\TT$, and $\QQ$, $\gamma_1\,\chi_1 = \gamma_2\,\chi_2$ and $\gamma_1^2\,\Upsilon_1=\gamma_2^2\,\Upsilon_2$.
\end{proposition}
We provide in Appendix \ref{App:sketch_chiscale} a sketch of the proof of Proposition \ref{prop:chiscale}.

Now we introduce the tensorization result for $\chi$ and $\Upsilon$:
\begin{proposition}[$\chi$ and $\Upsilon$for Tensorized Domains]\label{prop:chiupsilontendom} We introduce, for a given two-dimensional domain $\Omega^{\text{2D}} \subset \RR^2$, the tensorized (or extruded) three-dimensional domain $\Omega \subset \RR^3 \equiv \Omega^{\text{2D}}\times (0,L_3)$.
Then
\begin{align}
\Upsilon(1,1;\Omega) &= \Upsilon_{\text{2D}}(1,1;\Omega^{\text{2D}}) + \Upsilon_{\text{1D}}(1,1;\Omega^\text{1D}) \; , \label{eq:upsilontendom1}\\
\chi(1,1;\Omega) &= \chi_{\text{2D}}(1,1;\Omega^{\text{2D}}) + \chi_{\text{1D}}(1,1;\Omega^\text{1D}) \nonumber\\&+ \gamma_\text{1D} \Upsilon_\text{2D}(1,1;\Omega^\text{2D}) + \gamma_\text{2D}\Upsilon_\text{1D}(1,1;\Omega^\text{1D})\label{eq:chitendom1}
\end{align}
where $\chi_{\text{2D}}$, $\Upsilon_\text{2D}$, and $\gamma_\text{2D}$ respectively refer to the values of $\chi$, $\Upsilon$, and $\gamma$ for the $\text{2D}$ domain, $\Omega_\text{2D}$. The paired $\cdot_\text{1D}$ variables refer to the corresponding variables for the $\text{1D}$ domain, $\Omega_\text{1D}$.
For this relation, scale is important and therefore the extrusion length, $L_3$ is relevant in the evaluation of these quantities.
\end{proposition}
We provide a proof of Proposition \ref{prop:chiupsilontendom} in Appendix \ref{App:sketch_chiupsilontendom}.

We provide a table of normalized $\chi$ and $\Upsilon$ values for the canonical domains and selected triangles in Table \ref{table:canonicalchi}. Recall that $\gamma\chi$ and $\gamma^2\Upsilon$ are scale-invariant (Proposition \ref{prop:chiscale}). These values are derived by evaluation of the appropriate quadratic functionals of $\psipz$ on each of these domains (Appendix \ref{App:phispherecalc} for the case of sphere, and Appendix \ref{App:sketch_phirt} for the selected triangles). Note for the selected triangles, we must first solve for the constant $b_{00}$ to satisfy the normalization constraint $M(\psipz)=0$.

\input{Tables/canonfeval}

We now summarize the computation of $\phi$, $\Upsilon$, and $\chi$ for given domain $\Omega$ and coefficients $\sigma$, $\kappa$: Step 1, we find $\psipz_h (\approx \psipz)$ by finite element approximation of PDE (Appendix \ref{App:phisaddle}); Step 2, we evaluate $\phi_h$, $\Upsilon_h$, and $\chi_h$ -- each of which is a quadratic functional of $\psi'^0$. Clearly, Step 1 dominates the computational cost --- Step 2 is effectively free. It follows that the cost to compute $\Upsilon$ and $\chi$ --- which shall be required to form our error estimator for $\UavgLP$ --- is effectively the same as the cost to compute $\phi$ --- already required to form $\UavgLP$. It further follows that the marginal computational cost associated with extension from first-order to second-order, if we consider both prediction and error estimation, is negligible.

We now provide the error estimate associated with our second-order Pad\'e approximation:
\begin{proposition}[Asymptotic Error Bound for $\UavgLP$]\label{Prop:N3p}
We define
\begin{align}
\EavgLP(\Tf) \equiv \max_{T \in [0,T_{\text{final}}]} |U_\avg(T) - \UavgLP(T)|,
\end{align}
for $\UavgLP$ given by Definition \ref{def:slowtime} and the original QoI, $\uavgLP$, given by Definition \ref{def:lumpedapprox}.

Then the asymptotic error satisfies
\begin{align}
\EavgLP(\Tf) = \EavgL{2\text{P\,asymp}} + \pcalO(B^3) \text{ as } B \rightarrow 0 
\end{align}
where
\begin{align}
\EavgL{2\text{P\,asymp}}\equiv\left(\frac{|\gamma\chi-\gamma^2\Upsilon-\phi^2|}{e\gamma^2}+\Upsilon\right)B^2;
\end{align}
which may also be expressed with scale-invariant terms and $\Bdunk$ as
\begin{align}
\EavgL{2\text{P\,asymp}}=\left(|\gamma\chi-\gamma^2\Upsilon-\phi^2|e^{-1}+\gamma^2\Upsilon\right)\Bdunk^2.\label{eq:sinv2p}
\end{align}
\end{proposition}
Recall that the second-order asymptotic bound $\EavgL{2P\,\text{asymp}}$ requires no PDE solutions other than $\psipz$; since $\phi$, $\chi$, and $\Upsilon$ are all functional evaluations of $\psipz$. Note also that all scale-dependence in \eqref{eq:sinv2p} is reflected in $\Bdunkp$. We provide in Appendix \ref{App:PropN3p} a sketch of the proof of Proposition \ref{Prop:N3p}. 

Before proceeding to the error results for SART-1 and SART-2, we tabulate the normalized functional evaluations for each geometry, necessary for error estimation for $\UavgLP$, in Table \ref{tab:sartfeval}. We present the numerical results for the error $\EavgLP(2)$ for the two domains in Tables \ref{tab:EavgL2_sart1} and \ref{tab:EavgL2_sart2}, respectively.  We observe similar trends as for the first-order results; however we note that the bound for the second-order approximation is looser than that of the first-order approximation. This is due to uncertainty in the temporal behavior of the higher eigenmodes: requiring an additional bounding term \eqref{eq:sumbound2} which does not decay in time. However, for low $B$, we see the ratio is less than 2.5. for SART-1 and less than 2 for SART-2.

\input{Tables/sartfeval}

\input{Tables/sart1b}
\input{Tables/sart2b}

In addition to the second-order convergence of $\UavgLP$ shown in Table \ref{tab:EavgL2_sart1} and Table \ref{tab:EavgL2_sart2}, we see a remarkable reduction in error when compared against first-order $\UavgL{1}$ results from Table \ref{tab:NumRes1a} and Table \ref{tab:NumRes2a}. For the smallest $B$ we considered, $B=10^{-3}$, the error is reduced by a factor of greater than 500 for SART-1 and greater than 90 for SART-2. We note that the reduction is observed only as $B$ tends to 0 and indeed we see slight increase in the order $1$ regime of $B$. With these results, we emphasize the improvement in analysis of the dunking problem as a result of generalized framework presented in this work: computation of $\psipz$ yields a significant improvement in accuracy as well as providing an error bound which is missing from the `classical' analysis for general geometry.

%% file: Tables/canonfeval.tex
\begin{table}[H]
\begin{center}
\caption{Values of normalized functional evaluations for canonical domains and selected triangles.}  \label{table:canonicalchi}
\begin{tabular}{l|c|c|c}
\multicolumn{1}{c|}{$\Omega$} & $\phi$ & $\gamma\,\chi$ & $\gamma^2\,\Upsilon$\\
\hline\hline
$\Omega_\text{ref,1}$: Interval & 1/3 & 1/9 & 1/45\\
$\Omega_\text{ref,2}$: Disk  & 1/2 & 1/4 & 1/12\\
$\Omega_\text{ref,3}$: Sphere & 3/5 & 9/25 & 27/175\\
$\Omega_\text{right triangle}(W=1)$ & 4/3 & $\frac{4}{5}(3+2\sqrt{2})$ & $\frac{4}{15}(3+2\sqrt{2})$\\
$\Omega_\text{equilateral triangle}$ & 1 & 9/5 & 3/5\\
$\Omega_\text{right triangle}(W)$ & $\sim \frac{2}{3}W^{-2} $ & $\sim \frac{28}{15} W^{-4}$ & $ \sim \frac{28}{45} W^{-4}$
\end{tabular}
\end{center}
\end{table}

%% file: Tables/sartfeval.tex
\begin{table}[H]
\begin{center}
\caption{Values of normalized $\chi$ and $\Upsilon$ for SART-1 and SART-2 Domains.}  \label{tab:sartfeval}
\begin{tabular}{l|c|c|c|c}
\multicolumn{1}{c|}{$\Omega$} & $\phi$ & $\gamma\,\chi$ & $\gamma^2\,\Upsilon$ & $e_h^\phi$\\
\hline\hline
SART-1 & 9.14 & $\sci{4.65}{2}$ & $\sci{1.55}{2}$ & $\sci{2.07}{-13}$ \\
SART-2 & 161 & $\sci{1.21}{5}$ & $\sci{4.02}{4}$ & $\sci{1.57}{-11}$
\end{tabular}
\end{center}
\end{table}

%% file: Tables/sart1b.tex
\begin{table}[H]
	\centering
	\caption{Error in Approximation $\UavgLP$ to $U_\avg$ for $T_\text{final}=2$: SART-1 Domain.}\label{tab:EavgL2_sart1}
	\begin{tabular}{ l | c | c | c | c | c}
	\multicolumn{1}{c|}{$B$} & $\Bdunk\equiv B/\gamma$ & $\Bdunkp\equiv \phi B/\gamma$ & $\EavgLP(2)$ & $\EavgL{2\text{P\,asymp}}$ & $\EavgL{2\text{P\,asymp}}/\EavgLP(2)$  \\[.1em]
	\hline
	\hline
$\sci{1}{-3}$ & $\sci{5.48}{-5}$ & $\sci{5.01}{-4}$ & $\sci{3.24}{-7}$ & $\sci{7.16}{-7}$ & $2.21$ \\
$\sci{2}{-3}$ & $\sci{1.10}{-4}$ & $\sci{1.00}{-3}$ & $\sci{1.72}{-6}$ & $\sci{2.86}{-6}$ & $1.66$ \\
$\sci{5}{-3}$ & $\sci{2.74}{-4}$ & $\sci{2.50}{-3}$ & $\sci{1.19}{-5}$ & $\sci{1.79}{-5}$ & $1.50$ \\
$\sci{1}{-2}$ & $\sci{5.48}{-4}$ & $\sci{5.01}{-3}$ & $\sci{4.79}{-5}$ & $\sci{7.16}{-5}$ & $1.50$ \\
$\sci{2}{-2}$ & $\sci{1.10}{-3}$ & $\sci{1.00}{-2}$ & $\sci{1.87}{-4}$ & $\sci{2.86}{-4}$ & $1.53$ \\
$\sci{5}{-2}$ & $\sci{2.74}{-3}$ & $\sci{2.50}{-2}$ & $\sci{1.07}{-3}$ & $\sci{1.79}{-3}$ & $1.67$ \\
$\sci{1}{-1}$ & $\sci{5.48}{-3}$ & $\sci{5.01}{-2}$ & $\sci{3.73}{-3}$ & $\sci{7.16}{-3}$ & $1.92$ \\
$\sci{2}{-1}$ & $\sci{1.10}{-2}$ & $\sci{1.00}{-1}$ & $\sci{1.17}{-2}$ & $\sci{2.86}{-2}$ & $2.46$ \\
$\sci{5}{-1}$ & $\sci{2.74}{-2}$ & $\sci{2.50}{-1}$ & $\sci{4.23}{-2}$ & $\sci{1.79}{-1}$ & $4.24$ \\
$\sci{1}{ 0}$ & $\sci{5.48}{-2}$ & $\sci{5.01}{-1}$ & $\sci{9.35}{-2}$ & $\sci{7.16}{-1}$ & $7.66$
	\end{tabular}
\end{table}

%% file: Tables/sart2b.tex
\begin{table}[H]
	\centering
	\caption{Error in Approximation $\UavgLP$ to $U_\avg$ for $T_\text{final}=2$: SART-2 Domain.}\label{tab:EavgL2_sart2}
	\begin{tabular}{ l | c | l | c | l | c }
	\multicolumn{1}{c|}{$B$} & $\Bdunk\equiv B/\gamma$ & \multicolumn{1}{|c|}{$\Bdunkp\equiv \phi B/\gamma$} & $\EavgLP(2)$ & \multicolumn{1}{|c|}{$\EavgL{2\text{P\,asymp}}$} & $\EavgL{2\text{P\,asymp}}/\EavgLP(2)$  \\[.1em]
	\hline
	\hline
$\sci{1}{-3}$ & $\sci{1.51}{-5}$ & $\quad\sci{2.44}{-3}$ & $\sci{9.18}{-6}$ & $\sci{1.38}{-5}$ & $1.50$ \\
$\sci{2}{-3}$ & $\sci{3.03}{-5}$ & $\quad\sci{4.88}{-3}$ & $\sci{3.72}{-5}$ & $\sci{5.52}{-5}$ & $1.48$ \\
$\sci{5}{-3}$ & $\sci{7.57}{-5}$ & $\quad\sci{1.22}{-2}$ & $\sci{2.26}{-4}$ & $\sci{3.45}{-4}$ & $1.52$ \\
$\sci{1}{-2}$ & $\sci{1.51}{-4}$ & $\quad\sci{2.44}{-2}$ & $\sci{8.56}{-4}$ & $\sci{1.38}{-3}$ & $1.61$ \\
$\sci{2}{-2}$ & $\sci{3.03}{-4}$ & $\quad\sci{4.88}{-2}$ & $\sci{3.07}{-3}$ & $\sci{5.52}{-3}$ & $1.80$ \\
$\sci{5}{-2}$ & $\sci{7.57}{-4}$ & $\quad\sci{1.22}{-1}$ & $\sci{1.43}{-2}$ & $\sci{3.45}{-2}$ & $2.41$ \\
$\sci{1}{-1}$ & $\sci{1.51}{-3}$ & $\quad\sci{2.44}{-1}$ & $\sci{3.90}{-2}$ & $\sci{1.38}{-1}$ & $3.53$ \\
$\sci{2}{-1}$ & $\sci{3.03}{-3}$ & $\quad\sci{4.88}{-1}$ & $\sci{9.09}{-2}$ & $\sci{5.52}{-1}$ & $6.07$ \\
$\sci{5}{-1}$ & $\sci{7.57}{-3}$ & $\quad\sci{1.22}{ 0}$ & $\sci{2.18}{-1}$ & $\sci{3.45}{ 0}$ & $15.8$ \\
$\sci{1}{ 0}$ & $\sci{1.51}{-2}$ & $\quad\sci{2.44}{ 0}$ & $\sci{3.54}{-1}$ & $\sci{1.38}{ 1}$ & $38.9$
	\end{tabular}
\end{table}

%% file: udelta.tex
\subsection{Relative Domain-Boundary Difference QoI} \label{sec:udelta}

We now establish an error bound for the domain-boundary difference QoI.

\begin{proposition}[\ Error in Lumped $U_\Delta$ Approximation: Asymptotic Estimate] \label{prop:Deltalumpederror}

Given a (slow-)time ``cut-off'' $T_0$ for $0 < T_0 \le 1$, we define the relative error of $\UDeltaL$ (Proposition \ref{def:slowtime}) as
\begin{align}
\EDeltaLr(T_0,T;B) &\equiv \max_{T'\in[T_0,T]}\frac{\left|U_\Delta(T';B)-\UDeltaL(B)\right|}{\displaystyle \UDeltaL(B)} \;, \label{eq:eDeltareldef}
\end{align}
and introduce constants based on quadratic functional evaluations of $\psipz$ 
\begin{align}
C_0&\equiv \frac{\gamma^2\Upsilon}{e\, \phi},\\
C_1&\equiv \frac{|\gamma\chi-\gamma^2\Upsilon - \phi^2|}{\phi}.
\end{align}
Then, as $B\rightarrow 0$, we have,
\begin{align}
\EDeltaLr(T_0,\Tf;B) \le \left(\frac{C_0}{T_0} + C_1\right) \underbrace{(B/\gamma)}_{\Bdunk} +\, \pcalO(B^2) \; . \label{eq:eDeltarelbound}
\end{align}
We note that the terms encapsulated by $\pcalO(B^2)$ is independent of $\Tf$.
\end{proposition}
We provide a proof of Proposition \ref{prop:Deltalumpederror} in Appendix \ref{App:sketch_Deltalumpederror}.

We note that $\EDeltaLr(T_0,T;B)$ is the maximum error in the lumped approximation {\em relative} to the lumped approximation for $U_\Delta(T;B)$, $\UDeltaL(B)$ from time $T_0$ up to time $T$. Thus, as $B \rightarrow 0$, although $U_{\Delta}(T;B)$ is $\pcalO(B)$, the {\em absolute} error in $\UDeltaL(T;B)$ is order $\pcalO(B^2)$; thus, for $B \rightarrow 0$, we can accurately discriminate $\UDeltaL(T;B)$ as a function of $B$.

We can not hope to capture the short-time evolution of the field on the boundary with a lumped approximation. The cut-off $T_0$ restricts attention to times away from $T = 0$; as $T_0$ decreases, the error in $\UDeltaL(T;B)$ increases commensurately. However, we note that, in slow-time variables, our error estimate \eqref{eq:eDeltarelbound} is valid for $T \ge T_0$, and hence for $T_0$ small (say, $T_0 = 0.2$) the error estimate is relevant over most of the time interval of interest, $(0,\Tf]$ (say for $\Tf\ge 1$). However, we note that even for modest $T_0$, the error estimate contribution for the initial transient is too pessimistic.

We present in Table \ref{tab:EDeltaL_sart1} and Table \ref{tab:EDeltaL_sart2} the maximum relative error $\EDeltaLr(\Tf;T_0,B)$, error contributions $C_1 \Bdunk$ and $(C_0/T_0+C_1)\Bdunk$ vs. $B$ (with $\Bdunk$ and $\Bdunkp$) for SART-1 and SART-2, respectively. For both cases, we used parameters $T_0=0.2$, $\Tf=2$ and observe first-order convergence in $B$ as predicted. For SART-1, the maximum relative error is less than 8\% for $\UDeltaL$ as large as 0.025, and for SART-2 the maximum relative error is less than 5\% for $\UDeltaL$ as large as 0.025. We note that the full error bound does behave as bounds for all $B$ values considered; however, $C_1 \Bdunk$ provides very accurate error estimates as $B$ decreases.

We remark that the `classical' lumped approximation provides a $U_\Delta$ estimate of 0. Thus for applications where the QoI is $U_\Delta$, $\phi$ is an essential parameter to approximate $U_\Delta$ as a function of $B$, while $\chi$ and $\Upsilon$ are necessary to provide an error estimate with minimal marginal cost. 

\input{Tables/sart1c}
\input{Tables/sart2c}

%% file: Tables/sart1c.tex
\begin{table}[H]
	\centering
	\caption{ \ Relative Error in Approximation $\UDeltaL$ to $U_\Delta$: SART-1 Domain, $\Tf = 2$, $T_0=0.2$.}\label{tab:EDeltaL_sart1}
	\begin{tabular}{ l | c | c | l | l}

	\multicolumn{1}{c|}{$B$} & $\Bdunk\equiv B/\gamma$ & $\EDeltaLr(0.2,2;B)$ & \multicolumn{1}{|c|}{$C_1 \Bdunk$} & \multicolumn{1}{|c}{$(C_0/T_0 + C_1)\Bdunk$} \\[.1em]
	\hline
	\hline
   $\sci{1}{-3}$ & $\sci{5.48}{-5}$ & $\sci{1.36}{-3}$ & $\sci{1.81}{-3}$ & $\quad\sci{3.52}{-3}$ \\
   $\sci{2}{-3}$ & $\sci{1.10}{-4}$ & $\sci{2.71}{-3}$ & $\sci{3.61}{-3}$ & $\quad\sci{7.03}{-3}$ \\
   $\sci{5}{-3}$ & $\sci{2.74}{-4}$ & $\sci{6.74}{-3}$ & $\sci{9.03}{-3}$ & $\quad\sci{1.76}{-2}$ \\
   $\sci{1}{-2}$ & $\sci{5.48}{-4}$ & $\sci{1.33}{-2}$ & $\sci{1.81}{-2}$ & $\quad\sci{3.52}{-2}$ \\
   $\sci{2}{-2}$ & $\sci{1.10}{-3}$ & $\sci{2.64}{-2}$ & $\sci{3.61}{-2}$ & $\quad\sci{7.03}{-2}$ \\
   $\sci{5}{-2}$ & $\sci{2.74}{-3}$ & $\sci{8.32}{-2}$ & $\sci{9.03}{-2}$ & $\quad\sci{1.76}{-1}$ \\
   $\sci{1}{-1}$ & $\sci{5.48}{-3}$ & $\sci{2.08}{-1}$ & $\sci{1.81}{-1}$ & $\quad\sci{3.52}{-1}$ \\
   $\sci{2}{-1}$ & $\sci{1.10}{-2}$ & $\sci{3.80}{-1}$ & $\sci{3.61}{-1}$ & $\quad\sci{7.03}{-1}$ \\
   $\sci{5}{-1}$ & $\sci{2.74}{-2}$ & $\sci{5.82}{-1}$ & $\sci{9.03}{-1}$ & $\quad\sci{1.76}{ 0}$ \\
   $\sci{1}{ 0}$ & $\sci{5.48}{-2}$ & $\sci{6.79}{-1}$ & $\sci{1.81}{ 0}$ & $\quad\sci{3.52}{ 0}$
	\end{tabular}
\end{table}

%% file: Tables/sart2c.tex
\begin{table}[H]
	\centering
	\caption{ \ Relative Error in Approximation $\UDeltaL$ to $U_\Delta$: SART-2 Domain, $\Tf=2$, $T_0=0.2$.}\label{tab:EDeltaL_sart2}
	\begin{tabular}{ l | c | c | l | l }
	\multicolumn{1}{c|}{$B$} & $\Bdunk\equiv B/\gamma$ & $\EDeltaLr(0.2,2;B)$ & \multicolumn{1}{|c|}{$C_1 \Bdunk$} & $(C_0/T_0 + C_1)\Bdunk$ \\[.1em]
	\hline
	\hline
    $\sci{1}{-3}$ & $\sci{1.51}{-5}$ & $\sci{5.07}{-3}$ & $\sci{7.53}{-3}$ & $\quad\sci{1.45}{-2}$ \\
    $\sci{2}{-3}$ & $\sci{3.03}{-5}$ & $\sci{1.01}{-2}$ & $\sci{1.51}{-2}$ & $\quad\sci{2.89}{-2}$ \\
    $\sci{5}{-3}$ & $\sci{7.57}{-5}$ & $\sci{2.48}{-2}$ & $\sci{3.77}{-2}$ & $\quad\sci{7.24}{-2}$ \\
    $\sci{1}{-2}$ & $\sci{1.51}{-4}$ & $\sci{5.49}{-2}$ & $\sci{7.53}{-2}$ & $\quad\sci{1.45}{-1}$ \\
    $\sci{2}{-2}$ & $\sci{3.03}{-4}$ & $\sci{1.51}{-1}$ & $\sci{1.51}{-1}$ & $\quad\sci{2.89}{-1}$ \\
    $\sci{5}{-2}$ & $\sci{7.57}{-4}$ & $\sci{3.92}{-1}$ & $\sci{3.77}{-1}$ & $\quad\sci{7.24}{-1}$ \\
    $\sci{1}{-1}$ & $\sci{1.51}{-3}$ & $\sci{5.70}{-1}$ & $\sci{7.53}{-1}$ & $\quad\sci{1.45}{ 0}$ \\
    $\sci{2}{-1}$ & $\sci{3.03}{-3}$ & $\sci{6.98}{-1}$ & $\sci{1.51}{ 0}$ & $\quad\sci{2.89}{ 0}$ \\
    $\sci{5}{-1}$ & $\sci{7.57}{-3}$ & $\sci{7.93}{-1}$ & $\sci{3.77}{ 0}$ & $\quad\sci{7.24}{ 0}$ \\
    $\sci{1}{ 0}$ & $\sci{1.51}{-2}$ & $\sci{8.28}{-1}$ & $\sci{7.53}{ 0}$ & $\quad\sci{1.45}{ 1}$
	\end{tabular}
\end{table}

%% file: appa.tex
\section{Weak Formulations}\label{AppA}

\subsection{Bilinear Forms}\label{App:forms}

We  first define $L^2(\Omega)$, the space of functions $w$ which are square-integrable over $\Omega$, and associated norm $\| w\|_{L^2(\Omega)} \equiv (\int_\Omega w^2)^{1/2}$. We next introduce the space $H^1(\Omega) \equiv \{w \in L^2(\Omega)\,|\, |\bnabla w| \in L^2(\Omega)\}$, and associated norm
\begin{align}
\| w \|_{H^1(\Omega)} \equiv \int_\Omega w^2 + |\bnabla w|^2 \; ; \label{eq:Xnorm}
\end{align}
note that this norm shall only be applied to nondimensional quantities, and hence we do not need to include any constants for dimensional homogeneity. Finally, for convenience, we define $X \equiv X(\Omega) \equiv H^1(\Omega)$ and $\|\,\cdot\, \|_X \equiv \| \, \cdot \, \|_{H^1(\Omega)}$.

We then introduce, for all $w \in X$ and all $v \in X$,
\begin{align}
m(w,v;\sigma) \equiv & \int_\Omega \sigma\, w\,v \, , \label{eq:mdef} \\
a_0(w,v;\kappa) \equiv & \int_\Omega \kappa\, \boldsymbol{\nabla}w \cdot \boldsymbol{\nabla} v \, , \label{eq:a0def} \\
a_1(w,v) \equiv & \int_{\partial\Omega} w\, v \, , \label{eq:a1def}\\[.1em]
a(w,v;B;\kappa) \equiv & \,a_0(w,v;\kappa) + B a_1(w,v) \, .\label{eq:adef}
\end{align}
It is readily demonstrated that $m$, $a_0$, $a_1$, and $a$ are (i) continuous in $X$, (ii) symmetric, and (iii) positive semi-definite. Note we shall, as convenient, suppress explicit reference to the coefficient functions $\sigma$ and $\kappa$.

\subsection{Initial Value Problem}\label{App:ivp}

We shall  require $\sigma \in L^\infty(\Omega)$ and $\kappa \in L^\infty(\Omega)$. 
We then seek, for given $B \in \Rzp$, $u(\cdot,t;B)$ for $t \in (0,\tf]$ such that
\begin{align}
m(\partial_t u(\cdot,t;B),v;\sigma) + a(u(\cdot,t;B),v;B;\kappa) = 0, \forall v \in X \, , \label{eq:ivpweak_1}
\end{align}
subject to initial condition 
\begin{align}
u(\cdot,t = 0;B) = 1 \;. \label{eq:ivpweak_3}
\end{align} 
Recall that $\partial_t$ refers to the partial derivative with respect to time. 

The Robin heat equation admits a unique solution in Bochner space $C^0([0,\tf],L^2(\Omega))\cap L^2((0,\tf),H^1(\Omega))$; here $C^0([0,\tf])$ refers to the space of continuous functions over the interval $[0,\tf]$. It follows from $u \in L^2((0,\tf),H^1(\Omega))$ that $u(\cdot,t;B) \in X$  for almost all $t \in [0,\tf]$; the qualification ``almost all'' in particular addresses the issue of short-time singularity.

\subsection{Eigenproblem}\label{App:evp}

There is a very large literature on the Laplacian eigenproblem for Dirichlet, Neumann, and Robin boundary conditions; see \cite{GandN} for a review.

We seek $(\psi_j(B)\in X,\lambda_j(B) \in \RR_{0+})_{j = 1,2,\ldots}$ such that
\begin{align}
a(\psi_j(B),v;B;\kappa) & =  \lambda_j(B) m(\psi_j(B),v;\sigma), \forall v \in X \; , \label{eq:evpweak_1}
\end{align}
subject to normalization
\begin{align}
m(\psi_j(B),\psi_j(B);\sigma) & =  1 \; . \label{eq:evpweak_3}
\end{align}
Recall $\RR_{0+}$ refers to the non-negative reals.

We recall a few points from the main document: It is standard to note that  the eigenvalues are real and non-negative for $B \ge 0$. We shall enumerate the eigenvalues in increasing magnitude, hence $0 \le \lambda_1(B) <  \lambda_2(B) \le \lambda_3(B) \le \ldots$. It is also standard to note that the eigenfunctions are orthonormal:
\begin{align} 
m(\psi_i(B), \psi_j(B); \sigma) = \delta_{ij}\; , \; i,j = 1,2,\dots \; ; \label{eq:efortho2}
\end{align}
here $\delta_{ij}$ is the Kronecker-delta symbol. We can also demonstrate that $\lambda_1(B)$ is of multiplicity unity and that $\psi_1(B)$ is one-signed --- we shall choose $\psi_1(B)$ positive, which, with \eqref{eq:evpweak_3}, uniquely specifies $\psi_1(B)$.
We shall henceforth adopt the convention that $\psi(B)$ (respectively, $\lambda(B)$) without subscript shall refer to the first eigenfunction, $\psi_1(B)$ (respectively, the first eigenvalue, $\lambda_1(B)$), and that furthermore $\psiz(B)$ (respectively, $\lambdaz(B)$) shall refer to $\psi(B = 0)$ (respectively, $\lambda(B = 0))$.

\subsection{Sensitivity Equation: Weak Form}\label{App:phieqn}

We provide here the weak form for the sensitivity equation introduced in (strong form in) Proposition \ref{prop:EEforSensitivityDerivative}. We introduce space $Z_0(\sigma) \equiv \{w \in X\,|\, m(w,1;\sigma) = 0\} \subset H^1(\Omega)$; we imbue $Z_0(\sigma)$ with the standard $H^1(\Omega)$ inner product and associated norm. We first provide a weak form for test space $X$: Find $\psipz(\sigma,\kappa) \in Z_0(\sigma)$ unique solution  to
\begin{align}
a_0(\psipz,v;\kappa) =  |\Omega|^{-1/2} L(v;\sigma) \; , \forall v \in X \; . \label{eq:phiellp_0A}
\end{align}
where linear form
\begin{align}
L(v;\sigma) \equiv  \gamma \int_\Omega \sigma \,v  -   \int_{\partial\Omega}  v.\label{eq:Ldef}
\end{align}
We also provide an equivalent weak form for test space $Z_0(\sigma)$:
Find $\psipz(\sigma,\kappa) \in Z_0(\sigma)$ unique solution  to
\begin{align}
a_0(\psipz,v;\kappa) = - |\Omega |^{-1/2}  \int_{\partial\Omega}  v \; , \forall v \in Z_0(\sigma) \; . \label{eq:phiellp_0}
\end{align}
In either case we may then evaluate
\begin{align}
\phi(\sigma,\kappa)= a_0(\psipz,\psipz;\kappa)  \;; \label{eq:phiellp_a}
\end{align}
we recall that $a_0(w,v;\kappa) \equiv \int_\Omega \kappa(x) \bnabla w \cdot \bnabla v$. The derivation of these weak forms is given in Appendix \ref{App:sketch_EEforSensitivityDerivative}.

Note that, in addition to the constrained elliptic problems \eqref{eq:phiellp_0A} or \eqref{eq:phiellp_0}, we may also formulate $\psipz$ as a maximization principle --- useful for theoretical purposes, and a saddle problem --- useful for numerical purposes.
We discuss the former in context in Appendix \ref{App:sketch_phiprops}, and the latter in the next section, Appendix \ref{App:phisaddle}.

\subsection{Sensitivity Equation: Saddle Form}\label{App:phisaddle}

We now look for the pair $(\psipz \in X, p \in \RR)$ such that
\begin{align}
a_0(\psipz,v;\kappa) + p \, m(1,v;\sigma) & = -|\Omega|^{-1/2} a_1(1,v),\; \forall v \in X \label{eq:phisaddle1} \\
q\, m(1,\psipz;\sigma) & = 0 \; , \forall q \in \RR \; . \label{eq:phisaddle2}
\end{align}
It can readily be demonstrated that this saddle problem is inf-sup stable. Furthermore, we may explicitly evaluate $p$: we choose $v = 1$ in \eqref{eq:phisaddle1} to find $p = -|\Omega|^{-1/2} \gamma$; we may then move $p \, m(1,v;\sigma) $ in \eqref{eq:phisaddle1} to the right-hand side to directly recover \eqref{eq:phiellp_0A}. 

In the computational context, and in particular in the finite element context, the saddle problem is convenient: \eqref{eq:phisaddle1} -- \eqref{eq:phisaddle2} yields a non-singular, square, and sparse system amenable to fast direct solution techniques.

\subsection{Lumped Approximation from Galerkin Projection} \label{App:sketch_galerkinlumped}

We restrict space $X$ to $\PP_0$ in \eqref{eq:ivpweak_1} i.e, $v,u\in \PP_0\subset X$,  and proceed to simplify the bilinear forms for constant functions for all $t\in(0,t_\text{final}]$:

\begin{align} 
(v\partial_t u)\, m(1,1) + (vu)\,a_0(1,1) + (B v u)\, a_1(1,1) &= 0, \text{ hence}\\
(v\partial_t u)\, |\Omega| + (vu)\, 0 +  (B v u)\, |\partial\Omega| &= 0, \text{ hence}\\
v(\partial_t u + B\gamma u) &= 0, \text{ hence}\\
\partial_tu+B\gamma u &= 0,
\end{align}
subject to initial condition
\begin{align}
\u(0) = 1 \; . \label{eq:Uic}
\end{align}
The solution to this initial value problem is the lumped approximation, $\uavgL{1}$.

%% file: appb.tex
\section{Closed-Form Derivation: $\phi$ for Homogeneous Sphere}\label{App:phispherecalc}

We consider a sphere of radius $R$ and hence $|\Omega| = 4\pi R^3/3$ and $\gamma = 3/R$. In spherical coordinates with spherical (radial) coordinate $r$, our problem for $\psipz$ takes the form
\begin{align}
-\frac{d}{dr}\left(r^2\dfrac{d\psipz}{dr}\right) =  |\Omega|^{-1/2} \gamma r^2 \; , \label{eq:Appb1}
\end{align}
with boundary condition (boundedness at $r = 0$ and)
\begin{align}
\dfrac{d\psipz}{dr} = - |\Omega|^{-1/2} \; ; \label{eq:Appb2}
\end{align}
we observe that equation \eqref{eq:Appb1} -- \eqref{eq:Appb2} is solvable. Once we obtain $\psipz$, and particular $d\psipz/dr$, we can evaluate $\phi$ as
\begin{align}
\phi & \equiv 4\pi \int_0^R r^2 \left(\dfrac{d\psipz}{dr}\right)^2 dr   \; .\label{eq:Appb123} 
\end{align}
It remains to determine $d\psipz/dr$.
 
We first integrate \eqref{eq:Appb1} to obtain
\begin{align}
\dfrac{d\psipz}{dr}(r) = -\dfrac{1}{3}|\Omega|^{-1/2} \gamma r \; ; \label{eq:Appb3}
\end{align}
note that our boundary condition \eqref{eq:Appb2} is automatically satisfied thanks to the solvability compatibility condition.
We next integrate \eqref{eq:Appb3}  to find
\begin{align}
\psipz(r) = -\dfrac{1}{6}|\Omega|^{-1/2}\gamma r^2  +  C \; . \label{eq:Appb4}
\end{align}
Here $C$ is a constant which serves to satisfy our normalization condition; note $\phi$ does not depend on $C$.  We then insert \eqref{eq:Appb3} and our expressions for $|\Omega|$ and $\gamma$ into \eqref{eq:Appb123}  to obtain $\phi = 3/5$.

%% file: appc.tex
\section{Sketches of Proofs}\label{App:sketches}

We note here that the statements and proofs of Appendix sections \ref{App:unity} -- \ref{App:sketch_Bzero} and \ref{App:sketch_l1ltBg} are elementary or standard, but we include them here for completeness.

\subsection{Proof of Lemma \ref{lem:rou}} \label{App:sketch_rou}\label{App:unity}

To begin, we invoke $L^2(\Omega)$ completeness of the eigenfunctions to write 
\begin{align}
\lim_{N\rightarrow \infty}\, \norm{\; 1 - \sum_{k=1}^N b_k \psi_k(\cdot;B)\;}_{L^2(\Omega)} \rightarrow 0 \; , \label{eq:Mpf_10}
\end{align}
or in abbreviated form (to be interpreted as \eqref{eq:Mpf_10})
\begin{align}
\sum_{k=1}^\infty b_k \psi_k(\cdot;B) = 1\; , \label{eq:Mpf_1}
\end{align}
where
\begin{align}
b_k = \OmMeas M(\psi_k(B)), k = 1,2,\ldots \; ; \label{eq:Mpf_2}
\end{align}
note to derive \eqref{eq:Mpf_2}, we multiply \eqref{eq:Mpf_1} by $\sigma \psi_{k'}(B)$, integrate over $\Omega$, invoke orthonormality \eqref{eq:efortho2}, and apply the definition of $M$, \eqref{eq:Mdef}.
We now insert \eqref{eq:Mpf_2} in \eqref{eq:Mpf_1} to obtain \eqref{eq:Mpsi_one}: 
\begin{align}
\lim_{N\rightarrow \infty} \,\sum_{k=1}^\infty \OmMeas M(\psi_k(B)) \psi_k(\cdot;B) = 1\; . \label{eq:Mpf_3}
\end{align}
Finally, we observe that, for any constant function $C=c$, $M(C) = \frac{1}{\OmMeas}\int_\Omega \sigma C = c$, since the domain mean of $\sigma$ is, by construction, unity; application of (linear) functional $M$ to \eqref{eq:Mpf_3} thus yields \eqref{eq:Mpsisquare}.

\subsection{Proof of Proposition \ref{prop:sov}} \label{App:sketch_sov}

Each term in the sum \eqref{eq:sov_sum} satisfies the initial value problem \eqref{eq:ivpweak_1} and hence too a (finite) sum satisfies \eqref{eq:ivpweak_1}. Note that the sum of the initial conditions corresponding to each of the term approaches the initial condition \eqref{eq:ivpweak_3} as $N\rightarrow\infty$. Equation \eqref{eq:sov_lim} then follows from Lemma \ref{lem:rou} and the stability of the heat equation in $C^0([0,\tf],L^2(\Omega))\cap L^2((0,\tf),H^1(\Omega))$ with respect to $L^2(\Omega)$ perturbations in initial conditions.

\subsection{Proof of Proposition \ref{prop:sov_QoI}} \label{App:sketch_sov_QoI}

We first note $M$ and $H$ are continuous linear functionals in $X$, and since from Proposition \ref{prop:sov}, $u^N(\cdot,t;B)$ of \eqref{eq:sov_sum0} converges to $u(\cdot,t;B)$ in $X$, hence
\begin{align}
\lim_{N\rightarrow\infty}M(u^N(\cdot,t;B)) = M(u(\cdot,t;B)) = u_\avg(t;B)
\end{align}
and
\begin{align}
\lim_{N\rightarrow\infty}H(u^N(\cdot,t;B)) = H(u(\cdot,t;B)) = u_\paavg(t;B)
\end{align}
for all $t \in (0,\tf]$. We first consider $u_\avg$: from \eqref{eq:sov_sum}, we directly obtain \eqref{eq:sov_avgQoI}
\begin{align}
u_\avg(t;B) &= M\left(   \sum_{k = 1}^\infty M(\psi_k(B)) \psi_k(\cdot;B) \exp(-\lambda_k(B) t) \right) \nonumber\\
&= \sum_{k = 1}^\infty (M(\psi_k(B)))^2  \exp(-\lambda_k(B) t) )  \; . \label{eq:sov_uavg_pf}
\end{align}
We next consider $u_\paavg$: we first obtain
\begin{align}
u_\paavg(t;B) = H\left(   \sum_{k = 1}^\infty M(\psi_k(B)) \psi_k(\cdot;B) \exp(-\lambda_k(B) t) \right)\nonumber\\
= \sum_{k = 1}^\infty M(\psi_k(B)) H(\psi_k(B))  \exp(-\lambda_k(B) t) ) \; ; \label{eq:sov_upaavg_pf1}
\end{align}
we now note from \eqref{eq:evpweak_1}, \eqref{eq:mdef} -- \eqref{eq:adef}, and the definitions of $H$ \eqref{eq:Hdef} and $\gamma$ \eqref{eq:gammadef} that
\begin{align}
H(\psi_k) = \dfrac{M(\psi_k)\lambda_k}{B \gamma} \; ; \label{eq:sov_upaavg_pf2}
\end{align}
\eqref{eq:sov_paavgQoI} then directly follows from \eqref{eq:sov_upaavg_pf1} and \eqref{eq:sov_upaavg_pf2}.

\subsection{Proof of Lemma \ref{lem:Bzero}}\label{App:sketch_Bzero} The Rayleigh Quotient associated with our eigenproblem \eqref{eq:evpweak_1} is simply
\begin{align}
\calR(v;\kappa,\sigma,B) \equiv \frac{a_0(v,v;\kappa)+Ba_1(v,v)}{m(v,v;\sigma)}.\label{eq:RQB}
\end{align}
Thus, the minimizer of $\calR(v;B=0)$ over all functions $v$ in $X$ is the constant function, $C$; we then apply our normalization condition \eqref{eq:evpweak_3}, and recall that $\sigma$ is unity mean, to obtain $\psiz = \OmMeas^{-1/2}$. The minimum of the Rayleigh Quotient is zero, hence $\lambdaz = 0$. 

\subsection{Proof of Lemma \ref{lem:eigfcnamp}}\label{App:sketch_eigfcnamp}

To begin, we summarize the notation for the proofs of Lemma \ref{lem:eigfcnamp} here and also Proposition \ref{prop:eigvasymp}. We shall denote $\psi_1(B)$ by $\psi(B)$, the first and second derivatives of $\psi$ with respect to $B$ as $\psip(B)$ and $\psipp(B)$, respectively, and $\psi(B = 0)$, $\psip(B=0)$, and $\psipp(B=0)$ as $\psiz$, $\psipz$, and $\psippz$, respectively. Turning now to the eigenvalue, we denote $\lambda_1(B)$ by $\lambda(B)$, the first and second derivatives of $\lambda$ with respect to $B$ as $\lambdap(B)$ and $\lambdapp(B)$, respectively, and furthermore $\lambda(B = 0)$, $\lambdap(B=0)$, and $\lambdapp(B=0)$ as $\lambdaz$, $\lambdapz$, and $\lambdappz$, respectively.

In this section, for clarity, we will suppress the parametric argument $\sigma$ to the bilinear form $m$. We first recall that $\lambda = \lambda_1$ is of multiplicity unity; we further note that our (normalized) eigenfunction, $\psi = \psi_1$, is unique. We next note that $\psi$ and $\lambda$ admit derivatives of any order, as can be demonstrated by continuity arguments and induction. Note that, for the higher modes, we must be careful to consider a branch of the eigenproblem. However, no such difficulty exists for the first mode: the branch of the first mode, which we recall is of multiplicity unity, remains distinct from the branches of the higher modes --- there is no mode crossing.

Proceeding accordingly, we find from \eqref{eq:evpweak_3} that
\begin{align}
0 = \frac{d}{dB} (m(\psi,\psi)) = 2 m(\psi',\psi)\;,  \text{ for any }B \in \Rzp \; . \label{eq:dB1}
\end{align}
We then differentiate a second time to find
\begin{align}
0 = \frac{d}{dB} (m(\psi',\psi)) = m(\psi',\psi') + m(\psi'',\psi) \; , \text{ for any }B \in \Rzp \; ,\label{eq:dB2}
\end{align}
or
\begin{align}
m(\psi'',\psi) = - m(\psi',\psi')  \; , \text{ for any }B \in \Rzp \; .\label{eq:dB3}
\end{align}
We then invoke \eqref{eq:Mdef}, \eqref{eq:mdef}, Lemma \ref{lem:Bzero}, and \eqref{eq:dB1} (for $B = 0$)  to show
\begin{align}
M(\psipz) = \OmMeas^{-1} m(\psipz,1) = \OmMeas^{-1/2} m(\psipz,\psiz) =  0 \; . \label{eq:Mpsiprime}
\end{align}
In a similar fashion, now invoking \eqref{eq:dB3} for $B = 0$, we find that
\begin{align}
\begin{split}
M(\psippz) = \OmMeas^{-1} m(\psippz,1) & = \OmMeas^{-1/2} m(\psippz,\psiz) \\& =  -\OmMeas^{-1/2} m(\psipz,\psipz) = -\OmMeas^{-1/2}\Upsilon \; , \label{eq:dB4}
\end{split}
\end{align}
for $\Upsilon$ given by \eqref{eq:Upsilondef}. 

We now consider application of $M$ to a Taylor series for $\psi$ around $B = 0$:
\begin{align}
M(\psi(B)) & =  M(\psiz + B\psipz + \frac{B^2}{2}\psippz + \pcalO(B^3)) \; \label{eq:dB5a} \\
& = \OmMeas^{-1/2}(1 - \frac{B^2}{2}\Upsilon + \pcalO(B^3)) \; \label{eq:dB5b}
\end{align}
and hence
\begin{align}
(M(\psi(B)))^2 = \OmMeas^{-1}( 1 - B^2\Upsilon + \pcalO(B^3)) \text{ as } B \rightarrow 0\;. \label{eq:dB6}
\end{align}
This concludes the proof.

\subsection{Proof of Proposition \ref{prop:eigvasymp}}\label{App:sketch_eigvasymp}

 We first note, from Lemma \ref{lem:Bzero} and standard Taylor series analysis, that 
\begin{align}
\lambda(B) \sim 0 + \lambdapz B + \lambdappz \dfrac{B^2}{2} + \pcalO(B^3) \text{ as } B \rightarrow 0 \; .
\end{align}
We must thus demonstrate that $\lambdapz = \gamma$ and $\lambdappz = -2\phi$. In this section, for clarity, we suppress explicit reference to $\sigma$ and $\kappa$.

We first address $\lambdapz$. To begin, we differentiate \eqref{eq:evpweak_1} (for $j = 1$) with respect to $B$ to find
\begin{align}
a_0(\psip,v) + Ba_1(\psip,v) = \lambda \, m(\psip,v) + (\lambdap \,m(\psi,v) - a_1(\psi,v)), \;\forall v \in X\; \; . \label{eq:lprimepos2A}
\end{align}
We now evaluate \eqref{eq:lprimepos2A} for $B = 0$ and apply Lemma \ref{lem:Bzero}  to obtain
\begin{align}
a_0(\psipz,v) = \lambdapz m(\psiz,v) - a_1(\psiz,v)\; , \; \forall v \in X\; . \label{eq:lprimepos2B}
\end{align}
We now choose $v = \psiz$, invoke Lemma \ref{lem:Bzero}, the definition of $a_0$ \eqref{eq:a0def}, and \eqref{eq:evpweak_3} to obtain
\begin{align}
\lambdapz = a_1(\psiz,\psiz) \; . \label{eq:lprimepos3}
\end{align}
Finally, we apply Lemma \ref{lem:Bzero}, the definition of $a_1$ \eqref{eq:a1def}, and the definition of $\gamma$ \eqref{eq:gammadef}, to obtain $\lambdapz = \gamma$.

We next address $\lambdappz$. To begin, we return to \eqref{eq:lprimepos2B}  but now choose $v = \psipz$ to obtain, upon application of \eqref{eq:dB1} (for $B = 0$),
\begin{align}
a_0(\psipz,\psipz) = -a_1(\psiz,\psipz) \; . \label{eq:phiell_6}
\end{align}
We next consider the second derivative of \eqref{eq:evpweak_1} with respect to $B$, or alternatively the first derivative of \eqref{eq:lprimepos2A}: 
\begin{align}
a_0(\psipp,v) + Ba_1(\psipp,v) = \lambda m(\psipp,v) + 2 \lambdap m(\psip,v) + \lambda'' m(\psi,v) - 2a_1(\psip,v) \;, \forall v \in X \; ; \label{eq:phiell_7}
\end{align}
evaluation of \eqref{eq:phiell_7} at $B = 0$ (recall $\lambdaz = 0$ and $\lambdapz = \gamma$) then yields
\begin{align}
a_0(\psippz,v) = 2 \gamma\, m(\psipz,v) + \lambdappz m(\psiz,v) - 2a_1(\psipz,v) \;, \forall v \in X \; . \label{eq:phiell_8}
\end{align}
We now take $v = \psiz$ in \eqref{eq:phiell_8} and recall that $a_0(\cdot,\psiz) = 0$ from \eqref{eq:a0def}, $m(\psipz,\psiz) = 0$ from \eqref{eq:dB1} (for $B = 0$), and $m(\psiz,\psiz) = 1$ from \eqref{eq:evpweak_3}; we thereby obtain
\begin{align}
\lambdappz = 2 a_1(\psipz,\psiz) \; , \label{eq:phiell_9}
\end{align}
which then, with $\phi = -\lambda''_0/2$, \eqref{eq:phiell_6}, and symmetry yields \eqref{eq:phiellp_a}. Strict positivity of $\phi$ is established by Proposition \ref{prop:phiprop1}.

\subsection{Proof of Proposition \ref{prop:EEforSensitivityDerivative}}\label{App:sketch_EEforSensitivityDerivative}

In this section, for clarity, we suppress explicit reference to $\sigma$ and $\kappa$. Our point of departure is \eqref{eq:lprimepos2B} but now, from Proposition \ref{prop:eigvasymp}, we replace $\lambdapz$ with $\gamma$:
\begin{align}
a_0(\psipz,v) = \gamma m(\psiz,v) - a_1(\psiz,v) \;, \forall v \in X \; . \label{eq:phiell_3}
\end{align}
We note that \eqref{eq:phiell_3} is singular but solvable: $a_0(\psipz,1) = 0$; $\gamma m(\psiz,1) = \gamma |\Omega | ^{1/2} M(1) = |\partial \Omega | |\Omega |^{-1/2}$ and $a_1(\psiz,1) = |\partial \Omega | |\Omega |^{-1/2} H(1)  =  |\partial \Omega | |\Omega |^{-1/2}$ and thus $\gamma m(\psiz,1) - a_1(\psiz,1) = 0$. Note the definitions of $M$ and $H$ are provided in \eqref{eq:Mdef} and \eqref{eq:Hdef}, respectively. The condition \eqref{eq:dB1} evaluated for $B = 0$, equivalent to $\psipz \in Z_0(\sigma)$, renders $\psipz$ unique. We now note that $\gamma m(\psiz,v) =  \OmMeas^{-1/2} \gamma \int_\Omega  \sigma v$. 
We have thus derived the weak form for test space $X$, \eqref{eq:phiellp_0A}, and hence also (under smoothness hypotheses) the strong form of Proposition \ref{prop:EEforSensitivityDerivative}.
 
It remains to derive the weak form for test space $Z_0(\sigma)$, \eqref{eq:phiellp_0}. In fact, since $a_0(w,v=1) = \gamma m(\psiz,v=1) - a_1(\psiz,v=1) = 0\;, \forall w \in X$, we may directly state our problem over the constrained space $Z_0(\sigma) \equiv \{w \in X \,|\, m(w,1) = 0 \}$: find $\psipz \in Z_0(\sigma)$ (equivalent to condition \eqref{eq:dB1} evaluated for $B = 0$) such that
\begin{align}
a_0(\psipz,v) = - |\Omega |^{-1/2} \int_{\partial\Omega} v \;, \forall v \in Z_0(\sigma) \; . \label{eq:phiell_5}
\end{align}
The first term on the right-hand side of \eqref{eq:phiell_3} no longer appears in \eqref{eq:phiell_5} --- since $m(\psi_0,v) = |\Omega |^{-1/2} m(1,v) = |\Omega |^{-1/2} m(v,1) = 0$ for any $v \in Z_0(\sigma)$; the second term on the right-hand side of \eqref{eq:phiell_3} appears as $ |\Omega |^{-1/2} \int_{\partial\Omega}  \,v$ in \eqref{eq:phiell_5} --- from the definition of $a_1$, \eqref{eq:a1def}, and Lemma \ref{lem:Bzero}.  We have thus derived, as desired, \eqref{eq:phiellp_0}. We can now more rigorously demonstrate well-posedness: we note that $a_0$ is continuous and coercive (Poincar\'e-Wirtinger inequality) in $Z_0(\sigma)$, and that $\int_{\partial\Omega} v$ is continuous in $Z_0(\sigma) \subset H^1(\Omega)$ (trace theorem); Lax-Milgram thus provides existence and uniqueness.

\subsection{Proof of Propositions \ref{prop:phiprop1} -- \ref{prop:phiprop4}}\label{App:sketch_phiprops}

We begin with Proposition \ref{prop:phiprop1}. We observe from \eqref{eq:phidef1} that $\phi$ will be zero only if $\psipz =C$ over $\Omega$ for some constant $C$. Let us assume $\psipz = C$. But from our normalization condition, $\psipz \in Z_0(\sigma) \equiv \{w \in X\,|\, \int_{\Omega} \sigma w =0\}$, we conclude that $C = 0$. The latter is clearly inconsistent with  \eqref{eq:xieqn_strong1}, and our hypothesis thus leads to a contradiction.

We next turn to Proposition \ref{prop:phiprop2}. The invariance with respect to translation and rotation is classical, and we thus focus on the dilation term.
We introduce a new domain $\Omega_2$ given by $\Omega_1$ scaled by factor $\beta$, written $\Omega_2=\beta\Omega_1$. It immediately follows that
\begin{align}
  |\Omega_2|&=\beta^d|\Omega_1|,\\
  \gamma_2&=\frac 1 \beta \gamma_1,\label{eq:gammascale}
\end{align}
where $d$ refers to spatial dimensionality, $\Omega_1\subset\RR^d,\Omega_2\subset\RR^d$.
Next we present the sensitivity equation and boundary condition in this new domain $\Omega_2$ ($\xi_2(\cdot) \equiv \psi'^0(\cdot;\Omega_2)$):
\begin{align}
    -\nabla^2\xi_2&=-|\Omega_2|^{-1/2}\gamma_2\text{ in }\Omega_2,\\
    -\partial_n\xi_2&=-|\Omega_2|^{-1/2}\text{ on }\partial\Omega_2,
\end{align}
We now introduce transformed coordinates $\hat x = x/\beta$ on $\Omega_1$ and transformed solution $\hat \xi_2\left(x/\beta\right)=\xi_2(x)$. The spatial derivative in this transformed space is
\begin{align}
    \frac{\partial}{\partial x}\xi_2 = \frac{\partial}{\partial x} \hat \xi_2\left(\frac{x}{\beta}\right) = \frac 1 \beta \frac {\partial \hat \xi_2}{\partial \hat x}.
\end{align}
Therefore the equations can be written as
\begin{align}
-\frac{1}{\beta^2}\hat{\nabla}^2\,\hat\xi_2&=-\frac{\beta^{-d/2}}{\beta}|\Omega_1|^{-1/2}\gamma_1,\text{ hence}\nonumber\\
  -\hat{\nabla}^2\,\left(\beta^{d/2-1}\hat{\xi}_2\right) &= -|\Omega_1|^{-1/2}\gamma_1 \text{ on }\Omega_1,\label{eq:scaledeq}
\end{align}
with boundary condition
\begin{align}
    -\frac{1}{\beta}\partial_{\hat{n}}\hat{\xi}_2&=-\beta^{-d/2}|\Omega_1|^{-1/2},\text{ hence}\nonumber\\
    -\partial_{\hat{n}}\left(\beta^{d/2-1}\hat{\xi}_2\right)&=-|\Omega_1|^{-1/2}\text{ on $\partial\Omega_1$}.\label{eq:scaledbc}
\end{align}
We conclude from \eqref{eq:scaledeq} -- \eqref{eq:scaledbc} that $\beta^{d/2-1}\hat{\xi}_2 = \xi_1 \equiv \psi'^0(\cdot;\Omega_1)$. We may now evaluate
\begin{align}
  \phi_2 &= \int_{\Omega_2}\bnabla \xi_2\cdot\bnabla \xi_2=\beta^d \int_{\hat\Omega}\frac{1}{\beta^2}\left|\hat{\bnabla}\hat\xi_2\right|^2\nonumber\\
  &=\beta^{d-2}\cdot\beta^{2(-d/2+1)}\int_{\Omega_1}|\bnabla \xi_1|^2\nonumber\\
  &=\phi_1.
\end{align}

We next consider Proposition \ref{prop:phiprop4}.  We first introduce a functional $J(w;\kappa)$ defined by
\begin{align}
J(w;\kappa) \equiv \left[ - a_0(w,w;\kappa) - 2 a_1(w,\psiz)\right]\; ; \label{eq:Jdef}
\end{align}
we can then show from \eqref{eq:phiellp_0} and standard variational arguments that
\begin{align}
\phi = \max_{w \in Z_0(\sigma)} J(w;\kappa)\,,\; \psipz = \arg\max_{w \in Z_0(\sigma)}  J(w;\kappa) \; . \label{eq:maxprinciple}
\end{align}
Note the Lagrangian formulation of this maximization principle yields a saddle problem, described in Appendix \ref{App:phisaddle}, which is very well-suited for numerical approximation and solution of equation for $\psipz$ and then subsequent evaluation of $\phi$.
We next introduce $\psipz_1 \in Z_0(\sigma)$ and $\psipz_2 \in Z_0(\sigma)$ (note $\sigma$ is fixed here), solutions to \eqref{eq:phiellp_0} associated to $\kappa = \kappa_1(x)$ and $\kappa = \kappa_2(x)$, respectively. Now let us require $\kappa_2(x)\ge \kappa_1(x)$, $\forall x \in \Omega$. We may then conclude the proof:
\begin{align}
\phi_2 = J(\psipz_2;\kappa_2) \le J(\psipz_2;\kappa_1) \le \max_{w \in Z_0} J(w;\kappa_1) = \phi_1 \; ,
\end{align}
where the first step follows from the definition of $J(\cdot;\kappa)$ \eqref{eq:Jdef} and $a_0$ \eqref{eq:a0def}, and the second step follows from our maximization principle \eqref{eq:maxprinciple}.

\subsection{Proof of Proposition \ref{prop:phisigma}} \label{App:sketch_phisigma}

To begin, we introduce a (slightly) new problem: Find $\xi \in Z_0(1)$ such that
\begin{align}
a_0(\xi,v) = |\Omega|^{-1/2} L(v;\sigma)\; , \forall v \in X \; , \label{eq:ps1}
\end{align}
where $a_0(\cdot,\cdot)$ refers to $a_0(\cdot,\cdot;1)$ and recall $Z_0(b) \equiv \{w \in X \,|\, m(w,1;b) = 0$. We observe that \eqref{eq:ps1} is solvable thanks to the definitions of $L$, \eqref{eq:Ldef}, and $\gamma$, \eqref{eq:gammadef}. Note that the equation for $\xi$ differs from the equation for $\psipz$ \eqref{eq:phiellp_0A} only in the constraint: $\psipz \in Z_0(\sigma)$, whereas $\xi \in Z_0(1)$. We directly observe that $\psipz = \xi  - \dashint_{\Omega} \sigma \xi$ --- the two solutions different only by a constant --- and hence also
\begin{align}
\phi = a_0(\xi,\xi) \; ; \label{eq:ps2}
\end{align}
the latter follows from our definition of $\phi$, \eqref{eq:phidef1}.

We may now write $\xi = \overline{\xi} + \hat{\xi}$, where $\overline{\xi} \in Z_0(1)$ satisfies
\begin{align}
a_0(\overline{\xi},v) = |\Omega|^{-1/2} L(v;1)\; , \forall v \in X \; , \label{eq:ps3}
\end{align}
and $\hat{\xi} \in Z_0(1)$ satisfies
\begin{align}
a_0(\hat{\xi},v) = |\Omega|^{-1/2} \gamma \int_\Omega (\sigma - 1)  v \; , \forall v \in X \; . \label{eq:ps4}
\end{align}
It is crucial to note that both \eqref{eq:ps3} and \eqref{eq:ps4} are solvable: the former thanks to the definition of $\gamma$ (and the unity-mean property of $\sigma$); the latter thanks again to the unity-mean property of $\sigma$. We further observe that $\overline{\xi} = \psipz$ for $\sigma = 1$, and hence 
\begin{align}
a_0(\overline{\xi},\overline{\xi}) = \phi(1) \; ; \label{eq:ps5} 
\end{align}
we recall that we consider here $\kappa = 1$ and fixed domain and thus $\phi$ is solely a function of $\sigma$

We now focus on \eqref{eq:ps4}. We note from the Cauchy-Schwarz inequality, the Poincar\'e-Wirtinger inequality, and the definition of $\mu$, \eqref{eq:mudef}, that 
\begin{align}
a_0(\hat{\xi},\hat{\xi}) &  \le  |\Omega|^{-1/2} \gamma \left(\int_\Omega (\sigma-1)^2\right)^{1/2} \left(\int_\Omega \hat{\xi}^2\right)^{1/2} \\[.1em]
& \le  |\Omega|^{-1/2} \gamma  \left(\int_\Omega (\sigma-1)^2\right)^{1/2} \left(\sup_{w \in Z_0(1)} \dfrac{ \left(\int_\Omega w^2\right)^{1/2}}
{ a_0^{1/2}(w,w)}\right) a_0^{1/2}(\hat{\xi},\hat{\xi}) \\[.1em]
& \le  |\Omega|^{-1/2} \gamma  \left(\int_\Omega (\sigma-1)^2\right)^{1/2} \mu^{-1/2} a_0^{1/2}(\hat{\xi},\hat{\xi}) \; , \label{eq:ps6}
\end{align}
and thus we obtain the bound
\begin{align}
a_0(\hat{\xi},\hat{\xi}) \le \gamma^2 \mu^{-1} \dashint_{\Omega} (\sigma-1)^2  \; . \label{eq:ps7}
\end{align}
Finally, we note that
\begin{align}
\phi(\sigma) = a_0(\overline{\xi} + \hat{\xi},\overline{\xi} + \hat{\xi})  \; ; \label{eq:ps8}
\end{align}
to conclude the proof, we then apply \eqref{eq:ps5}, \eqref{eq:ps7}, 
the Cauchy-Schwarz inequality, and the definition of $\delta$, \eqref{eq:deltadef}.

\subsection{Proof of Proposition \ref{prop:evfsigma}}\label{App:sketch_evfsigma}

We denote the value of $\sigma$ over the first domain by $\sigma_1$, and the value of $\sigma$ over the second domain by $\sigma_2$. To begin, we note that 
\begin{align}
 \dashint_{\Omega} (\sigma - 1)^2 & = \frac{1}{2}(\sigma_1-1)^2 + \frac{1}{2}(\sigma_2-1)^2  \; . \label{eq:evfs1}
 \end{align}
It then follows from \eqref{eq:evfs1} and our condition $\dashint_{\Omega}\sigma = 1$ --- hence $\sigma_1 + \sigma_2 = 2$ --- that
\begin{align}
 \dashint_{\Omega} (\sigma - 1)^2 & = (\sigma_1-1)^2 \; .
 \end{align}
We now note that $\sigma_1 > 0$ and thus $0 < \sigma_1 \le 2$. The maximum of $(b - 1)^2$ for $b \in [0,2]$ is 1, which concludes the proof.

\subsection{Proof of Proposition \ref{prop:phirt}}\label{App:sketch_phirt}

We shall assume that $\psipz$ takes the form
\begin{align}
\psipz(x_1,x_2) = b_{00} + b_{10}x_1 + b_{01}x_2 + b_{20}x_1^2 + b_{11}x_1 x_2 + b_{02} x_2^2 \;  \label{eq:pfphirt1}
\end{align}
for $b_{00}, b_{10}, b_{01}, b_{20}, b_{11}$, and $b_{02}$ real constants. We now impose five constraints: our equation for $\psipz$ --- the sensitivity equation in the interior of our triangle, \eqref{eq:xieqn_strong1}; the boundary condition \eqref{eq:xieqn_strong2} on each of the three edges of our triangle; and our normalization condition \eqref{eq:xieqn_strong3}. These conditions suffice to uniquely determine our six real constants. We should emphasize that the three constraints associated with the boundary condition must be satisfied at all points on the respective edges, and hence there is no {\it a priori} justification for our ansatz \eqref{eq:pfphirt1}. Once we obtain our coefficients, we evaluate $\phi$, from \eqref{eq:phidef1}, as
\begin{align}
\int_0^W\left(  \int_0^{1-\frac{x_1}{W}} \left(\frac{d\psipz}{dx_1}\right)^2 + \left(\frac{d\psipz}{dx_2}\right)^2 dx_2 \right) dx_1 \; \label{eq:pfphirt2}
\end{align}
note the integrand is polynomial in $(x_1,x_2)$.

We provide a few more details. The coefficient $b_{00}$ will not appear in \eqref{eq:pfphirt2}; furthermore, it can be shown that $b_{11} = 0$. The remaining coefficients are given by
\begin{align}
b_{10} = b_{01} = \sqrt{\dfrac{2}{W}} \;\; ; \quad b_{20} = b_{02} = -\sqrt{\dfrac{1}{2W}} \left( \dfrac{1 + W + \sqrt{1 + W^2}}{W} \right) \; .
\end{align}
As $W \rightarrow 0$ the terms $b_{02} = b_{20} \sim -2^{-1/2}W^{-3/2}$ dominate, and the integral \eqref{eq:pfphirt2} is then very easily evaluated.

For an isosceles right triangle, we assume the same solution form and substitute $W = 1$, which yields constants
\begin{align}
b_{10} = b_{01} = \sqrt{2} \;\; , \quad b_{20} = b_{02} = -1-\sqrt{2} \; .
\end{align}
For an equilateral right triangle, we define $\Omega$ as the union of $\Xi^\text{right triangle}(\sqrt{3})$ and its reflection about the $x_1$ axis. We assume the same solution form as before but now for the upper-half of the domain, and we also require $\frac{\partial\psi'^0}{dx_2}(x_1,0) = 0$ (from symmetry); the substitution $W = \sqrt{3}$ then yields constants
\begin{align}
b_{10} = 3^{-1/4}\;\; , \quad b_{01} = 0 \;\; , \quad b_{20} = b_{02} = -3^{1/4}/2 \;;
\end{align}
we then evaluate $\phi$ over the half-domain and double the result.

\subsection{Proof of Proposition \ref{prop:phitendom}}\label{App:sketch_phitendom}

We shall denote in this proof the measure of $\Omega^{\text{2D}}$ and $\partial\Omega^{\text{2D}}$ by $A$ and $P$, respectively. Then $|\Omega|$ and $|\partial\Omega|$ are given by $AL_3$ and $PL_3 + 2A$, respectively. We first recall our problem for $\psipz$ \eqref{eq:xieqn_strong1} -- \eqref{eq:xieqn_strong3} over the tensorized domain (recall we consider $\kappa = \sigma = 1$):
\begin{align} 
-\nabla^2 \psipz = (AL_3)^{-3/2}(PL_3 + 2A) \text{ in } \Omega \; , \label{eq:phitpf1}
\end{align}
with boundary condition
\begin{align}
\partial_n \psipz = -(AL_3)^{-1/2} \text{ on } \partial\Omega \; . \label{eq:phitpf2}
\end{align}
We then introduce corresponding 2D and 1D problems: find $\psipz_{\text{2D}}$ such that
\begin{align} 
-\nabla^2 \psipz_{\text{2D}} = (A)^{-3/2}(P) \text{ in } \Omega^{\text{2D}} \; , \label{eq:phitpf3}
\end{align}
with boundary condition
\begin{align}
\partial_n \psipz_{\text{2D}} = -(A)^{-1/2} \text{ on } \partial\Omega^{\text{2D}} \;  \label{eq:phitpf4}
\end{align}
and zero-mean normalization; 
 find $\psipz_{\text{1D}}$ such that
\begin{align} 
-\nabla^2 \psipz_{\text{1D}} = 2 (L_3)^{-3/2} \text{ in } \Omega^{\text{1D}} \; , \label{eq:phitpf5}
\end{align}
with boundary condition
\begin{align}
\partial_n \psipz_{\text{1D}} = -(L_3)^{-1/2} \text{ on } \partial\Omega^{\text{1D}} \; \label{eq:phitpf6}
\end{align}
and zero-mean normalization.  Note $\Omega^{\text{1D}} \equiv (0,L_3)$ and $\partial\Omega^{\text{1D}} \equiv \{0,L_3\}$. We note that the problems for $\psipz_{\text{2D}}$ and $\psipz_{\text{1D}}$ are solvable.
We further introduce
\begin{align}
\phi = \int_{\Omega} |\bnabla \psipz|^2\; ,\quad \phi_{\text{2D}} = \int_{\Omega_{\text{2D}}} |\bnabla_{\text{2D}} \psipz_{\text{2D}}|^2\;, \quad \phi_{\text{1D}} = \int_{\Omega_{\text{1D}}} |\bnabla_{\text{1D}} \psipz_{\text{1D}} |^2 \; , \label{eq:phitpf7} \
\end{align}
as provided by the definition \eqref{eq:phidef1}.  Note $\bnabla$, $\bnabla_{\text{2D}}$, and $\bnabla_{\text{1D}}$ refer to the gradient operator in $\RR^{d = 3}$, $\RR^2$, and $\RR^1$, respectively.

We now observe that
\begin{align}
\psipz = L_3^{-1/2} \psipz_{\text{2D}} + A^{-1/2} \psipz_{\text{1D}} \; , \label{eq:phitpf8}
\end{align}
and hence
\begin{align}
\phi & = L_3 \int_{\Omega_{\text{2D}}} L_3^{-1}  |\bnabla_{\text{2D}} \psipz_{\text{2D}}|^2 + A \int_{\Omega_{\text{1D}}} A^{-1} |\bnabla_{\text{1D}} \psipz_{\text{1D}}|^2 \\
& = \phi_{\text{2D}} + \phi_{\text{1D}} \; . 
\end{align} 
To conclude the proof we note from Table \ref{table:canonicalphi}  that $\phi_{\text{1D}} = 1/3$.

\subsection{Proof of Proposition \ref{prop:philb}}\label{App:sketch_philb}

We start with the weak form of the sensitivity equation, Equation \eqref{eq:phiellp_0A}. We know, by duality, that the $L^2$ norm of the function $\sqrt{k}\nabla\psipz$ may be written as
\begin{align}
\left\Vert\sqrt{\kappa}\,\nabla\psipz\right\Vert_{L^2(\Omega)} = \sup_{w\in(L^2(\Omega))^d}\frac{\left|\int_\Omega (\sqrt{\kappa}\,\nabla\psipz) \cdot w\right|}{\sqrt{\int_{\Omega}|w|^2}}.
\end{align}
Here, $d$ denotes spatial dimension (2 or 3) of $\Omega$. It follows, restricting the supremum to functions $w$ of the form $w=\sqrt{\kappa}\,\nabla v$ for $v\in H^1(\Omega)$, that
\begin{align}
\left\Vert\sqrt{\kappa}\,\nabla\psipz\right\Vert_{L^2(\Omega)}\ge \sup_{v\in H^1(\Omega)}\frac{\left|\int_\Omega (\sqrt{\kappa}\,\nabla\psipz)\cdot (\sqrt{\kappa}\,\nabla v)\right|}{\sqrt{\int_{\Omega}\kappa|\nabla v|^2}}.\label{eq:ressup}
\end{align}

Using now \eqref{eq:phiellp_0A} to evaluate the numerator of the right-hand side of \eqref{eq:ressup}, we obtain
\begin{align}
\left\Vert\sqrt{\kappa}\,\nabla\psipz\right\Vert_{L^2(\Omega)}\ge \sup_{v\in H^1(\Omega)}\frac{\left||\Omega|^{-1/2}\gamma\int_\Omega \sigma v - |\Omega|^{-1/2}\int_{\partial\Omega}v\right|}{\sqrt{\int_{\Omega}\kappa|\nabla v|^2}};\label{eq:ressup2}
\end{align}
which, given the definition \eqref{eq:phidef1} of $\phi$, yields the following bound:
\begin{align}
\phi(\sigma,\kappa,\Omega)\ge \left(\sup_{v\in H^1(\Omega)}\frac{\left||\Omega|^{-1/2}\gamma\int_\Omega \sigma v - |\Omega|^{-1/2}\int_{\partial\Omega}v\right|}{\sqrt{\int_{\Omega}\kappa|\nabla v|^2}}\right)^2.\label{eq:ressup3}
\end{align}

Given domain $\Omega$, we consider a largest possible ball, $B(x_0,R)\in\Omega$, where $x_0$ denotes its center and $R$ its radius. We are now going to choose, as test function for the supremum of the right-hand side, the function
\begin{align}
v^{\star}(x) =
\begin{cases}
R^2-|x-x_0|^2,\ x \in B(x_0,R)\\
0,\ x\in \Omega \setminus B(x_0,R)
\end{cases}
\end{align}
Thus, we obtain a bound, based on our test function which is zero-valued on and outside the ball $B(x_0,R)$,
\begin{align}
\phi(\sigma,\kappa,\Omega)\ge \frac{\left(|\Omega|^{-1/2}\gamma\int_{B(x_0,R)} \sigma v^\star\right)^2}{\int_{B(x_0,R)}\kappa|\nabla v^\star|^2}.\label{eq:ressup4}
\end{align}
We now consider the integral in the right-hand-side numerator of \eqref{eq:ressup4}. Noting that $v^\star$ and $\sigma$ are non-negative, we obtain the inequality
\begin{align}
\int_{B(x_0,R)}\sigma v^\star &\geq \underline{\sigma}\int_{B(x_0,R)}R^2-|x-x_0|^2 ,
\end{align}
where $\underline{\sigma} \equiv \essinf_{x\in B(x_0,R)}\sigma(x)$. Likewise, the integral in the denominator yields the inequality
\begin{align}
\int_{B(x_0,R)} \kappa |\nabla v^{\star}|^2 &\leq \overline{\kappa} \int_{B(x_0,r)}|\nabla(R^2-|x-x_0|^2)|^2,
\end{align}
where $\overline{\kappa}\equiv \sup_{x\in B(x_0,R)}\kappa(x)$.

We combine these results to simplify \eqref{eq:ressup4}:
\begin{align}
\phi(\sigma,\kappa,\Omega)\ge \frac{\underline{\sigma}^2}{\overline{\kappa}}\frac{|\partial\Omega|^2}{|\Omega|^3}\frac{\left(\int_{B(x_0,r)}R^2-|x-x_0|^2)\right)^2}{\int_{B(x_0,r)}|\nabla(R^2-|x-x_0|^2)|^2},
\end{align}
For $d=2$, the inequality reduces to
\begin{align}
\phi(\Omega^\text{2D})\ge\frac{\underline{\sigma}^2}{\overline{\kappa}} \frac{\pi}{8}\frac{|\partial\Omega^\text{2D}|^2}{|\Omega^\text{2D}|^3} R^4.
\end{align}
For $d=3$, the inequality reduces to
\begin{align}
\phi(\Omega^\text{3D})\ge \frac{\underline{\sigma}^2}{\overline{\kappa}}\frac{4\pi}{45}\frac{|\partial\Omega^\text{3D}|^2}{|\Omega^\text{3D}|^3} R^5.
\end{align}
This concludes the proof.

\subsection{Proof of Proposition \ref{prop:l1ltBg}}\label{App:sketch_l1ltBg}

\begin{proposition}[Rayleigh Quotient: Uniform Candidate Upper Bound] \label{prop:l1ltBg}
For any problem data $\{\sigma,\kappa,B,\Omega\}$,
\begin{align}
\lambda(B) \le \gamma B \label{eq:lam1gtr}
\end{align}
for $\gamma$ given by \eqref{eq:gammadef}. Note there are no additional restrictions on the properties and $\Omega$. Furthermore, for $B>0$, we obtain a strict inequality.
\end{proposition}
We first note that from \eqref{eq:RQB},
\begin{align}
\calR(1;\sigma,\kappa,B) = \frac{B a_1(1,1)}{m(1,1;\sigma)} = B \frac{\dOmMeas}{\OmMeas} \; ,\label{eq:l1Ray1}
\end{align}
since $\sigma$ is a unity domain mean function. Thus
\begin{align}
\gamma B = \calR(1;\sigma,\kappa,B) \ge \min_{w \in X}\calR(w;\sigma,\kappa,B) = \lambda_1(B),
\end{align}
We note that for $B>0$, $1$ does not satisfy the boundary condition and therefore cannot be the minimizer of the Rayleigh quotient. This concludes the proof.

\subsection{Proof of Lemma \ref{lem:UDeltaL}}\label{App:UDeltaL}

First we consider the difference between the series representations for $u_\avg$ \eqref{eq:sov_avgQoI} and $u_\paavg$ \eqref{eq:sov_paavgQoI}:
\begin{align}
u_\avg(t)-u_\paavg(t) &= \sum_{k=1}^\infty c_k\, e^{-\lambda_k t},\label{eqn:udeltack}\\
c_k&\equiv\OmMeas(M(\psi_k))^2\left(1-\frac{\lambda_k}{B\gamma}\right).\label{eqn:ckdef}
\end{align}
Because $\lambda_1 < B\gamma$ (Proposition \ref{prop:l1ltBg}), there exists $n>1$ such that
\begin{align}
\lambda_k < B\gamma,\ \forall\ k < n,\label{eqn:lamklen}\\
\lambda_k \ge B\gamma,\ \forall\ k \ge n;\label{eqn:lamkgtn}
\end{align}
therefore,
\begin{align}
c_k > 0,\ \forall\ k < n,\label{eqn:cklen}\\
c_k \le 0,\ \forall\ k \ge n.\label{eqn:ckgtn}
\end{align}
We note from the prescribed uniform initial condition \eqref{eq:ivpstrong_3} that
\begin{align}
u_\avg(0) - \u_\paavg(0) = (1) - (1) = 0\label{eqn:u0equal}.
\end{align}
With that relation and evaluation of \eqref{eqn:udeltack} at $t=0$, we have
\begin{align}
\sum_{k=1}^{n-1} c_k = -\sum_{k=n}^\infty c_k.\label{eqn:ckequal}
\end{align}
We now consider henceforth $t>0$ only.  From \eqref{eqn:lamklen} -- \eqref{eqn:ckgtn}, we obtain inequalities
\begin{align}
\sum_{k=1}^{n-1} c_k \left(e^{-\lambda_k t} -e^{-\lambda_{n-1} t}\right) &\ge 0,\\
\sum_{k=n}^\infty c_k \left(e^{-\lambda_k t} -e^{-\lambda_{n} t}\right) &\ge 0;
\end{align}
recognizing that each term in both series is non-negative. Summing these two inequalities while using \eqref{eqn:udeltack}, \eqref{eqn:lamklen}, \eqref{eqn:lamkgtn}, \eqref{eqn:cklen}, and \eqref{eqn:ckequal}; we show the positivity of the averages difference
\begin{align}
u_\avg(t)-u_\paavg(t) \ge \left(\sum_{k=1}^{n-1} c_k\right) \left(e^{-\lambda_{n-1} t}-e^{-\lambda_{n} t}\right) > 0\,.
\end{align}
Finally, dividing by a strictly positive quantity $u_\avg(t)$ (all terms in \eqref{eq:sov_avgQoI} are positive) we obtain \eqref{eq:UDeltaL}.  \eqref{eq:UDeltaLT0} follows from \eqref{eqn:u0equal} and \eqref{eq:UDeltaLB0} follows from the solution of the Neumann problem, $u(\cdot,t;B=0)=\OmMeas^{1/2}\psiz(\cdot)$ --- from inspection. This concludes the proof.

\subsection{Proof of Proposition \ref{prop:UavgBound}}\label{App:UavgBound}

We start by introducing a thermal energy equality based on the slow time-scale formulation:
\begin{lemma}[Thermal Energy Equality]\label{lem:TE1}
The weak form \eqref{eq:ivpweak_1} for $u(\cdot,t)$ transforms to the weak form for $U(\cdot,t)$:
\begin{align}
B \gamma\, m(\partial_T U,v;\sigma) + a(u,v;B,\kappa) = 0\ \forall v \in X,\label{eq:ivpweak_U}
\end{align}
subject to $U(\cdot,0) = 1$.
The weak form in turn yields the thermal energy equality,
\begin{align}
\frac{d}{dT} U_\avg =   - U_\paavg.\label{eqn:Uevo}
\end{align}
\end{lemma}

Proof: The weak form follows from standard change of variables arguments. To obtain the thermal energy equality we take $v = 1 \in X$ in the weak form.

We now turn to the proof of Proposition \ref{prop:UavgBound}. The thermal energy equality \eqref{eqn:Uevo} is our point of departure:
\begin{align}
\frac{d}{dT} U_\avg =  -U_\avg + \left(U_\avg - U_\paavg\right).
\end{align}
Now specifying an integration factor and invoking the definition of $L(v;\sigma)$, \eqref{eq:Ldef}, to obtain
\begin{align}
\frac{d}{dT}\left(e^T U_\avg\right) = e^T |\partial \Omega|^{-1}L(U;\sigma).
\end{align}
We then integrate over time and apply the initial condition $U(\cdot,0) = 1$ to obtain the time evolution of the average temperature
\begin{align}
U_\avg(T) =  e^{-T} + |\partial\Omega|^{-1}\int_0^T e^{-(T-T')}L(U(\cdot,T');\sigma)\,dT'.\label{eqn:Uavg_int}
\end{align}
From Lemma \ref{lem:UDeltaL}, we have $U_\avg > U_\paavg$ and hence $L(U; \sigma) = \dOmMeas (U_\avg - U_\paavg)>0$ for all $T > 0$; thus, the integral in \eqref{eqn:Uavg_int} is positive for $T>0$. For $T=0$, the prescribed uniform initial condition leads to $U_\avg(0) = \UavgL{1}(0) = 1$. For $B=0$, $T\equiv B\gamma\,t = 0$ which concludes the proof.

\subsection{Proof of Proposition \ref{Prop:N1}}\label{App:PropN1}

This proof, when restricted to a spherical domain and homogeneous material properties, is equivalent to the error analysis found in \cite{ODD}.  We first start with an elementary lemma for comparing two exponentials:
\begin{lemma}\label{lem:N1a}
Let $g(z;\epsilon) \equiv \exp(-z(1-\epsilon)) - \exp(-z)$ for $z \in [0,\infty)$  and $\epsilon \in (-1,1)$. Then, $|g(z;\epsilon)|$ attains a unique supremum at
\begin{align}
z_{\sup}(\epsilon) \equiv \argsup{z\,\in\, [0,\infty)} |g(z;\epsilon)| = \dfrac{-\ln(1-\epsilon)}{\epsilon},
\end{align}
and furthermore
\begin{align}
g_{\sup}(\epsilon) \equiv |g(z_{\sup};\epsilon)| \le \dfrac{|\epsilon|}{e} + \pcalO(\epsilon^2) \text{ as } \epsilon \rightarrow 0 .
\end{align}
Note here $e = \exp(1)$.
\end{lemma}

Proof of Lemma \ref{lem:N1a}: We consider both cases of $\epsilon > 0$ and $\epsilon < 0$ (we note $g_\text{sup}(0)=0$). In this proof $g'$ and $g''$ shall denote the first and second derivatives of $g$, respectively.

We first note from observation that $g(z;\epsilon) > 0$ for $z,\epsilon > 0$ and $g(z;\epsilon) < 0$ for $z > 0$ and $\epsilon < 0$; we may then take $|g(z)| = g(z)$ for $\epsilon>0$ or $|g(z)| = -g(z)$ for $\epsilon < 0$, which admits all derivatives. The first two derivatives of $g(z;\epsilon)$ are
\begin{align}
g'(z;\epsilon) &= -(1-\epsilon) \exp(-z(1-\epsilon)) + \exp(-z),\label{eq:gp1}\\
g''(z;\epsilon)&=\exp(-z)((1-\epsilon)^2 \exp(z\epsilon) - 1).\label{eq:gp2}
\end{align}
We first consider \eqref{eq:gp1} and conclude that $g'(z_0;\epsilon) = 0$ has unique solution
\begin{align}
z_0 =  \dfrac{-\ln(1-\epsilon)}{\epsilon} .
\end{align}
We then evaluate
\begin{align}
g''(z_0;\epsilon) = -\epsilon \exp(-z_{\sup}),
\end{align}
and conclude $|g(z_0;\epsilon)|$ is maximal for positive and negative $\epsilon$ since $g(0) = 0$ and $g(z\rightarrow\infty)\rightarrow 0$. Finally, we evaluate $g$ at $z_{\sup}$, expanding in terms of perturbation $\epsilon$:
\begin{align}
g(z_{\sup};\epsilon) &=  \exp(-z_{\sup})(\exp(\epsilon z_{\sup})-1)\nonumber\\
&= \exp\left(\dfrac{-\epsilon+\pcalO(\epsilon^2)}{\epsilon} \right)\left(\exp\left(\epsilon+\pcalO(\epsilon^2)\right)-1\right)  \nonumber \\
&= \exp\left(-1  \right)\left(1+\pcalO(\epsilon)\right)  \left(\epsilon+\pcalO(\epsilon^2)\right)\nonumber \\
&= \dfrac{\epsilon}{e} + \pcalO(\epsilon^2).
\end{align}
We also note that
\begin{align}
z_{\sup} = 1 + \epsilon/2 + \pcalO(\epsilon^2),\label{eq:zmaxasymp}
\end{align}
and hence $z_{\sup} \rightarrow 1$ as $\epsilon \rightarrow 0$, which concludes the proof of Lemma \ref{lem:N1a}.

Now we introduce another elementary lemma:
\begin{lemma}\label{lem:N1b}
We may express $u_{\avg}(t;B) - \uavgL{1}(t;B)$ as
\begin{align}
u_{\avg}(t;B) - \uavgL{1}(t;B) =  \sum_{k=1}^\infty \OmMeas (M(\psi_k(B)))^2 y_k(t;B)\; ,\label{eq:ebarsv2}
\end{align}
for
\begin{align}
y_k(t;B) \equiv  \exp(-\lambda_k(B) t) - \exp(-B\gamma t)\; , \forall k \in \NN \; .\label{eq:Ydef2}
\end{align}
Furthermore,
\begin{align} 
y_1(t;B) \ge 0, \forall t \in \calT \; ,  \label{eq:Y1gtr02}
\end{align}
and
\begin{align}
|y_k(t;B)| \le 1, \forall t \in \calT, \forall k \in \NN \; . \label{eq:Yle12}
\end{align}
The eigenpairs $(\psi_k(B),\lambda_k(B))$ are defined in Appendix \ref{App:evp}.
\end{lemma}

Proof of Lemma \ref{lem:N1b}: We first demonstrate \eqref{eq:ebarsv2}. We recall \eqref{eq:sov_sum} and \eqref{eq:u1avg} and invoke \eqref{eq:Mpsisquare}:
\begin{align}
u_{\avg}(t;B) - \uavgL{1}(t;B) &= \left( \sum_{k=1}^\infty \OmMeas (M(\psi_k(B)))^2 \exp(-\lambda_k(x;B)t) \right) - \uavgL{1}(t;B) \\
&= \sum_{k=1}^\infty \left( \OmMeas (M(\psi_k(B)))^2 (\exp(-\lambda_k(x;B)t - \uavgL{1}(t;B)) \right)\\
&= \sum_{k=1}^\infty \OmMeas (M(\psi_k(B)))^2 y_k(t;B)\hspace{2.28in}\
\end{align}
for $y_k(t,B), k \in \NN$, given by \eqref{eq:Ydef2}. We next consider \eqref{eq:Y1gtr02}: we note that $y_1(t;B) \equiv \exp(-\lambda_1(B)t) - \exp(-B\gamma t)$; however, from Proposition \ref{prop:l1ltBg}, $B\gamma \ge \lambda(B)$, and hence $y_1 \ge 0$. Finally, we consider \eqref{eq:Yle12}: the result follows from the monotonicity of $\exp(-z)$ as $z$ increases. This concludes this proof for Lemma \ref{lem:N1b}.

We can now turn to the proof of Proposition \ref{Prop:N1}. We first appeal to Lemma \ref{lem:N1a}, \ref{lem:eigfcnamp} and \eqref{eq:ebarsv2} to obtain
\begin{align}
u_\avg(t;B) - \uavgL{1}(t;B) &=   (1 - B^2\Upsilon) y_1(t;B) +\nonumber \\
&\sum_{k=2}^\infty \OmMeas (M(\psi_k(B)))^2 y_k(t;B) + \pcalO(B^3)\; .\label{eq:ebarsv_Z2}
\end{align}
We then note from \eqref{eq:Mpsisquare}, \eqref{eq:eigamp}, and \eqref{eq:Yle12} that
\begin{align}
\left| \; \sum_{k=2}^\infty \OmMeas (M(\psi_k(B)))^2 y_k(t;B) \; \right|  \le \Upsilon B^2 + \pcalO(B^3)\; , \label{eq:sumbound2}
\end{align}
which yields
\begin{align}
\left|u_\avg(t;B) -\uavgL{1}(t;B)\right| \le \left|1-B^2\Upsilon\right|y_1(t;B) + \Upsilon B^2 +  \pcalO(B^3) \; . \label{eq:ebarA12}
\end{align}
We now introduce

\begin{align}
Y_1(T;B) \equiv y_1\left(\dfrac{T}{B\gamma}\right),
\end{align}
in terms of which we may write
\begin{align}
\left|U_\avg(T;B) -\UavgL{1}(T;B)\right| &\le \left|1-B^2\Upsilon\right||Y(T;B)| + \Upsilon B^2 +  \pcalO(B^3),\label{eq:U1avgLE_full}\\
&\le |Y(T;B)| + \pcalO(B^2)\label{eq:U1avgLE_abbrev},
\end{align}
where we have abbreviated $Y_1$ simply as $Y$. Thus in all subsequent analysis we may consider just $Y(T;B)$, as the remaining contributions to the error are $\pcalO(B^2)$.

We may write $Y(T;B)$ as
\begin{align}
Y(T;B) = \exp\left(-\dfrac{\lambda(B)}{B\gamma}T\right)  - \exp(-T),\label{eq:YTB0}
\end{align}
where we recall that $\lambda(B) \equiv \lambda_1(B)$. We now define
\begin{align}
\epsilon \equiv 1 - \dfrac{\lambda(B)}{B\gamma}\label{eq:deltadef2}
\end{align}
and note that
\begin{align}
\epsilon= \phi B / \gamma + \pcalO(B^2),\label{eq:deltasub}
\end{align}
due to Proposition \ref{prop:eigvasymp}. We now substitute the $\epsilon$ definition into \eqref{eq:YTB0} to obtain
\begin{align}
Y(T;B) = \exp(-T(1-\epsilon)) - \exp(-T) .\label{eq:YTB}
\end{align}
We now directly apply Lemma \ref{lem:N1a} and \eqref{eq:deltasub} to obtain
\begin{align}
\max_{T\in[0,\Tf]} |Y(T;B)|\le\max_{T \in [0,\infty)} |Y(T;B)| = \frac{\phi B}{\gamma\,e} + \pcalO(B^2).
\end{align}
Substitution of this inequality into the the maximal evaluation of the inequality in \eqref{eq:U1avgLE_abbrev}, we obtain \eqref{eq:easymp}.
Finally, we obtain the $\argmax$ result, $\eqref{eq:Tmax}$, from substitution of $\epsilon$ into \eqref{eq:zmaxasymp}. This concludes the proof.

\subsection{Proof of Proposition \ref{Prop:N2}}\label{App:PropN2}

We introduce the quadratic energy equality in addition to the thermal energy equality (Lemma \ref{lem:TE1}) to bound the error in $\UavgL{1}$.

\begin{lemma}[Quadratic Energy Equality]\label{lemTE}
The weak form in \eqref{eq:ivpweak_U} yields the “quadratic” energy equality,
\begin{align}
B\gamma \frac{1}{2}( m(U,U) - |\Omega|) &+ \int_0^T a_0(U(T’),U(T’))dT’\nonumber\\&+ B \int_0^T a_1(U(T’),U(T’))dT’= 0\label{eqn:quadraticenergy}\,.
\end{align}
\end{lemma}

Proof: Starting from \eqref{eq:ivpweak_U}, we take $v = U \in X$ in the weak form, integrate in time $T’$ from $0$ to $T$, and invoke the initial condition. $\Box$

We now turn to the proof of Proposition \ref{Prop:N2}. Note, all inequalities assume $T\ge0$. The thermal energy equality \eqref{eqn:Uavg_int} is our point of departure. We let $v=U$ in \eqref{eq:phiellp_0A}, then invoke the Cauchy-Schwarz inequality, and apply \eqref{eq:phiellp_a} to deduce that
\begin{align}
|\Omega|^{-1/2}L(U)\le \phi^{1/2}a_0^{1/2}(U,U).
\end{align}
Next we invoke Cauchy-Schwarz inequality to bound the time integral as
\begin{align}
&|\Omega|^{-1/2}\int_0^T e^{-(T-T')}L(U)\,dT'\nonumber\\
&\quad\le \phi^{1/2}\left(\int_0^T e^{-2(T-T')}\,dT'\right)^{1/2}\left(\int_0^T a_0\left(U(T'),U(T')\right)\,dT'\right)^{1/2}.\label{eqn:cstime}
\end{align}
The first integral on the right hand side is easily evaluated as
\begin{align}
\int_0^T e^{-2(T-T')}\,dT' = \frac{1 - e^{-2T}}{2}\label{eqn:Tint}\,,
\end{align}
and furthermore
\begin{align}
\left|\frac{1 - e^{-2T}}{2}\right| \le \frac{1}{2}\;\label{eqn:Tint};
\end{align}
the second integral can be bound based on the energy balance \eqref{eqn:quadraticenergy},
\begin{align}
\int_0^T a_0\left(U(T'),U(T')\right)\,dT' \le B\gamma\frac{|\Omega|}{2}\label{eqn:Uavg_ineq}\,.
\end{align}

Assembling \eqref{eqn:Uavg_int}, \eqref{eqn:cstime}, \eqref{eqn:Tint}, and \eqref{eqn:Uavg_ineq}; we obtain
\begin{align}
\left|U_\avg(T)-e^{-T}\right|\le \frac 1 2 \left(\frac{\phi B}{\gamma}\right)^{1/2},\ \forall\, T\, \ge 0\,.\label{eqn:Uavg1_err}
\end{align}

The relative error bound is obtained through multiplication of \eqref{eqn:Uavg1_err} by $\exp(T)$.
\subsection{Proof of Lemma \ref{lem:OB3}}\label{App:OB3}

Compared to the procedure for finding $\lambdappz(=-2\phi)$ (Appendix \ref{App:sketch_eigvasymp}), we opt for an explicit perturbation approach to solve for $\lambdapppz$. To start, we consider the equations for $\psi(\equiv\psi_1)$ and $\lambda(\equiv\lambda_1))$:
\begin{align}
m(\psi,\psi) &=1,\label{eqn:psinorm}\\
\lambda\, m(\psi,\psi) &= a_0(\psi,\psi)+Ba_1(\psi,\psi),\label{eqn:psi2}\\
\lambda\, m(\psi,v) &= a_0(\psi,v)+Ba_1(\psi,v)\label{eqn:psi1}.
\end{align}
We then expand $(\psi,\lambda)$ around $B=0$: $\psi=\sum_{k=0}^\infty B^k[(\frac{d}{dB})^k\psi](0)/k!$ and $\lambda=\sum_{k=0}^\infty B^k[(\frac{d}{dB})^k\lambda](0)/k!$. We then substitute these expansions into the equations above to obtain equations containing $\pcalO(1)$, $\pcalO(B)$, $\pcalO(B^2)$, ... which are then grouped ordered from the low to high order in $B$. Thus, we start with $\pcalO(1)$ terms:
\begin{align}
m(\psi^0,\psi^0) &=1,\label{eqn:B0_1}\\
\lambdaz m(\psi^0,\psi^0) &=a_0(\psiz,\psiz),\\
\lambdaz m(\psiz,v) &= a_0(\psiz,v),\ \forall v\in\, X.
\end{align}
By inspection $\psiz=\OmMeas^{-1/2}$ and therefore 
\begin{align}
\lambdaz &= 0,\label{eqn:B0_2}\\
a_0(\psiz,v) &= 0,\ \forall v\in\, X.\label{eqn:B0_3}
\end{align}
Now proceeding to $\pcalO(B)$ terms using \eqref{eqn:B0_2} and \eqref{eqn:B0_3}, we start with
\begin{align}
2m(B\psi'^0,\psiz) &=0,\\
B\lambdapz m(\psiz,\psiz) &= B a_1(\psiz,\psiz),\\
B\lambdapz m(\psiz,v) &= a_0(B\psipz,v) + Ba_1(\psiz,v),\ \forall v\in\, X.
\end{align}
Simplifying the terms with \eqref{eqn:B0_1} (and recognizing that $a_1(\psiz,\psiz)=\gamma$), we obtain
\begin{align}
m(\psi'^0,\psiz) &=0,\label{eqn:B1_1}\\
\lambda'^0 &= \gamma,\label{eqn:B1_2}\\
a_0(\psi'^0,v) &= \gamma\, m(\psiz,v) -a_1(\psiz,v),\ \forall v\in\, X.\label{eqn:B1_3}
\end{align}
We note that we re-identify the equations for $\psipz$ and $\lambdapz$ obtained earlier. We now proceed to $\pcalO(B^2)$ terms using \eqref{eqn:B0_2}, \eqref{eqn:B1_1}, and \eqref{eqn:B0_3}; we start with
\begin{align}
m(B\psi'^0,B\psi'^0) + 2m\left(\frac{B^2}{2}\psi''^0,\psiz\right) &=0,\\
\frac{\lambda''^0}{2}B^2m(\psiz,\psiz) &= a_0(B\psipz,B\psipz)+2Ba_1(B\psipz,\psiz),\\
\frac{\lambda''^0}{2}B^2m(\psiz,v) + B\gamma\, m(B\psipz,v)&= a_0\left(\frac{B^2}{2}\psippz,v\right)+Ba_1(B\psipz,v),\ \forall v\in\, X.
\end{align}
Simplifying using \eqref{eqn:B1_1}, \eqref{eqn:B1_3} with $v=\psipz$, and recalling that $\phi\equiv a_0(\psipz,\psipz)$; we obtain
\begin{align}
m(\psi'^0,\psi'^0) + m(\psi''^0,\psiz) &=0,\label{eqn:B2_1}\\
-\frac{\lambda''^0}{2} &= \phi,\label{eqn:B2_2}\\
\frac{1}{2}a_0(\psi''^0,v) &= \gamma\, m(\psi'^0,v) - a_1(\psi'^0,v)-\phi\, m(\psi^0,v),\ \forall v\in\, X.\label{eqn:B2_3}
\end{align}
We note the re-identification of $\lambdappz$. Now proceeding to $\pcalO(B^3)$ terms (note we only consider \eqref{eqn:psi2}) using \eqref{eqn:B1_1}, \eqref{eqn:B2_1}, and \eqref{eqn:B0_3}; we start with
\begin{align}
\frac{\lambda'''^0}{6}B^3m(\psiz,\psiz) &= 2 a_0 \left(\frac{B^2}{2}\psippz,B\psipz\right),\text{ hence} \nonumber\\
& + Ba_1(B\psi'^0,B\psi'^0)\nonumber\\
& + 2 B a_1 \left(\frac{B^2}{2}\psippz,\psiz\right),
\end{align}
which simplifies, using \eqref{eqn:B0_1}, \eqref{eqn:B1_3} with $v=\psippz$, and \eqref{eqn:B2_1}, to
\begin{align}
\frac{\lambda'''^0}{6} &= \chi-\gamma\, \Upsilon.\label{eqn:B3_2}
\end{align}
We finally show compatibility for \eqref{eqn:B1_3} and \eqref{eqn:B2_3} for constant test functions, say $v=\psiz$. For \eqref{eqn:B1_3}, using \eqref{eqn:B0_1} and \eqref{eqn:B0_3}, we have
\begin{align}
a_0(\psipz,\psiz) &= \gamma\, m(\psiz,\psiz) - a_1(\psiz,\psiz)\\
0 &= \gamma - \gamma;
\end{align}
and for \eqref{eqn:B2_3}, using \eqref{eqn:B0_1}, \eqref{eqn:B0_3}, \eqref{eqn:B1_1}, and \eqref{eqn:B1_3} with $v=\psiz$; we have
\begin{align}
\frac{1}{2}a_0(\psiz,\psi''^0) &= \gamma\, m(\psiz,\psi'^0) - a_1(\psiz,\psi'^0)-\phi\, m(\psiz,\psi^0),\text{ hence}\\
0 &= 0 - (-\phi) - \phi.
\end{align}
This concludes the proof.

\subsection{Proof of Proposition \ref{prop:chiscale}} \label{App:sketch_chiscale}

We start the proof from an intermediate result from Appendix \ref{App:sketch_phiprops}: combining the result from  \eqref{eq:scaledeq} and \eqref{eq:scaledbc} to establish the scaling relation of $\psipz$ w.r.t. dilation factor $\beta$,
\begin{align}
  \beta^{d/2-1}\hat{\xi}_2 = \xi_1 .
\end{align}
recall that $\xi_1(\cdot)\equiv \psipz(\cdot;\Omega_1)$, $\xi_2\equiv \psipz(\cdot;\Omega_2)$, and $\hat{\xi}_2(x)= \xi_2(\beta x),\ \forall x \in \Omega_1$. We now proceed similarly with a different set of quadratic functionals on $\psipz$; we first consider $\Upsilon$ defined in \eqref{eq:Upsilondef}:
\begin{align}
  \Upsilon_2 &= \int_{\Omega_2} \sigma\, \xi_2\, \xi_2\nonumber\\
  &= \beta^d \int_{\Omega_1} \sigma\, (\beta^{1-d/2}\xi_1)^2\nonumber\\
  &= \beta^2 \Upsilon_1.\label{eq:Upsilonscale}
\end{align}
For the evaluation of $\chi$ defined in \eqref{eqn:chidef}, we note that the scaling of the domain-boundary measure is $|\partial\Omega_2| = \beta^{d-1}|\partial\Omega_1|$. Following an analogous procedure, we obtain
\begin{align}
  \chi_2 &= \int_{\partial\Omega_2} \xi_2\, \xi_2\nonumber\\
  &= \beta^{d-1} \int_{\partial\Omega_1} (\beta^{1-d/2}\xi_1)^2\nonumber\\
  &= \beta\, \chi_1.\label{eq:chiscale}
\end{align}
Substitution of $\beta=\gamma_1/\gamma_2$ from \eqref{eq:gammascale} into \eqref{eq:Upsilonscale} and \eqref{eq:chiscale} concludes the proof.

\subsection{Proof of Proposition \ref{prop:chiupsilontendom}}\label{App:sketch_chiupsilontendom}
Proceeding from Equation \eqref{eq:phitpf8} from the proof of Proposition \ref{prop:phitendom}, we evaluate $\Upsilon$ as
\begin{align}
\Upsilon &= \int_{\Omega} \left( L_3^{-1} (\psipz_{\text{2D}})^2 + 2 L_3^{-1/2}A^{-1/2} \psipz_\text{2D}\psipz_\text{1D}+  A^{-1} (\psipz_{\text{1D}})^2 \right)\\
&= \Upsilon_{\text{2D}} + 2 L_3^{1/2}A^{1/2} \left(\int_{\Omega_\text{2D}} \psipz_{\text{2D}}\int_{\Omega_\text{1D}} \psipz_{\text{1D}}\right) + \Upsilon_{\text{1D}}.
\end{align}
We recognize that both $\int_{\Omega_\text{2D}} \psipz_{\text{2D}}$ and $\int_{\Omega_\text{1D}}\psipz_\text{1D}$ evaluate to zero from the zero-mean condition of Proposition \ref{prop:EEforSensitivityDerivative}; thus we obtain \eqref{eq:upsilontendom1}.

For $\chi$, we follow the same steps, however we must recognize that for the tensorized domain, we must consider integration on both the lateral sides, $\int_{\Omega^\text{1D}}\int_{\partial\Omega^\text{2D}}$, and the two faces, $\int_{\Omega^\text{2D}}\int_{\partial\Omega^\text{1D}}$, for the boundary integration:
\begin{align}
\chi &= \int_{\partial\Omega} \left( L_3^{-1} (\psipz_{\text{2D}})^2 + 2 L_3^{-1/2}A^{-1/2} \psipz_\text{2D}\psipz_\text{1D}+  A^{-1} (\psipz_{\text{1D}})^2 \right)\\
&= \chi_\text{2D} + \gamma_\text{1D}\Upsilon_{\text{2D}} \nonumber\\
&+2 L_3^{-1/2}A^{-1/2} \left(\int_{\partial\Omega^\text{2D}} \psipz_{\text{2D}}\,\int_{\Omega_\text{1D}}\psipz_\text{1D}+\int_{\partial\Omega^\text{1D}} \psipz_{\text{1D}}\,\int_{\Omega^\text{2D}}\psipz_\text{2D}\right) \nonumber\\
&+ \gamma_\text{2D}\Upsilon_{\text{1D}} + \chi_\text{1D} .
\end{align}
We once again recognize that both $\int_{\Omega_\text{2D}} \psipz_{\text{2D}}$ and $\int_{\Omega_\text{1D}}\psipz_\text{1D}$ evaluate to zero from the zero-mean condition of Proposition \ref{prop:EEforSensitivityDerivative}; thus we obtain \eqref{eq:chitendom1}.

\subsection{Proof of Proposition \ref{Prop:N3p}}\label{App:PropN3p}

The proof of Proposition \ref{Prop:N3p} requires only small modifications to address the second-order approximation of Proposition \ref{Prop:N2}. First, in Lemma \ref{lem:N1b} we replace coefficients $y_k$ with
\begin{align}
\hat{y}_k(t) \equiv \exp(-\lambda_k(B)t) - \exp\left(\dfrac{-B\gamma}{1 + \phi B/\gamma}\right).
\end{align}
Next, in Proposition \ref{Prop:N1} proper, we introduce $\hat{T} \equiv T/(1 + \Bdunkp)$, in terms of which we may write
\begin{align}
\hat{Y}(\hat{T}) = \exp\left(-\lambda(B)\dfrac{1 + \Bdunkp}{B\gamma}\hat{T}\right) -\exp(\hat{T}),
\end{align}
and subsequently
\begin{align}
\hat{Y}(\hat{T}) = \exp(-\hat{T}(1-\hat{\epsilon})) - \exp(\hat{T})
\end{align}
where $\hat{\epsilon}$, previously defined by \eqref{eq:deltadef2}, is now given by
\begin{align}
\hat{\epsilon} \equiv 1 - \lambda(B)\dfrac{1 + \Bdunkp}{B\gamma} .\label{eq:deltahatdef}
\end{align}
Finally, expansion in $B$ yields $\hat{\epsilon} = \frac{\phi^2-\gamma\chi+\gamma^2\Upsilon}{\gamma^2}B^2 +\pcalO(B^3)$. To conclude the proof, we now apply \eqref{eq:ebarsv2} and note that $B^2\hat{\epsilon}$ is $\pcalO(B^4)$.

\subsection{Proof of Proposition \ref{prop:Deltalumpederror}}\label{App:sketch_Deltalumpederror}

We start with two elementary lemmas:
\begin{lemma}\label{lem:EdeltaL_A}
For $B \rightarrow 0$,
\begin{align}
\left| 1 - \frac{\lambda_1(B)}{B \gamma} - \phi\frac{B}{\gamma}\left(1+\phi\frac{B}{\gamma}\right)^{-1}\right| = |\gamma\chi-\gamma^2\,\Upsilon-\phi^2|\frac{B^2}{\gamma^2} + O(B^3) \;  \label{eq:LemmaB}
\end{align}
where $\phi$, $\chi$ and $\Upsilon$ are defined in \eqref{eq:phidef1}, \eqref{eqn:chidef} and \eqref{eq:Upsilondef}; respectively.
\end{lemma}
The proof of result \eqref{eq:LemmaB} follows from Lemma \ref{lem:OB3}: rearranging the terms, dividing by $B\gamma$, expanding the rational polynomial expression for $\UDeltaL$ about $B=0$, and taking the absolute value yields the desired result. We now introduce the second lemma:
\begin{lemma}\label{lem:EdeltaL_B}
Let $g(\nu,z) \equiv \nu \exp(-\nu z)$. Then for any $\nu \in \Rzp$, and any $z \ge z_0$ for $z_0 \in \Rp$, $g(\nu,z) \le \exp(-1)/z_0$. \end{lemma}
Proof: To demonstrate this result, we note that for a given non-negative $z$, $g(\nu,z)$ attains a global maximum of $\exp(-1)/z$ at $\nu_{\max} = 1/z$; thus for any non-negative $\nu$, and $z \ge z_0$, $g(\nu,z) \le \exp(-1)/z_0$.

We now proceed to the main proof. It shall prove useful to introduce
\begin{align}
A_1(T) &\equiv |\Omega| M^2(\psi_1) \exp(-\lambda_1 T/(B\gamma)), \; \label{eq:A1Delta}\\
\EDeltaL(T;T_0,B) &\equiv \max_{T'\in[T_0,T]}| U_\Delta(T';B) - \UDeltaL(B)|. \label{eq:eDeltadef}
\end{align}
We first note from \eqref{eq:udelta},  \eqref{eq:sov_avgQoI}, \eqref{eq:sov_paavgQoI}, and \eqref{eq:eDeltadef} that
\begin{align}
\EDeltaL(T;T_0,B) = \max_{T'\in[T_0,T]}\displaystyle
\frac{\left|\, \displaystyle\sum_{k = 1}^\infty |\Omega| M^2(\psi_k) \left( 1 - \frac{\lambda_k}{B\gamma} -\UDeltaL(B) \right) \exp(-\lambda_k T'/(B\gamma))  \,\right|}
{\displaystyle \sum_{k = 1}^\infty |\Omega| M^2(\psi_k) \exp(-\lambda_k T'/(B\gamma)) } \; ;\label{eq:A2Delta}
\end{align}
note we may bring the $\UDeltaL$ term into the sum thanks to the denominator.
We now take \eqref{eq:A2Delta} as our point of departure: we bring the absolute value in the numerator inside the sum --- to provide and upper bound for the numerator; we next break the sum in the numerator into two parts, $k = 1$, and $k > 1$; we next neglect  the (positive) summands for $k > 1$ in the denominator --- to provide a lower bound for the denominator; finally, we divide through (numerator and denominator) by $A_1(T)$. We thereby obtain
\begin{align}
\EDeltaL(T;T_0,B) &\le  \left| 1 - \frac{\lambda_1}{B\gamma} - \UDeltaL \right| \nonumber\\ &+ \max_{T'\in[T_0,T]}A_1^{-1}(T') \displaystyle\sum_{k = 2}^\infty |\Omega| M^2(\psi_k) \left| 1 - \frac{\lambda_k}{B\gamma} - \UDeltaL(B) \right| \exp(-\lambda_k T'/(B\gamma)) \; . \label{eq:A3Delta}
\end{align}

Then, we first expand $A_1^{-1}$ in the second term in \eqref{eq:A3Delta} which yields
\begin{align}
&A_1^{-1}(T') \displaystyle\sum_{k = 2}^\infty |\Omega| M^2(\psi_k) \left| 1 - \frac{\lambda_k}{B\gamma} - \UDeltaL(B) \right| \exp(-\lambda_k T'/(B\gamma)) \nonumber\\
&\quad\quad\quad= \sum_{k = 2}^\infty \frac{\displaystyle|\Omega| M^2(\psi_k)}{\OmMeas M^2(\psi_1)} \left| 1 - \frac{\lambda_k}{B\gamma} - \UDeltaL(B) \right| \exp(-(\lambda_k-\lambda_1) T/(B\gamma))\,.\label{eq:twoninetynine}
\end{align}
We next consider the series terms, noting that, from the triangle inequality,
\begin{align}
\left|1-\frac{\lambda_k}{B\gamma}-\UDeltaL(B)\right| &\le \left|\frac{\lambda_k-\lambda_1}{B\gamma}\right| + \left|1-\frac{\lambda_1}{B\gamma}-\UDeltaL(B)\right|.
\end{align}
We note that the second term on the right-hand side is $\pcalO(B^2)$ (Lemma \ref{lem:EdeltaL_A}) and the first term on the right-hand side multiplied by $\exp(-(\lambda_k-\lambda_1) T'/(B\gamma))$ is of the form $g(\nu,z)$, where $\nu=(\lambda_k-\lambda_1)/(B\gamma)$ and $z=T'$.
Using this inequality in conjunction with Lemma \ref{lem:EdeltaL_B}, we obtain
\begin{align}
\max_{T'\in[T_0,T]}\left| 1 - \frac{\lambda_k}{B\gamma} - \UDeltaL(B) \right| e^{-(\lambda_k-\lambda_1) T'/(B\gamma)}&\le\frac{1}{e\,T_0}+\pcalO(B^2).
\end{align}
Now we take this factor outside the summation and recognize that, from Lemma \ref{lem:eigfcnamp},
\begin{align}
\sum_{k=2}^\infty \frac{|\Omega| (M(\psi_k))^2}{\OmMeas (M(\psi_1))^2} &= \frac{1 - \OmMeas (M(\psi_1))^2}{\OmMeas (M(\psi_1))^2}\nonumber\\
&= \frac{1}{\OmMeas M^2(\psi_1)} - 1\nonumber\\
&= \Upsilon B^2 + \pcalO(B^3).\label{eq:sumsimp}
\end{align}

Now returning to \eqref{eq:A3Delta} and incorporating components \eqref{eq:twoninetynine} -- \eqref{eq:sumsimp}, we find
\begin{align}
\EDeltaL(T;T_0,B) \le|\gamma\chi-\gamma^2\Upsilon-\phi^2|\frac{B^2}{\gamma^2} + \frac{B^2\Upsilon}{e T_0}+\pcalO(B^3)\,.
\end{align}
Finally, we divide by $\UDeltaL(B)\equiv \phi \frac{B}{\gamma}\left(1+\phi\frac{B}{\gamma}\right)^{-1}$ and factor out $B/\gamma$ to obtain
\begin{align}
\EDeltaLr(T;T_0,B) \le \biggl(\,\underbrace{\frac{|\gamma\,\chi-\gamma^2\,\Upsilon-\phi^2|}{\phi}}_{C_1} + \underbrace{  \frac{\gamma^2\,\Upsilon}{e\phi}}_{C_0} \frac{1}{T_0}\biggr)(B/\gamma) + \pcalO(B^2)  \text{ as } B \rightarrow 0 \, . \label{eq:A8Delta}
\end{align}
Evaluation at $T=\Tf$ concludes the proof.

%% file: appd.tex
\section{Solution Bounds for Non-Uniform Heat-Transfer Coefficient}\label{sec:htc_sketch}

In the main body of this paper we consider a heat transfer coefficient $\dv{h}$ which is uniform: independent of space and time. Here we consider a generalization more relevant in actual applications. As we shall see, the new results take good advantage of the framework and theory already developed in the main body of the paper for the case of uniform heat transfer coefficient.

We shall consider forced convection, in which case the heat transfer coefficient does not depend on the wall-ambient temperature difference. (In particular, we shall not address here natural convection or radiation, though both can be incorporated in the framework, in particular through linearization approaches.) However, in general, the (local) forced-convection heat transfer coefficient may depend on both space and time, in which case in equation (12) we replace $B$ with $\calB(x,t)$. In the case in which the heat transfer coefficient depends only on space, $\calB(x)$, then the first-order lumped approximation for $\uavg$ remains valid if we now interpret $B$ as the average of $\calB(\cdot)$ over $\partial\Omega$; however, the {\em error estimate} $\EavgL{1\,\text{asymp}}$ of \eqref{eq:eavgla}  --- premised on time-independent spatially uniform heat transfer coefficient --- is no longer valid. In the case in which the heat transfer coefficient depends on both space and time, even the first-order lumped approximation for $\uavg$ must be re-examined, and more generally our error analysis framework in the main body of the paper is no longer directly applicable.

\subsection{A Two-Sided Bound}\label{sec:twobound}

We first consider a perturbed dunking problem:

\begin{proposition}[Weak Statement: Perturbed Problem] \label{prop:ws_ppR} We are provided (i) a function $\calB\in L^\infty( \calT,L^\infty(\Omega))$, (ii) general initial condition $\tiu_0 \in L^2(\Omega)$, and (iii) ``perturbation'' linear functional $p \in
L^1(\calT,L^2(\Omega)) + L^2(\calT,X')$. Here, $X\equiv H^1(\Omega)$, $X'$ is the dual space to $X$, and $\calT=[0,t_\text{final}]$.

Then there exists a unique field $\tiu(\cdot,t;\calB) \in \LiLt\, \cap\, \LtX$ solution to
\begin{align}
m(\partial_t \tiu,v;\sigma) + a_0(\tiu,v;\kappa) + a_1^\star(\tiu,v;\calB)  & = p(v) \; , \; \forall v \in X,\; \forall t \in \calT  \; , \label{eq:ws_ppR}
\end{align}
subject to initial condition
\begin{align}
\tiu(x,t = 0;\calB) & = \tiu_0(x),\; \forall x \in \Omega \; .\label{eq:ws_ppR2}
\end{align}
Here, $a_1^\star(w,v;\calB) = \int_{\partial\Omega} \calB w v$.
\end{proposition}
We note that $a_1^\star(w,v;B) = B a_1(w,v)$ for $B \in \Rzp$ and we thus recover the original dunking initial value problem of Section \ref{App:ivp} when $p(v) = 0$. We can now proceed to the positivity statement for non-negative initial conditions and non-negative perturbations:
\begin{proposition}[Positivity of Solution to Perturbed Problem] \label{prop:MPGSR} Let $\tiu(\cdot,\cdot;\calB)$ satisfy the perturbed problem of Proposition \ref{prop:ws_ppR}; we shall further restrict, for the purposes of the current proposition, linear functional $p$ to the space $L^2(\calT,L^2(\Omega))$. We take $\tiu_0$ non-negative, and $p(v) \ge 0$ for all non-negative $v \in X$. Then $\tiu(x,t) \ge 0$, $\forall x \in \Omega, \forall t \in \calT$.
\end{proposition}

Proof: We start the proof by introducing $\wm(\cdot,t;\calB) \equiv -\min(\tiu(\cdot,t;\calB),0)$; note that (for any $x \in \Omega$) $\wm(x,t;\calB)  = 0$ for $\tiu(x,t;\calB) \ge 0$,  $\wm(x,t;\calB) = -\tiu(x,t;\calB)$ for $\tiu(x,t;\calB) \le 0$, and $\wm(x,t;\calB)\ge 0$. We first observe that $\| \wm(\cdot,0;\calB) \| = 0$, where $\| \cdot \| \equiv m^{1/2}(\cdot,\cdot)$. We next note that, for any $t \in \calT_+\equiv\calT\setminus\{0\}$, $\wm(\cdot,t;\calB) \in X \equiv H^1(\Omega)$; we may thus choose $\wm(\cdot,t;\calB)$ as test function $v$ in \eqref{eq:ws_ppR} to obtain
\begin{align}
m(\partial_t \tiu,\wm) + a^\star(\tiu,\wm;\calB)  = p(\wm) \;, \forall t \in \calT_+ \; \label{eq:mpR_1}
\end{align}
where
\begin{align}
a^\star(w,v;\calB) = a_0(w,v) + a_1^\star(w,v;\calB).
\end{align}
We now claim that
\begin{align}
m(\partial_t \tiu,\wm) & =  -m(\partial_t \wm,\wm) = -\frac{1}{2}\frac{d}{dt} \| \wm \|^2 \;, \forall t \in \calT_+\; ,  \label{eq:clb1R} \\
a^\star(\tiu,\wm;\calB) & = -a^\star(\wm,\wm;\calB) \;, \forall t \in \calT_+ \; . \label{eq:clb2R}
\end{align}
We then insert \eqref{eq:clb1R} -- \eqref{eq:clb2R} into \eqref{eq:mpR_1} to obtain
\begin{align}
\frac{1}{2} \frac{d}{dt}\| \wm \| ^2 + a^*(\wm,\wm,\calB)   = -p(\wm) \;, \forall t \in \calT_+ \; . \label{eq:mpR_2}
\end{align}
We now recall that $a^*$ is positive-semidefinite. We can further show that $-p(\wm) \le 0$: $p(v) \ge 0$ for all non-negative $v$, and hence $-p(v) \le 0$ for all non-negative $v$; by construction, $\wm$ is non-negative. We thus deduce from \eqref{eq:mpR_2} that $\| \wm(t) \|$ is non-increasing. It then follows from the initial condition $\| \wm(\cdot,0) \| = 0$ --- and since $\| \cdot \|$ is a norm ---  that $\| \wm(\cdot,t) \| = 0$  for all $t$ in $\calT$. But $\| \wm(\cdot,t) \| = 0$ implies $\| \min(\tiu(\cdot,t),0)\, \| = 0$ implies $\tiu(x,t) \ge 0, \forall x \in \Omega$, which concludes the proof. $\square$

Now we can make the following statement:
\begin{proposition}[$\calB$ Dependence: $\tiu_1{(\cdot,t;\calB)}$ is Non-Increasing in $\calB$] \label{prop:Bdep} Let $\tiu_1(\cdot,t;\sigma,\kappa,\calB^1)$ be the solution to \eqref{eq:ws_ppR} -- \eqref{eq:ws_ppR2} with parameters $\sigma$, $\kappa$, and $\calB^1$, and $\tiu_2(\cdot,t;\sigma,\kappa,\calB^2)$ be the solution to \eqref{eq:ws_ppR} -- \eqref{eq:ws_ppR2} with parameters $\sigma$, $\kappa$, and $\calB^2$. We now suppose $\calB^1 \le \calB^2$ for all $x \in \partial\Omega$ and $t \in \calT$. Then $\tiu_1(x,t) \ge \tiu_2(x,t)$ for all $x \in \Omega$ and $t \in \calT$.
\end{proposition}

Proof: We start with the substitution $\tiu_1(\cdot,t) = \tiu_2(\cdot,t) + w(\cdot,t)$. Then we have
\begin{align}
 m(\partial_t (\tiu_2+w),v) + a_0(\tiu_2+w,v) + a_1^\star(\tiu_2+w,v;\calB^1) &= 0\,.
\end{align}
By expanding the terms and invoking the equation for $u^*$, we have
\begin{align}
 m(\partial_t w,v) + a_0(w,v) + a_1^\star(w,v;\calB^1) &= a_1^\star(\tiu_2,v;\calB^2-\calB^1)\,.
\end{align}
We first apply Proposition \ref{prop:MPGSR} (with $p(v) = 0$) to the equation for $\tiu_2(\cdot,t)$ to deduce $\tiu_2\ge0$. We then apply Proposition \ref{prop:MPGSR} for $p(v) = a_1^\star(\tiu_2,v;\calB^2-\calB^1)$ to the equation for $w(\cdot,t)$: since $\tiu_2 > 0$ and $\calB^2-\calB^1 > 0$ then $p(v) > 0$ for all non-negative $v$. Thus $w(x,t) \ge 0$ for all $x \in \Omega$ and $t \in \calT,$ and hence $\tiu_1(x,t)\ge \tiu_2(x,t)$ for all $x\in\Omega$ and $t \in \calT$. This concludes the proof. $\square$

Now we can construct the upper and lower QoI bounds for the solution to a time-dependent non-uniform heat-transfer coefficient problem. Note in this appendix, QoI shall refer to the domain average QoI. In the remainder of this appendix, $p(v)$ of Proposition \ref{prop:ws_ppR} is zero.

\begin{definition}[Extreme Values of the Heat Transfer Coefficient]
We define 
\begin{align}
\calB_\text{inf} &\equiv \essinf_{(x,t)\in\Omega\times\calT} \calB(x,t),\\
\calB_\text{sup} &\equiv \esssup_{(x,t)\in\Omega\times\calT} \calB(x,t),\text{ and}\\
r&\equiv\frac{\Bsup-\Binf}{\Bsup}\label{eq:rdef},
\end{align}
where we recall that $\dashint_{\partial\Omega} \calB$ is the average of $\calB$ over $\partial\Omega$.
\end{definition}

\begin{proposition}[Domain-Average Bounds for Non-Uniform Heat Transfer Coefficient] \label{prop:htcbounds}
  Let $\calB \in L^\infty(\calT,L^\infty(\partial\Omega))$ be a non-uniform heat-transfer coefficient and fix $\sigma$ and $\kappa$. Then, the domain mean of the solution to the perturbed problem is bounded by the domain mean of the solutions to the perturbed problem for $B=\Binf$ and $B=\Bsup$, respectively $u(x,t;\Binf)$ and $u(x,t;\Bsup)$: 
\begin{align}
u_\avg(t;\Bsup) \le \tiu_\avg(t;\calB) \le u_\avg(t;\Binf),\label{eq:uhtc_bounds}
\end{align}
for all $t\in\calT$.  Here $\tiu_\avg=\dashint_{\Omega}\sigma \tiu$.
\end{proposition}

Proof: The solution bounds are obtained directly by application of Proposition \ref{prop:Bdep} to $u(\cdot,t;\calB)$, $u(\cdot,t;\Bsup)$, and $u(\cdot,t;\Binf)$:
\begin{align}
u(x,t;\Bsup) \le \tiu(x,t;\calB) \le u(x,t;\Binf),\ \forall\,(x,t)\in\Omega\times\calT.
\end{align}
To obtain the QoI bound, we apply $M$, defined in Equation \eqref{eq:Mdef}, to the left and right inequalities. Since the weight function $\sigma$ is strictly positive, the inequality holds for the integrals:
\begin{align}
M(u(\cdot,t;\Bsup)) \le M(\tiu(\cdot,t;\calB)) \le M(u(\cdot,t;\Binf)).\label{eq:uhtc_bounds2}
\end{align}
 This concludes the proof. $\square$

We now introduce bounding variables:
\begin{definition}[QoI Bounding Variables for Non-Uniform Heat Transfer Coefficient]\label{def:qoibound} We define
\begin{align}
\uavgUB(t;\calB)&\equiv \uavgL{1}(t;\Binf) + \EavgLa{1}(\Binf),\label{eq:uavgub}\\
\uavgLB(t;\calB)&\equiv \uavgL{1}(t;\Bsup),\\
\uavgMID(t;\calB)&\equiv \frac 1 2 (\uavgLB + \uavgUB),
\end{align}
where $\EavgLa{1}$ is given by \eqref{eq:eavgla}.
\end{definition}
 With these definitions we can now claim:
\begin{proposition}[Two-Sided Bound]\label{prop:twosidedbound} The solution QoI satisfies
\begin{align}
\uavgLB(t;\calB) \le \tiuavg(t;\calB)\le \uavgUB(t;\calB) + \pcalO((\Binf)^2)\text{ as }\Binf\rightarrow 0
\end{align}
for all $t\in\calT$.
\end{proposition}

Proof: We start with the $\tiuavg(t;\calB)$ bound of Proposition \ref{prop:htcbounds}. Applying Proposition \ref{prop:UavgBound} to the left inequality for $B=\Bsup$ and Proposition \ref{Prop:N1} to the right inequality for $B=\Binf$, we obtain the result. $\square$

We can then provide

\begin{proposition}[A Priori Estimate for Bound Gap] \label{prop:ubounds_eta} We assume that $r$ is bounded as $\Bsup\rightarrow 0$. For a non-uniform heat transfer coefficient $\calB$, $\uavgMID(t;\calB)$ satisfies
\begin{align}
\max_{t\in[0,\infty)}|\tiu_\avg(\cdot,t;\calB)  - \uavgMID(t;\calB)| &\le G + \pcalO((\Binf)^2)\text{ as }\Binf\rightarrow 0,
\label{eq:rbound}
\end{align}
where
\begin{align}
G\equiv\frac{1}{2}\Biggl((1-r)^{\frac{1-r}{r}} - (1-r)^{\frac{1}{r}}+ \EavgLa{1}(\Binf)\Biggl)
,\label{eq:gdef}
\end{align}
and $r$ is defined in Equation \eqref{eq:rdef}.  We may further apply the asymptotic expansion
\begin{align}
(1-r)^{\frac{1-r}{r}} - (1-r)^{\frac{1}{r}} = r/e + \pcalO(r^2)
\end{align}
to obtain $G = \frac{r}{2e} + \pcalO(r^2) + \pcalO(\Binf)$ as $r\rightarrow 0$, $\Binf\rightarrow 0$.
\end{proposition}

Proof: We start by recalling the inequality
\begin{align}
\uavgLB(t;\Bsup) \le \tiu_\avg(t;\calB) \le \uavgUB(t;\Binf)
\label{eq:ubar_bounds}
\end{align}
from Proposition \ref{prop:twosidedbound}. We now subtract $\uavgMID$ from each term in \eqref{eq:ubar_bounds} to obtain
\begin{align}
|\tiu_\avg(t;\calB) - \uavgMID(t;\calB)| \le \frac 1 2 (\uavgUB - \uavgLB),\ \forall\,t\in\calT.\label{eq:umidbounds}
\end{align}
We now consider the right-hand-side difference of \eqref{eq:umidbounds} under our hypothesis that the difference parameter $r$ is bounded as $\Bsup \rightarrow 0$. Taking advantage of the asymptotic bounds on the first-order lumped approximation (Proposition \ref{Prop:N1}), we obtain
\begin{align}
\uavgUB(t;\calB) - \uavgLB(t;\calB) &\le (\exp(-\Binf\gamma t) + \EavgLa{1}(\Binf) + \pcalO((\Binf)^2))\nonumber\\
&- \exp(\Bsup\gamma t)\label{eq:ublbdiff1}.
\end{align}
We now recognize $\Binf = (1-r)\Bsup$ and rearrange the terms to obtain the inequality
\begin{align}
\uavgUB(t;\calB) - \uavgLB(t;\calB) &\le \exp(-(1-r)\Bsup\gamma t) -\exp(-\Bsup\gamma t)\nonumber\\
&+ \EavgLa{1}(\Binf) + \pcalO((\Binf)^2)\label{eq:ublbdiff2}.
\end{align}
(We note that the lower bound, $\uavgLB$, does not contain the error term, $\EavgLa{1}(\Bsup)$, as $\uavgL{1}$ is a lower bound for $u_\avg$ (Proposition \ref{prop:UavgBound}).) We now invoke Lemma \ref{lem:N1a} with change of variables $z = \Bsup \gamma t$ and $r=\epsilon$ to bound the difference in the exponential terms,
\begin{align}
g(z;r)\equiv \exp(-(1-r)z) - \exp(-z)\left(=\exp(-(1-r)\Bsup\gamma t) - \exp(-\Bsup\gamma t)\right),
\end{align}
with
\begin{align}
g(z_{\sup}(r);r) =  (1-r)^{\frac{1-r}{r}} - (1-r)^{\frac 1 r},\label{eq:gsup}
\end{align}
where the $\sup$ is attained at $z_{\sup}(r)=-\frac{\ln(1-r)}{r}$.  Finally, inserting \eqref{eq:gsup} into \eqref{eq:ublbdiff2} we obtain
\begin{align}
\uavgUB - \uavgLB &\le g(z_{\sup}(r);r)+\EavgLa{1}(\Binf) + \pcalO((\Binf)^2)\\
&= (1-r)^{\frac{1-r}{r}} - (1-r)^{\frac{1}{r}}
+ \EavgLa{1}(\Binf) + \pcalO((\Binf)^2)\\
&= 2 G + \pcalO((\Binf)^2),
\end{align}
where $G$ is defined in Equation \eqref{eq:gdef}. Substitution into \eqref{eq:umidbounds} concludes the proof. $\square$

As an alternative to invoking Lemma \ref{lem:N1a}, the result of Equation \eqref{eq:ublbdiff1} can be obtained more directly through factorization and the application of the limit expression for the exponential function. We start with
\begin{align}
(1-r)^{\frac{1-r}{r}}-(1-r)^{\frac 1 r} &= (1-r)^{\frac{1}{r}}\left((1-r)^{-1}-1\right)\\
&= (1+(-1)r)^{\frac{1}{r}}\frac{r}{1-r}.
\end{align}
Now, noting that $\lim_{r\to 0}(1-r)^{1/r}\to e$, we recover the desired result.

In Table \ref{tab:rres}, we provide numerical results for $g(z_{\sup}(r);r)$.
\input{Tables/rres.tex}
\noindent Even for the larger $r$, the contribution of $g(z_{\sup}(r);r)$ to the bound remains relatively small. We emphasize that Proposition \ref{prop:ubounds_eta} serves only a theoretical purpose: to prove convergence as $r\rightarrow 0$. In actual practice, since $G$ is pessimistic, we will always evaluate our bounds for $\tiuavg(t;\calB)$ directly from $\uavgLB$ and $\uavgUB$.

\subsection{A One-Sided Bound}\label{sec:onebound}
In this section we shall assume that $\calB$ depends on space but {\em not} on time. We start by introducing new variables based on the spatial variation of $\calB$:
\begin{definition}[Spatial Dependence on $\calB$] We introduce the variables
\begin{align}
\overline{\calB}\equiv \dashint_\Omega \calB(x)\label{eq:bbardef},\\
\eta(\cdot)\equiv \calB(\cdot) / \overline{\calB};
\end{align}
note that $\dashint_{\partial\Omega}\eta = 1$.
\end{definition}

We may then claim
\begin{proposition}[A Lower Bound for {$\tiuavg(\cdot;\calB)$} for $\calB$ Time-Independent] \label{prop:eta_lb} The QoI $u_\avg(t;\calB)$ satisfies
\begin{align}
\tiu_\avg(t;\calB) \ge \uavgL{1}(t;\Bbar),
\end{align}
for all $t\in\calT$. Recall that $\tiuavg$ is defined in Proposition \ref{prop:htcbounds}.
\end{proposition}

Note here $\uavgL{1}$ is the first-order lumped approximation for $B = \Bbar$. The proof for this proposition follows from the same arguments developed in Lemma \ref{lem:UDeltaL} and Proposition \ref{prop:UavgBound}. The key changes are the definitions of $H$, $a_1$, and $L$, and $u_\paavg$, now
\begin{align}
H^\star(v;\eta)\equiv\dashint_{\partial\Omega} \eta\,v,\\
a_1^\star(w,v;\eta)\equiv \int_{\partial\Omega} \eta\, w \,v,\\
L^\star(v;\sigma,\eta)\equiv \gamma \int_\Omega \sigma \,v - \int_{\partial\Omega} \eta \,v,\\
u_\paavg^\star =\dashint_{\partial\Omega} \eta \,\tiu(=H^\star(\tiu;\eta)),\label{eq:upaavgstar}
\end{align}
and the substitution $\Bbar$ for $B$. Effecting these changes (and associated downstream changes), $\tiu_{\paavg}$ has the same relation to $\tiu_\avg$ as before: we retain the series representations of Proposition \ref{prop:sov_QoI}, and in particular, the constant $\lambda_k / (B\gamma)$ factor in \eqref{eq:sov_paavgQoI}. Furthermore, the modified definitions affect neither Proposition \ref{prop:l1ltBg}, since $\dashint_{\partial\Omega}\eta=1$, nor the eigenvalue ordering required in Lemma \ref{lem:UDeltaL}. Proposition \ref{prop:UavgBound} then directly applies, due to the new definition of $\eta$-dependent $\tiu_\paavg$ \eqref{eq:upaavgstar}. This concludes the proof. $\square$

\subsection{An Alternative Two-Sided Bound}\label{sec:newtwobound}

In this section, we continue to assume that $\calB$ is independent of time. We can now establish a sharper two-sided bound based on the lower bound developed in Section \ref{sec:onebound}. We first introduce
\begin{definition}[Bound Variables for $\calB$ Independent of Time] We define
\begin{align}
\uavgLBp(t;\calB)&\equiv\uavgL{1}(t;\Bbar),\label{eq:uavglbp}\\
\uavgMIDp(t;\calB)&\equiv \frac 1 2 (\uavgLBp(t;\calB) + \uavgUB(t;\calB)),\text{ and }\\
r'&\equiv \frac{\Bbar - \Binf}{\Bbar}.
\end{align}
\end{definition}
Armed with this improved lower bound (valid for time-independent $\calB$), we can re-purpose the framework of Proposition \ref{prop:twosidedbound} and Proposition \ref{prop:ubounds_eta} with $\uavgLB$, $\uavgMID$, and $r$ respectively replaced by $\uavgLBp$, $\uavgMIDp$, and $r'$. Hence
\begin{align}
\uavgLBp(t;\calB)\le\tiuavg(t;\calB)\le\uavgUB(t;\calB).\label{eq:newtwosided}
\end{align}
We further note that the bound gap
\begin{align}
G'\equiv\frac{1}{2}\Biggl((1-r')^{\frac{1-r'}{r'}} - (1-r')^{\frac{1}{r'}}+ \EavgLa{1}(\Binf)\Biggl)\label{eq:newtwosidedgap}
\end{align}
satisfies $G'\le G$. In short, time independence of $\calB$ permits us to improve our bounds. We shall refer to Equation \eqref{eq:newtwosided} as the \textit{alternative} two-sided bound, where $\uavgLBp$ and $\uavgUB$ are given in \eqref{eq:uavglbp}  and \eqref{eq:uavgub}, respectively.

\section{Connection to Real Flows}\label{sec:realflows}

We now relate our results of Appendix \ref{sec:htc_sketch} to a real fluid flow representative of forced convection more generally. We consider incompressible flow past a cylinder at $\mathrm{Re} \equiv \frac{V D}{\nu} = 150$ and $\mathrm{Pr} \equiv \frac{\nu}{\alpha_f} = 0.71$; here $V$ is the characteristic flow speed, $D(\equiv\ell, \text{ our extrinsic length scale})$ is the cylinder diameter, $\nu$ is the kinematic viscosity, and $\alpha_f$ is the thermal diffusivity of the working fluid. We consider the full conjugate heat transfer problem \cite{Perelmand} and hence require no assumption on heat transfer coefficient.

The fluid and solid properties are denoted by $\cdot_f$ and $\cdot_s$, respectively. We note in particular that the fluid-solid effusivity ratio, $\sqrt{\rho_f c_f k_f}/\sqrt{\rho_s c_s k_s}$ is $0.013$, small compared to unity; the latter is, in fact, a requirement --- typically satisfied for working fluid air --- for application of the engineering heat transfer coefficient approach. Note also that the fluid-solid diffusivity ratio, $\equiv \frac{\alpha_f}{\alpha_s}$, is $11/13$, hence order unity.\footnote{Our values for physical properties do not correspond to any particular fluid-solid material pairs; rather, these values are chosen as representative but with some accommodation to moderate the disparate fluid-solid time scales --- and thereby facilitate computation.}

We apply the spectral element code Nek5000 \cite{Nek5000} to solve the conjugate heat transfer problem for the velocity field and temperature field. We can then post-process the solution to obtain the truth heat transfer coefficient, $h^\text{th}(x,t)\equiv -\frac{k_f\partial_n T_f(x,t)}{T_f(x,t)-T_{\infty}}$ for $x\in\partial\Omega$, and associated truth Nusselt number:
$\mathrm{Nu}^\text{th}\equiv h^\text{th} \ell / k_f$;
note superscript th refers to the truth (conjugate heat transfer) result.
We could now replicate the conjugate heat transfer solution with a much simpler heat equation formulation solely in the solid --- but with the truth heat transfer coefficient reflected in the Robin boundary condition: $u_{\avg}^\text{th}(t)=\tiuavg(t;\calB)$. (See \cite{VandZandK} for an earlier comparison of the lumped approximation to the full conjugate heat transfer formulation.) However, our interest is to illustrate the theories of Sections \ref{sec:twobound} -- \ref{sec:newtwobound}.

\begin{figure}[H]
\centering
\includegraphics[width=0.8\textwidth]{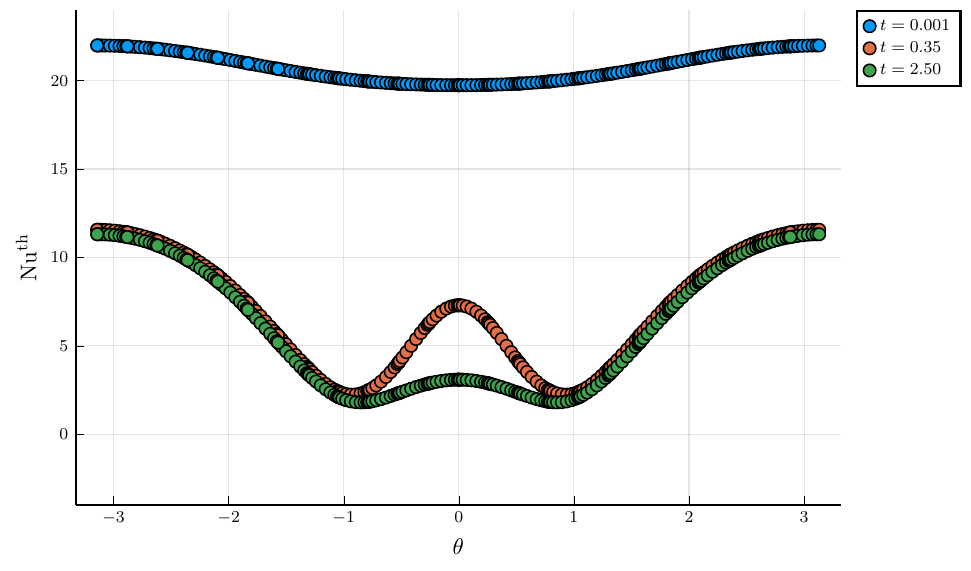}
\caption{Snapshots of $\mathrm{Nu}^\text{th}(\theta,t)$ at several times; here $\theta$ is the angle (in radians) such that $\theta = -\pi$ is the leading edge and $\theta=0$ is the trailing edge.}
\label{fig:Z1}
\end{figure}

\begin{figure}[H]
\centering
\includegraphics[width=0.8\textwidth]{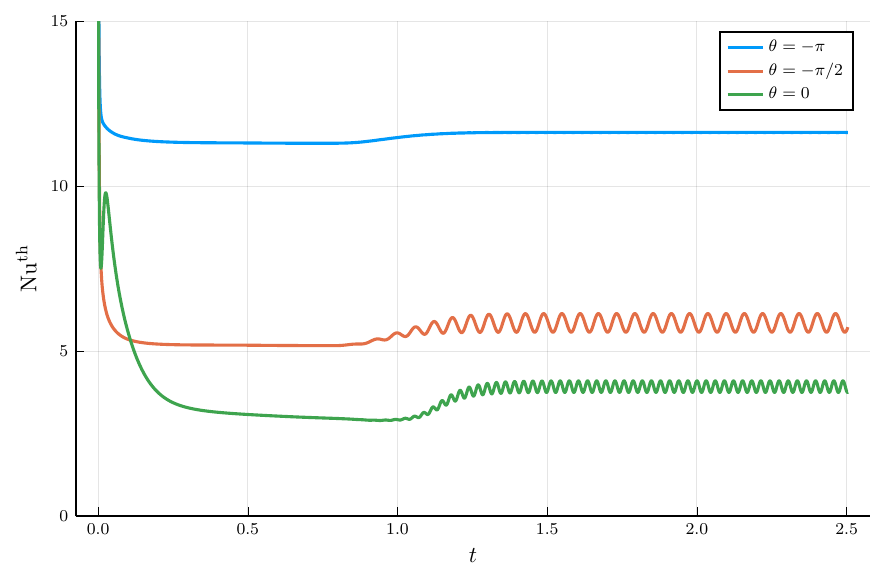}
\caption{Distribution of $\mathrm{Nu}^\text{th}(\theta,t)$ at several locations $\theta$; here $\theta$ is the angle (in radians) such that $\theta = -\pi$ is the leading edge and $\theta=0$ is the trailing edge.}
\label{fig:Z2}
\end{figure}

\begin{figure}[H]
\centering
\includegraphics[width=0.8\textwidth]{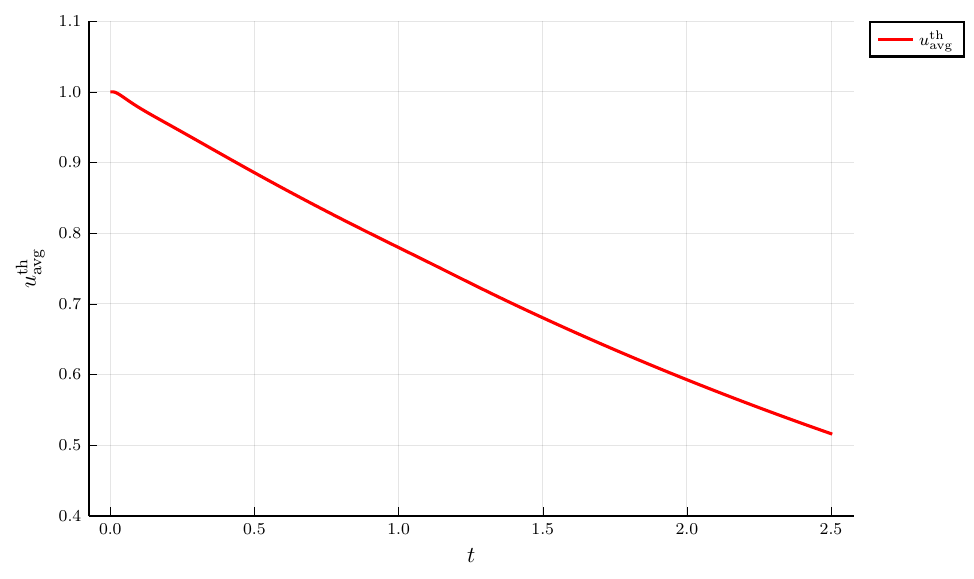}
\caption{Truth solution (CHT) for the domain average QoI, $u^\text{th}_\avg(t)$.}
\label{fig:Z3}
\end{figure}

To begin, we evaluate our quantity $\calB^\text{th}$ as $\mathrm{Nu}^\text{th}(k_f/k_s)$. For our fluid-solid pair, $(k_f/k_s) = 0.011$; it is typically the case that, for working fluid air and most common solids (in particular metals), $(k_f/k_s)$ is quite small, and hence ever for larger Nusselt number the Biot number can remain quite small. We plot in Figure \ref{fig:Z1} $\mathrm{ Nu }^\text{th}$ at several times as a function of space, and in Figure \ref{fig:Z2} $\mathrm{ Nu }^\text{th}$ at several spatial locations along the fluid-solid interface --- leading edge, side edge, and trailing edge --- as a function of time. We further plot, in Figure \ref{fig:Z3}, $\uavg^\text{th}$ as a function of time.  We note the presence of three time scales: a convective time scale, $\tau_{conv} \ll 1$ (in our nondimensionalization), a (solid) diffusive time scale, $\tau_{diff} = 1$ (by construction), and an equilibration time, $\tau_{eq} \gg 1$ --- which can be identified in Figure \ref{fig:Z2} and Figure \ref{fig:Z3}. (For a uniform heat transfer coefficient, $\tau_{eq} = 1/(B \gamma)$.)

In the first stage, as seen in Figure \ref{fig:Z2}, $\mathrm{Nu}^\text{th}$ is singular as $t \rightarrow 0$ . This singularity arises from the initial formation of a diffusive boundary layer --- subsequently replaced, for $t\gtrsim t_\text{conv}$, by a convective boundary layer. We observe from Figure \ref{fig:Z3} that the integrable singularity at time t = 0 has relatively little effect on the long-term temporal dependence of $\uavg^\text{th}$; $\uavg^\text{th}$ is approximately exponential in time, with only very small deviations for very small times. If we agree to neglect this short-time singularity, we can then redefine $\Binf$ and $\Bsup$ over the time interval $t > 0.25$ (say); applying the $\inf$ and $\sup$ operations yields $\Binf$ = 0.0185, $\Bsup$ = 0.126, and $G$ (of Proposition \ref{prop:ubounds_eta}) = 0.308.  We note $\Bsup$ (as re-defined) is quite small, and we thus expect the small-Biot approximation to be relevant.  

\begin{figure}[H]
\centering
\includegraphics[width=0.8\textwidth]{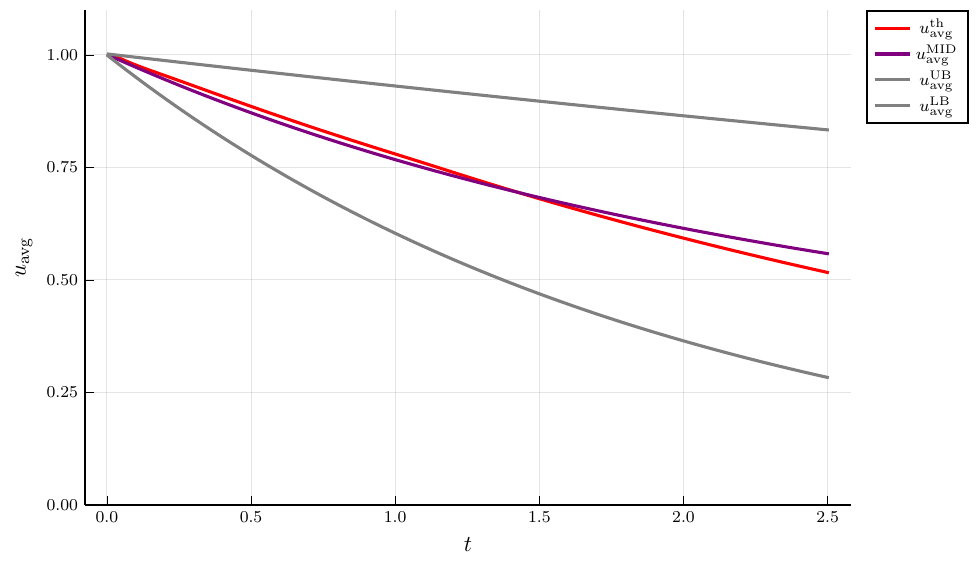}
\caption{Bounds provided by Proposition \ref{prop:twosidedbound}.}
\label{fig:Z4}
\end{figure}

We now turn to the two-sided bound from Proposition \ref{prop:ubounds_eta}, valid for $\calB$ which depends on space and time. We plot, in Figure \ref{fig:Z4}, $\uavgLB$, $\uavgUB$, and $\uavgMID(t)$, all provided in Definition \ref{def:qoibound}. We observe that, as might be expected from (our arguments related to the singularity at $t = 0$) and Proposition \ref{prop:ubounds_eta}, $\uavg^{\text{th}}(t)$ is indeed included within the bounding envelope. However, for this particular flow, the gap is not small, and hence the bounds are rather loose.

Proceeding to the one-sided bound, we again need some less rigorous arguments before we can apply our rigorous estimates. We note from Figure \ref{fig:Z2} that, following the fluid transient stage $0 \le t \le t_\text{diff} = 1$, the flow --- and $\mathrm{Nu}^{\text{th}}$ --- exhibits a stationary periodic oscillation in time, and that furthermore the period of oscillation is $\pcalO(t_\text{conv}) \ll t_\text{diff}$ (and of course then also $\ll t_\text{eq}$, where we recall that $t_\text{eq}$ is associated with the {\em thermal} transient). If we neglect the fluid transient period --- as not contributing substantially to the decay of $\uavg$, only partially supported by Figure \ref{fig:Z3} --- we can then apply homogenization in time.\footnote{Preliminary results have already been obtained in this direction. Our ongoing work on the topic will be reported in an upcoming separate publication \cite{CLB}.} (As $\Bsup\rightarrow 0$, the justification to neglect the fluid transient period is stronger.) In particular, we can then re-define our $\overline{\calB}$ of Equation \eqref{eq:bbardef} of Section \ref{sec:onebound} as
\begin{align}
(k_f/k_s) \dfrac{1}{|\partial\Omega|}  \int_{\partial\Omega}\left(\lim_{t \rightarrow \infty} \dfrac{1}{t}\int_0^t \mathrm{Nu}^\text{th}(x,t)\, dt\right) \, dx;
\end{align}
in our case, $\Bbar = 0.069$.

\begin{figure}[H]
\centering
\includegraphics[width=0.8\textwidth]{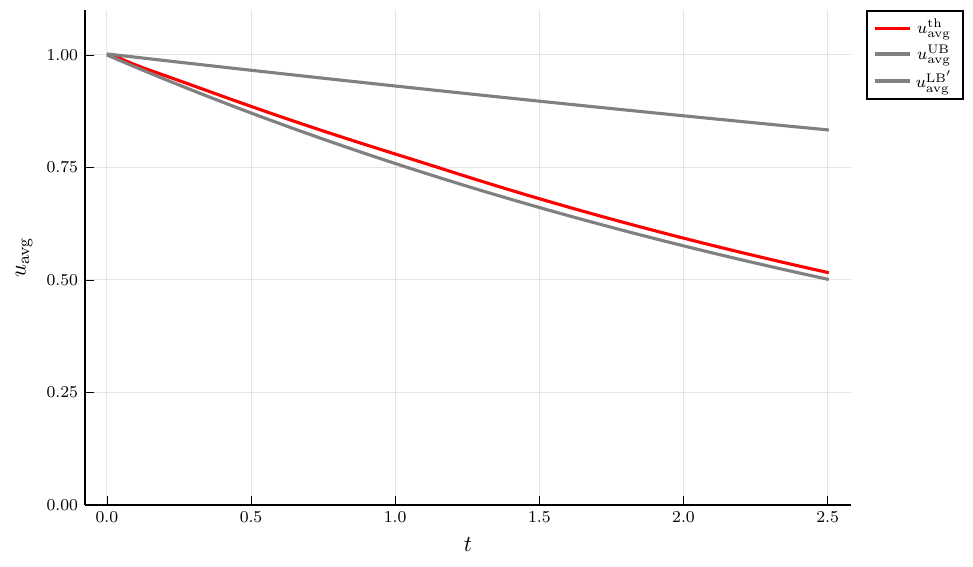}
\caption{Alternative two-sided bound.}
\label{fig:Z4p}
\end{figure}

We proceed directly to the alternative two-sided bound from Section \ref{sec:newtwobound} in this context. We plot, in Figure \ref{fig:Z4p}, $\uavgLBp$, and $\uavgUB$. We again observe that $\uavg^\text{th}(t)$ is indeed included within the improved bounding envelope.
In Figure \ref{fig:Z4p}, both $\uavg^\text{th}(t)$ and $u_\avg^\text{LB'}\equiv\uavgL{1}(t;\overline{\calB})$ are shown; we observe that our lower bound conjecture is confirmed.  We acknowledge that, in practice, experimental measurement of $\Binf$ and $\Bsup$ is difficult and hence rare. In some cases $\Binf$ and $\Bsup$ may be estimated. Also, for computational evaluation of the heat transfer coefficient --- now more common --- $\Binf$ and $\Bsup$ are indeed readily available. In any event, our theoretical framework, and associated realization in Figure \ref{fig:Z4} and Figure \ref{fig:Z4p}, does provide a sense of stability, and furthermore justifies (to a certain extent) our framework of the main body of the paper --- in which the heat transfer coefficient is treated as uniform.

%% file: Tables/rres.tex
\begin{table}[H]
\centering
\caption{Error contribution from the exponential difference. }
\label{tab:rres}
\begin{tabular}{c | c | c}
$r$ & $g(z_{\sup}(r);r)$ & $r\exp(-1)$\\
\hline
0.01   & 0.00370  & 0.00368 \\
0.02   & 0.00743  & 0.00736 \\
0.05   & 0.01887  & 0.01839 \\
0.10   & 0.03874  & 0.03679 \\
0.20   & 0.08192  & 0.07358 \\
0.50   & 0.25000  & 0.18394
\end{tabular}
\end{table}

%% file: main.bbl
\begin{thebibliography}{10}

\bibitem{LandL}
J.~H. {Lienhard IV} and J.~H. {Lienhard V}, {\em A Heat Transfer Textbook}.
\newblock Phlogiston Press, {V}ersion 5.10, 14 August 2020.

\bibitem{Cengel}
Y.~A. Cengel, {\em Heat and Mass Transfer (a Practical Approach)}.
\newblock McGraw-Hill, 2007.

\bibitem{IandD}
F.~P. Incropera and D.~P. DeWitt, {\em Fundamentals of Heat and Mass Transfer}.
\newblock New York City, New York: John Wiley \& Sons, Inc., 4th edition~ed.,
  1996.

\bibitem{ODD}
M.~Gockenbach and K.~Schmidtke, ``{N}ewton's law of heating and the heat
  equation,'' {\em Involve}, vol.~2, no.~4, pp.~419--437, 2009.

\bibitem{Joseph}
D.~D. Joseph, ``Parameter and domain dependence of eigenvalues of elliptic
  partial differential equations,'' {\em Archive for Rational Mechanics and
  Analysis}, vol.~24, pp.~325--351, 1967.

\bibitem{Kato}
T.~Kato, {\em A Short Introduction to Perturbation Theory for Linear
  Operators}.
\newblock Springer, 1982.

\bibitem{Yano}
M.~Yano, ``Private communication.''
\newblock Acknowledgments for adaptive FE software and discussions.

\bibitem{Arendt}
W.~Arendt, A.~F.~M. {ter Elst}, and J.~Gl\"{u}ck, ``Strict positivity for the
  principal eigenfunction of elliptic operators with various boundary
  conditions,'' {\em Advanced Nonlinear Studies}, vol.~20, no.~3, pp.~633--650,
  2020.

\bibitem{Evans}
L.~C. Evans, {\em Partial Differential Equations}, vol.~19 of {\em Graduate
  Studies in Mathematics}.
\newblock American Mathematical Society, 2nd~ed., 2010.

\bibitem{BandGM}
G.~A. Baker and P.~Graves-Morris, {\em Pad\'e Approximants}.
\newblock Encyclopedia of Mathematics and its Applications, Cambridge
  University Press, 2~ed., 1996.

\bibitem{Ostro2008}
A.~G. Ostrogorsky and B.~B. Mikic, ``Explicit solutions for boundary value
  problems in diffusion of heat and mass,'' {\em Journal of Crystal Growth},
  vol.~310, no.~11, pp.~2691--2696, 2008.

\bibitem{PandW}
L.~E. Payne and H.~F. Weinberger, ``An optimal {P}oincar\'e inequality for
  convex domains,'' {\em Archive for Rational Mechanics and Analysis}, vol.~5,
  pp.~286--292, 1960.

\bibitem{FaberKrahn}
D.~Daners, ``A {F}aber-{K}rahn inequality for {R}obin problems in any space
  dimension,'' {\em Mathematische Annalen}, vol.~335, pp.~767--785, 2006.

\bibitem{RockWets98}
R.~Rockafellar and R.~J.-B. Wets, {\em Variational Analysis}.
\newblock Heidelberg, Berlin, New York: Springer Verlag, 1998.

\bibitem{BandLR}
S.~C. Brenner and L.~{Ridgway Scott}, {\em The Mathematical Theory of Finite
  Element Methods}.
\newblock Springer Science+Business Media, LLC, third~ed., 2008.

\bibitem{SandB}
B.~Szab{\'o} and I.~Babu{\v{s}}ka, {\em Introduction to Finite Element
  Analysis: Formulation, Verification and Validation}.
\newblock John Wiley \& Sons, Ltd, first~ed., 2011.

\bibitem{GandN}
D.~S. Grebenkov and B.-T. Nguyen, ``Geometrical sructure of {L}aplacian
  eigenfunctions,'' {\em SIAM Review}, vol.~55, no.~4, pp.~601--667, 2013.

\bibitem{Perelmand}
T.~L. Perelman, ``{On Conjugated Problems of Heat Transfer},'' {\em
  International Journal of Heat and Mass Transfer}, vol.~3, pp.~293--303, 1961.

\bibitem{Nek5000}
P.~Fischer, J.~Kruse, J.~Mullen, H.~Tufo, J.~Lottes, and S.~Kerkemeier,
  ``Nek5000: Open source spectral element cfd solver,'' {\em Argonne National
  Laboratory, Mathematics and Computer Science Division, Argonne, IL, see
  https://nek5000. mcs. anl. gov/index. php/MainPage}, vol.~2, 2008.

\bibitem{VandZandK}
Z.~Virag, M.~Živić, and S.~Krizmanić, ``Cooling of a sphere by natural
  convection – the applicability of the lumped capacitance method,'' {\em
  International Journal of Heat and Mass Transfer}, vol.~54, no.~11,
  pp.~2303--2309, 2011.

\bibitem{CLB}
T.~Guo, K.~Kaneko, C.~Le~Bris, and A.~Patera, ``Homogenization in time applied
  to the {R}obin heat equation.''
\newblock In preparation.

\end{thebibliography}
